%% file: main.tex
\documentclass[final,hidelinks,onefignum,onetabnum]{siamart220329}

\input{ex_shared}

\ifpdf
\hypersetup{
  pdftitle={An Example Article},
  pdfauthor={D. Doe, P. T. Frank, and J. E. Smith}
}
\fi

\usepackage[margin=1in]{geometry}
\usepackage{algpseudocode}
\usepackage{mathfont}
\usepackage{booktabs}

\newtheorem{assumption}{Assumption}

\algrenewcommand\algorithmicrequire{\textbf{Input:}}
\algrenewcommand\algorithmicensure{\textbf{Output:}}

\newcommand{\utt}{u_{\theta_{t}}}
\newcommand{\ut}{u_{\theta}}

\newcommand{\barOmega}{\bar{\Omega}}
\newcommand{\barTheta}{\bar{\Theta}}
\newcommand{\barvarepsilon}{\bar{\varepsilon}}
\newcommand{\hatvarphi}{\hat{\varphi}}

\newcommand{\hatthetat}{\hat{\theta}_{t}}
\newcommand{\tildethetat}{\tilde{\theta}_{t}}

\definecolor{brickred}{rgb}{0.8, 0.25, 0.33}

\date{}

\begin{document}
\maketitle

\begin{abstract}
    We develop a novel computational framework to approximate solution operators of evolution partial differential equations (PDEs). By employing a general nonlinear reduced-order model, such as a deep neural network, to approximate the solution of a given PDE, we realize that the evolution of the model parameters is a control problem in the parameter space. Based on this observation, we propose to approximate the solution operator of the PDE by learning the control vector field in the parameter space. From any initial value, this control field can steer the parameter to generate a trajectory such that the corresponding reduced-order model solves the PDE. This allows for substantially reduced computational cost to solve the evolution PDE with arbitrary initial conditions. We also develop comprehensive error analysis for the proposed method when solving a large class of semilinear parabolic PDEs. Numerical experiments on different high-dimensional evolution PDEs with various initial conditions demonstrate the promising results of the proposed method. 
\end{abstract}

\section{Introduction}
\label{sec:intro}
Partial differential equations (PDEs) are ubiquitous in modeling and are vital in numerous applications from finance, engineering, and science \cite{evans1998partial}. As the solutions of many PDEs lack analytical form, it is necessary to use numerical methods to approximate the solutions \cite{atkinson1989introduction, evans1998partial}. Traditional numerical methods such as finite difference and finite element methods rely upon the discretization of problem domains, which does not scale to high-dimensional problems due to the so-called ``curse of dimensionality''.

In recent years, deep neural networks (DNNs), which can be thought of as a type of nonlinear reduced order models, have emerged as powerful tools for solving high-dimensional PDEs \cite{raissi2019physics-informed, bao2020numerical,han2017overcoming,e2018deep,cuomo2022scientific,han2018high-pde,han2018deep,kumar2011multilayer}. For example, in \cite{raissi2019physics-informed, bao2020numerical,e2018deep,cuomo2022scientific,zang2020weak}, the solution of a given PDE is parameterized as a DNN, and the network parameters are trained to minimize potential violations (in various definitions) to the PDE. 
These methods have shown numerous successes in solving a large variety of PDEs empirically. Their successes are partly due to the provable universal approximation power of DNNs  \cite{hornik1991approximation, yarotsky2017error, liang2017deep}. 
On the other hand, these methods aim at solving specific instances of PDEs, and as a consequence, they need to start from scratch for the same PDE whenever the initial and/or boundary value changes. 

%One drawback of many state-of-the-art deep learning PDE methods is they can be expensive, and potentially time-consuming, to train. This leads to issues when a user wants to retrain a PDE solution to account for varying initial conditions. For example, in many applications, the underlying dynamics of a PDE problem are known and fixed (that is the differential operator of the PDE is well-defined), if a user wanted to vary the initial conditions in order to see the evolution of different solutions, many current neural network methods would require an expensive retrain to account for the changed initial condition \cite{raissi2019physics-informed,e2018deep,bao2020numericalarxiv,zang2020weak}. This imposes a large training cost for every tested initial condition. 

There have also been recent studies to find solution operators of PDEs \cite{li2020fourier,lu2019deeponet}. These methods aim at finding the map from the problem's parameters to the corresponding solution. Finding solution operators has substantial applications as the same PDE may need to run many times with different initial or boundary value configurations. However, existing methods fall short in tackling high-dimensional problems as many require spatial discretization to represent the solution operators using DNNs. 

In this paper, we propose a new approach to find solution operators of high-dimensional evolution PDEs. For a given PDE, we first parameterize its solution as a general reduced-order model, such as a DNN, whose parameters denoted as $\theta$ are to be determined. Then we seek to find a vector field on the parameter space which describes how $\theta$ evolves in time. This vector field essentially acts as a controller on the parameter space, steering the parameters so that the induced DNN evolves and approximates the PDE solution for all time. Once such a vector field is found, we can easily change the initial conditions of the PDE by simply starting at a new point in the parameter space. Then we follow the control vector field to find the parameters trajectory which gives an approximation of the time-evolving solution. Thus, different initial conditions can be considered for the same PDE without solving it repeatedly. 
Our contributions can be summarized as follows.
\begin{enumerate}
    \item We develop a new computational framework to find the solution operator of any given initial value problem (IVP) defined by high-dimensional nonlinear evolution PDEs. This framework is purely based on the evolution PDE itself and does not require any solutions of the PDE for training. Once we find the solution operator, we can quickly compute solutions of the PDE with any initial value at a low computational cost.
    
    \item We provide comprehensive theoretical analysis to establish error bounds for the proposed method when solving linear PDEs and some special nonlinear PDEs.
    
    \item We conduct a series of numerical experiments to demonstrate the effectiveness of the proposed method in solving a variety of linear and nonlinear PDEs.
\end{enumerate}

% \subsection{Paper organization}
The remainder of this paper is organized as follows. 
In Section \ref{sec:related}, we provide an overview of recent neural network based numerical methods for solving PDEs. We outline the fundamentals of our proposed approach in Section \ref{subsec:model} and provide details of our method and its key characteristics in Section \ref{subsec:method}. We conduct comprehensive error analysis in Section \ref{subsec:error-analysis}. We demonstrate the performance of the proposed method on several linear and nonlinear evolution PDEs in Section \ref{sec:numerical-results}. Some variations and generalizations of the proposed approach are given in Section \ref{sec:variation}. Finally, Section \ref{sec:conclusion} concludes this paper.

\section{Related Work}
\label{sec:related}

\subsection{Classical methods for solving PDEs}
Classical numerical methods for solving PDEs, such as finite difference \cite{thomas2013fdm} and finite element methods \cite{johnson2012numerical}, discretize the spatial domain using mesh or triangulation. These methods convert a PDE to its discrete counterpart, which is a system of algebraic equations with finite number of unknowns, and solve the system to obtain approximate solution on the grid points
\cite{ames2014numerical,thomas2013numerical,evans2012numerical,quarteroni2008numerical}.
These methods have been significantly advanced in the past decades, and they are able to handle complicated situations such as irregular domains. However, they severely suffer ``curse of dimensionality'' when applied to high-dimensional problems---the number of unknowns increases exponentially fast with respect to spatial dimension, which renders them computationally intractable for many problems.

\subsection{Neural network based methods for solving PDEs}
%
% In recent years, DNNs \cite{goodfellow2016deep}, as a class of powerful nonlinear reduced-order models, have demonstrated unprecedented successes in a variety of scientific computing fields. 
% 
Early attempts using neural networks to solve PDEs can be seen in \cite{dissanayake1994neural-network-based,lagaris1998artificial,lee1990neural,kumar2011multilayer}.
DNNs emerged in recent years and demonstrated striking power in solving PDEs through various approaches 
%\cite{raissi2019physics-informed,e2018deep,bao2020numerical,zang2020weak,sirignano2018dgm:,ramabathiran2021spinn,nusken2021solving,yang2020physics,wang2022respecting,liang2022finite,wang2022isl2,zhang2020physics,berg2018unified}
\cite{raissi2019physics-informed,e2018deep,bao2020numerical,zang2020weak,sirignano2018dgm:,nusken2021solving,yang2020physics,berg2018unified}.
%
%For example, neural networks have been used to improve the efficiency and accuracy of finite difference methods \cite{lee1990neural,gobovic1994analog, yentis1996vlsi} and Galerkin method \cite{rudd2015constrained}, to solve ordinary differential equations (ODEs) \cite{meade-jr1994numerical,meade-jr1994solution,lagaris1998artificial}, matrix algebra problems \cite{wang1990structured}, and nonlinear differential equations \cite{tompson2017accelerating, suzuki2017neural}.
%
DNNs, which are the key machinery of deep learning, have demonstrated extraordinary potential in solving many high-dimensional nonlinear PDEs, which were considered computationally intractable using classical methods. For example, a variety of DNN based methods have been proposed based on strong form \cite{raissi2019physics-informed,nabian2018deep,dissanayake1994neural-network-based,berg2018unified,magill2018neural,pang2019fpinns,kharazmi2020hp,pang2020npinns,ramabathiran2021spinn}, variational form \cite{e2018deep}, and weak form \cite{zang2020weak,bao2020numerical} of PDEs.
They are considered with adaptive collocation strategy \cite{anitescu2019artificial}, adversarial inference procedure \cite{yang2019adversarial}, oscillatory solutions \cite{cai2020phase}, and multiscale methods \cite{liu2020multi,wang2020multi,cai2019multi}.
Improvements of these methods with adaptive activation functions \cite{jagtap2020adaptive}, networks structures \cite{gu2020selectnet,gu2020structure,huang2020int}, boundary conditions \cite{lyu2020enforcing,dong2020method}, structure probing \cite{huang2020int}, as well as their convergence \cite{luo2020two,shin2020convergence}, are also studied. Readers interested in these methods can also refer to \cite{ramabathiran2021spinn,yang2020physics,wang2022respecting,liang2022finite,wang2022isl2,zhang2020physics}.
Further, there are methods that can solve inverse problems such as parameter identifications 

For a class of high-dimensional PDEs which have equivalent backward stochastic differential equation (SDE) formulations due to Feynman-Kac theory, deep learning methods have been applied by leveraging such correspondences \cite{beck2017machine,fujii2017asymptotic,han2017overcoming,e2017deep,han2018high-pde, han2020solving,pham2021neural,hure2020deep,hutzenthaler2020proof}.
These methods are shown to be good even in high dimensions \cite{han2018high-pde,hure2020deep,pham2021neural}, however, they are limited to solving the special type of evolution equations whose generator function has a corresponding SDE.
% 
% Many of the algorithms then require doing multiple runs to get a full picture, using a large amount of computation. Therefore, if initial conditions change, a large amount of retraining needs to be done.

%\ye{There is a class of methods that solve PDE using backward SDEs (related to Feymann Kac formula). An early work is at \url{https://arxiv.org/abs/1707.02568} and there are many followups too. They need to use the relation between BSDE and PDE. It's good at even higher dimension but can only solve $u(x)$ at a given $x$ at a time.}

For evolution PDEs, parameter evolution algorithms \cite{du2021evolutional, bruna2022neural-galerkin,anderson2022evolution} have also been considered. These methods parameterize the PDE solution as neural network \cite{du2021evolutional,bruna2022neural-galerkin} or an adaptively chosen ansatz as discussed in \cite{anderson2022evolution}. In these methods, the parameters are evolved forward in time through a time marching scheme, where at each step a linear system \cite{bruna2022neural-galerkin,du2021evolutional} or a constrained optimization problem \cite{anderson2022evolution} needs to be solved.

\subsection{Learning solution operator of PDEs}
The aforementioned methods aim at solving specific instance of a given PDE, and they need to be rerun from scratch when any of the problem configuration (e.g., initial value, boundary value, problem domain) changes. In contrast, the solution operator of a PDE directly maps a problem configuration to its corresponding solution. To this end, several methods have been proposed to approximate Green's functions for some linear PDEs \cite{boulle2022learning-green,teng2022green,nicolas2022data-driven,lin2022bi-green}, as solutions to such PDEs have explicit expression based on their Green's functions.
However, this approach only applies to a small class of linear PDEs whose solution can be represented using Green's functions.
Moreover, Green's functions have singularities and it requires special care to approximate them using neural networks.
For example, rational functions are used as activation functions of DNNs to address singularities in \cite{nicolas2022data-driven}. In \cite{boulle2022learning-green}, the singularities are represented with the help of fundamental solutions. 

For general nonlinear PDEs, DNNs have been used for operator approximation and meta-learning for PDEs
%\cite{chen2020meta,fan2019bcr,li2020neural,cai2020deepm,clark2020deep,mao2020deepm,pang2020npinns,lu2019deeponet,lu2019deepxde,li2020fourier,kovachki2021universal,wen2022u-fno,wang2021learning}.
\cite{mao2020deepm,guo2018data-driven,lu2019deeponet,lu2019deepxde,li2020fourier,wen2022u-fno,wang2021learning,regazzoni2019machine}.
For example, the work \cite{guo2018data-driven} considers solving parametric PDEs in low-dimension ($d\le 3$ for the examples in the paper). Their method requires discretization of the PDE system and needs to be supplied by many full-order solutions for different combinations of time discretization points and parameter selections for their network training. Then their method applies proper orthogonal decomposition to these solutions to obtain a set of reduced basis to construct solutions for new problems.
The work \cite{regazzoni2019machine} requires a massive amount of pairs of ODE/PDE control and the corresponding system outputs, which are produced by solving the original ODE/PDE system; then the DNN is trained on such pairs to learn the mapping between these two subjects which are discretized as vectors by evaluating the functions only at grid points in the domain.
DeepONets \cite{lu2019deeponet,lu2019deepxde,wang2021learning} seek to approximate solution mappings by use of a ``branch" and ``trunk" network. FNOs \cite{li2020fourier,wen2022u-fno} use Fourier transforms to map a neural network to a low dimensional space and then back to the solution. 
In addition, several works that apply spatial discretization of the problem or transform domains and use convolutional neural networks (CNNs) \cite{raonic2023convolutional,guo2016cnn-operator,zhu2018bayesian} or graph neural networks (GNNs) \cite{kovachki2023neural,li2020neural,lotzsch2022learning}. 
Interested readers may also refer to generalizations and extensions of these methods in \cite{chen2020meta,fan2019bcr,li2020neural,cai2020deepm,clark2020deep,mao2020deepm,pang2020npinns,lu2019deepxde,kovachki2021universal}.
%by taking generic (GNN) or uniform (CNN) meshes of the relevant input domain or data and then passing this mesh as an input to a neural operator that seeks to map to the solution on a, potentially different, mesh. 
%One advantage  methods seek to leverage the local properties of many PDEs by having the input mesh form a graph or as the nodes of an image and then applying convolutional layers and other structures used in the image classification and generation communities.
A key similarity of all these methods is they require certain domain discretization and often a large number of labeled pairs of IVP initial conditions (or PDE parameters) and the corresponding solution obtained through other methods for training. This limits their applicability to high dimensional problems where such training data is unavailable or the mesh is prohibitive to generate due to cures of dimensionality. 
%s above 4 as either the training data or the meshes needed as inputs become prohibitive to generate. This is a particular drawback we seek to overcome in our approach to be outlined later. Nonetheless, The curious 

\subsection{Differences between our proposed approach and existing ones}
Different from all existing approaches, we propose to approximate solution operators of evolution PDEs in a control framework in parameter spaces induced by general reduced-order models such as DNNs. Unlike the existing solution operator approximation methods (e.g., DeepONet \cite{lu2019deeponet} and FNO \cite{li2020fourier}) which seek to directly approximate the infinite-dimensional operator, our approach is based on the relation between evolving solutions and their projected trajectories in the parameter space. This leads us to convert the problem of finding a solution operator over infinite-dimensional function space into a control vector field optimization problem over a finite-dimensional parameter space. As a result, the problem of solving an evolution PDE in continuous space is reduced to numerically solving a system of ODEs, which can be done accurately with very low computation complexity. Moreover, our approach does not require spatial discretization in any problem or transformed domain nor needs any basis function representation throughout problem formulation and computation. 
We provide mathematical insights into the parameter submanifold and its tangent spaces and establish their connection to the finite-dimensional parameter space. These
new insights led us to the proposed approach which approximates solution operators of PDEs by controlling network parameters in the parameter space.
These new features also enable our approach to solve evolution PDEs in high-dimensional cases. This is a significant advantage over existing operator learning methods such as DeepONet or FNOs as their spatial discretization schemes, which are used to generate the training data, hinder their application to high-dimensional cases. 
%\blue{Our approach can be thought of as a nonlinear autoencoder which reduces the dimension of an infinite-dimensional space to a finite-dimensional latent space. However, unlike traditional autoencoders which are often a blackbox, we shall provide mathematical insights about the finite-dimensional latent space, including defining and giving meaning to tangent spaces and trajectories on this space.  }
% \gaby{Seems we need more discussion here to satisfy reviewer 2, not sure what more needs to be added? Thoughts on this?}

\section{Proposed Method}
\label{sec:proposed}

The main goal of this paper is to develop a new computational framework to approximate the \emph{solution operator} for IVPs of high-dimensional evolution PDEs. The solution operator is a procedure that, once known, can efficiently map an arbitrarily given initial value $g$ to the solution of the IVP without solving the PDE again. 
We first propose to parameterize $u$ as a \emph{nonlinear reduced-order model}, such as a DNN, which is denoted by $\ut$ with parameters $\theta$, i.e., $\ut$ is a parametric function determined by the value of its finite-dimensional parameters $\theta$, and $\ut$ is used to approximate $u$.
%
% Then these nonlinear reduced-order models are parametric functions with controlled time-varying parameters to approximate the solutions of the evolution PDEs.

To find the solution operator, we propose to build a control vector field $V$ in the parameter space $\Theta$ where $\theta$ resides.
Then the solution operator can be implemented as a fast numerical solver of the ODE defined by $V$.
More precisely, we first find the parameters $\theta_{0}$ such that $u_{\theta_0}$ approximates $g$, then we follow the control vector field $V$ to obtain a trajectory $\{\theta_t\, | \, 0 \le t \le T\}$ in $\Theta$ with very low computational cost, which automatically induces a trajectory $\utt$ to approximate the true solution $u$ of the IVP with the initial value $g$.
We provide details of these constructions in the following subsections.

\subsection{Nonlinear reduced-order models and parameter submanifold} 
\label{subsec:model}
DNNs, which can be viewed as nonlinear reduced-order models, have emerged as powerful tools to solve high-dimensional PDEs in recent years \cite{raissi2017machine,raissi2017physics,raissi2019physics-informed,han2017overcoming,e2018deep,zang2020weak,bao2020numerical}.
Mathematically, a DNN can be expressed as the composition of a series of simple linear and nonlinear functions.
In the deep learning context, a typical building block of DNNs is called a {layer}, which is a mapping $h: \Rbb^{d} \to \Rbb^{d'}$ for some compatible input dimension $d$ and output dimension $d'$:
\begin{equation}
\label{eq:h}
    h(z; W, b) := \sigma(Wz + b),
\end{equation}
where $z \in \Rbb^d$ is the input variable of $h$, the matrix $W\in \Rbb^{d'\times d}$ and vector $b \in \Rbb^{d'}$ are called the weight and bias respectively, and $\sigma: \mathbb{R} \to \mathbb{R}$ is a nonlinear function that operates componentwise on its $d'$-dimensional argument vector $Wz + b$ (hence $\sigma$ is effectively a mapping from $\Rbb^{d'}$ to $\Rbb^{d'}$).
Common choices of activation functions include the hyperbolic tangent (tanh) and rectified linear unit (ReLU)  $\sigma(z)=\max(0,z)$. 
%Common choices of activation functions include the sigmoid $\sigma(z) = (1+e^{-z})^{-1}$, hyperbolic tangent (tanh)  $\sigma(z)=(e^{z} - e^{-z})/(e^{z} + e^{-z})$, rectified linear unit (ReLU)  $\sigma(z)=\max(0,z)$, rectified power unit (RePU) $\sigma(z)=\max(0,z)^p$ for $p>0$, swish function $\sigma(z) = z(1+e^{-z})^{-1}$, rational functions, and many others.
%
% Some activation functions can also have learnable parameters.
%
We only consider smooth activation functions $\sigma$ hereafter. % like those just listed.
A commonly used DNN structure $u_{\theta}$, often called feed-forward network (FFN), is defined as the composition of multiple layer functions of form \eqref{eq:h} as follows:
\begin{align}
& u_{\theta}(x) := u(x;\theta) = w^{\top}z_{L} + b, \label{eq:u_param} \\
\mbox{where} \quad & z_0 = x, \quad z_{l} = h_{l}(z_{l-1}):=h(z_{l-1};W_l, b_l),\quad l=1,\dots,L, \nonumber
\end{align}
and the $l$th hidden layer $h(\cdot;W_l,b_l): \Rbb^{d_{l-1}} \to \Rbb^{d_l}$ is determined by its weight and bias parameters $W_l\in \mathbb{R}^{d_{l} \times d_{l-1}}$ and $b_l \in \mathbb{R}^{d_{l}}$ for $l=1,\dots,L$ and $d_0 = d$. 
Here the output of $u_{\theta}$ is set to the affine transform of the last hidden layer $z_{NN} = h_L(z_{L-1})$ using weight $w\in\Rbb^{d_L}$ and bias $b \in \Rbb$. 
The \emph{network parameters} $\theta$ 
%\ye{So we use parameter or parameters? Need to be consistent.}\gaby{Parameters. I just missed this one so good catch. Thought I had caught them all.} 
refers to the collection of all learnable parameters (stacked as a vector in $\Rbb^{m}$) of $u_{\theta}$, i.e.,
\begin{equation}
\label{eq:theta}
\theta := (w, b, W_{L}, b_{L}, \dots, W_1, b_1) \in \mathbb{R}^{m},
\end{equation}
and training the network $u_\theta$ refers to finding the minimizer $\theta$ of some properly designed loss function.
%
% In practice, the choice of network depth $L$, layer widths $d_1,\dots,d_L$, and activation functions $\sigma$ (they can also vary across different hidden layers) are rather flexible. 
% %
% Indeed, a large variety of network architectures have been proposed to work more effectively for different applications.
% %
% The FCN structure \eqref{eq:u_param} is merely a basic DNN architecture, but it will help us to illustrate our main ideas and demonstrate the technical challenges in using DNNs as reduced-order models to solve PDEs in our discussion below. 

\begin{remark}
\label{remark:guhring}
DNNs are shown to be very powerful in approximating high-dimensional functions in a vast amount of studies in recent years, 
%
%The approximation ability of DNNs have been studied extensively, 
%
see, e.g., \cite{hornik1991approximation,hornik1989multilayer,liang2017deep,petersen2018optimal,guhring2020error,guhring2020approximation,yarotsky2017error,li2020repu}.
For example, it is shown in \cite{guhring2020error} that for any $M,\varepsilon > 0$, $k \in \mathbb{N}$, $p \in [1,\infty]$, and $\Omega=(0,1)^{d}\subset \Rbb^{d}$, denote $\Fcal:= \{f \in W^{k,p}({\Omega};\Rbb)\,|\, \|f\|_{W^{k,p}(\Omega)} \le M\}$, then there exists a DNN structure $\ut$ of form \eqref{eq:u_param} with sufficiently large $m$ and $L$ (which depend on $M$, $\varepsilon$, $d$ and $p$ only), such that for any $f\in \Fcal$, there is $\|u_{\theta}-f\|_{W^{k,p}(\Omega)} \le \varepsilon$ for some $\theta \in \Rbb^{m}$.
%
% Such approximation guarantees have been established for sets of $C^k$ functions as well as Sobolev spaces.
%
%This result implies that $\{\ut\,|\, \theta \in \mathbb{R}^{m}\}$ is an $\varepsilon$-net of the $M$-ball in the Sobolev space $W^{k,p}(\Omega)$, suggesting that DNNs are suitable to approximate solutions of PDEs.
%
This result suggests that DNNs are suitable to approximate solutoins of PDEs.
We note that this is one of the many error bounds of DNN approximations established in recent years, and such bounds are still being continuously improved nowadays.
\end{remark}

% \blue{
% \begin{remark}
% \label{remark:hornik}
% DNNs are shown to be very powerful in approximating high-dimensional functions in a vast amount of studies in recent years, 
% %
% %The approximation ability of DNNs have been studied extensively, 
% %
% see, e.g., \cite{hornik1991approximation,hornik1989multilayer,liang2017deep,petersen2018optimal}.
% %
% For example, it is shown in \cite{hornik1991approximation} that for any $M,\varepsilon > 0$, $k \in \mathbb{N}$, and any bounded open $\Omega \subset \Rbb^d$, denote $\mathcal{F}:= \{f \in C^{k}(\Omega)\cap C(\bar{\Omega})\,:\, \|f\|_{W^{k,\infty}(\Omega)} \le M\}$, then there exists a DNN architecture $\ut$ of form \eqref{eq:u_param} with sufficient large neuron number $m$ and depth $L$, which depend on $M$, $\varepsilon$, $\Omega$, and $d$ only, such that for any $f\in \mathcal{F}$, there is $\|u_{\theta}-f\|_{C^{k}(\Omega)} \le \varepsilon$ for some $\theta \in \Rbb^{m}$.
% %
% % Such approximation guarantees have been established for sets of $C^k$ functions as well as Sobolev spaces.
% %
% This result implies that $\{\ut\,:\, \theta \in \mathbb{R}^{m}\}$ is an $\varepsilon$-net of the $M$-ball of the Sobolev space $W^{k,p}(\Omega)$, suggesting that DNNs are suitable to approximate solutions of PDEs.
% %
% This is one of the many error bounds established in recent years, and such bounds are still being continuously improved nowadays. We refer the reader to \cite{guhring2020error,guhring2020approximation,yarotsky2017error,li2020repu} for some of these advances.
% \end{remark}
% }

Our approach relies on the key relation between the parameters $\theta$ and the reduced-order model $\ut$. 
More specifically, we identify the finite-dimensional parameter space $\Theta \subset \Rbb^{m}$ where $\theta$ belongs to and the submanifold $\Mcal$ of functions defined by 

\begin{equation}
    \label{eq:M}
    \Mcal := \set{u_{\theta}: {\Omega} \to \Rbb \ \vert \ \theta \in \Theta}.
\end{equation}
%
% First of all, we realize that $u_{\theta}$ \eqref{eq:u_param} renders an approximation to a solution using a finte-dimensional vector $\theta$. 
%
As we can see, $u_{\theta}$ defines a mapping from the parameter space $\Theta$ to the submanifold $\Mcal$ of the infinite-dimensional function space.
We call $\Mcal$ the \emph{parameter submanifold} determined by $u_\theta$.

To approximate a time-evolving function $u^{*}(\cdot,t)$, e.g., the solution of an evolution PDE, over time horizon $[0,T]$ using the reduced-order model $u_{\theta}$, we need to find a trajectory $\{\theta_t \in \Theta\,|\, 0 \le t \le T\}$ in the parameter space $\Theta$ so that $u_{\theta_t}(\cdot)$ is close to $u^{*}(\cdot,t)$ in the function space for every $t \in [0,T]$. 
%\zhou{Is it better to use $\theta(t)$ than $\theta_t$?} \gaby{Initially we had done $\theta(t)$ but I think it ends up being a bit cluttered notationally in that case.}
%
For example, if we consider $L^{2}(\Omega)$ as the function space, by closeness we mean $\| u_{\theta_{t}} - u^{*}(\cdot,t)\|_{L^2(\Omega)}$ is small for all $t$ (hereafter we denote $\|\cdot\|_{p} = \|\cdot\|_{L^{p}(\Omega)}$ for notation simplicity).
Notice that $\{u_{\theta_t}\,|\, 0 \le t \le T\}$ is a trajectory on $\Mcal$, whereas $u^{*}(\cdot,t)$ is a trajectory in the full space $L^{2}(\Omega)$. 

%
% , and thus the projection error plays an important role in the estimation of the difference between these two functions as time evolves.
%
% In most scenarios, we can derive the error bound in $W^{k,\infty}$, but we may only need $L^{\infty}$.

% Assume that the network structure is good enough, and we can represent $u^*(\cdot, t)$ using $u_{\theta_t}$ for some $\theta_t \in \Theta$ at any time $t \in [0,T]$, then the trajectory $\{\theta_t: 0 \le t \le T\}$ determines a path $\{u_{\theta_t}: 0 \le t \le T\}$ in $\Mcal$. However, since $u_{\cdot}: \Theta \to \Mcal$ is not injective, the inverse mapping may not render a single path in $\Theta$.

% Computational wise, the simplest approach is to minimize the least squares at different times. However it is  important to note that the solution is not unique since $\nabla_{\theta} u$ is not linearly dependent. Control should be done from accumulated sense.

\subsection{Proposed methodology}
\label{subsec:method}

Let $\Omega$ be an open bounded set in $\mathbb{R}^{d}$ and $F$ a \emph{nonlinear differential operator} of functions $u: \Omega \to \Rbb$ with necessary regularity conditions, we consider the IVP of the evolution PDE defined by $F$ with arbitrary initial value as follows:
\begin{equation}
\begin{cases}
\partial_t u(x,t) = F [u](x,t), & \ x \in \Omega,\ t \in (0,T],\\
% u(x,t)=0, & \ x \in \partial \Omega,\ t \in [0,T]\\
u(x,0)=g(x), & \ x \in \Omega,
\end{cases}
\label{eq:pde}
\end{equation}
where $T>0$ is some prescribed terminal time, and $g:\Rbb^d \to \Rbb$ stands for an initial value. For ease of presentation, we assume zero Dirichlet boundary condition $u(x,t) = 0$ for all $x\in \barOmega$ and $t \in [0,T]$ (for compatibility we henceforth assume $g(x)$ has zero trace on $\partial \Omega$) throughout this paper.
%\ye{``For ease of presentation, we assume zero Dirichlet boundary condition $u(x,t) = 0$ for all $x\in \barOmega$ and $t \in [0,T]$ (for compatibility we henceforth assume $g(x)$ has zero trace on $\partial \Omega$) throughout this paper''. Same statement at the beginning of Sec 3.3.}. 
We denote $u^{g}$ the solution to the IVP \eqref{eq:pde} with this initial $g$. 
The solution operator $\Scal_{F}$ of the IVP \eqref{eq:pde} is thus the mapping from the initial $g$ to the solution $u^{g}$ :
\begin{equation}
\label{eq:so}
    \Scal_{F}: C^{2}(\bar{\Omega}) \to C^{2,1}(\bar{\Omega}\times[0,T]), \quad \mbox{such that} \quad g \mapsto \Scal_{F}(g) := u^{g},
\end{equation}
where $C^{2}(\bar{\Omega}):= C(\bar{\Omega}) \cap C^2(\Omega)$ for short.
Our goal is to find a numerical approximation to $\Scal_{F}$. Namely, \emph{we want to find a fast computational scheme $\Scal_{F}$ that takes any initial $g$ as input and accurately estimate $u^{g}$ with low computation complexity.}

It is important to note the substantial difference between solving \eqref{eq:pde} for any given but fixed initial value $g$ and finding the solution operator \eqref{eq:so} that maps any $g$ to the corresponding solution $u^{g}$.
In the literature, most methods are developed for solving IVP \eqref{eq:pde} with a fixed $g$, such as traditional finite difference and finite element methods, as well as many state-of-the-art machine learning based methods.
%
% In particular, the mesh free approach of the machine learning methods have been shown to be more suitable to approximately solving the PDE problem when $\Omega$ is a subset of a high dimensional space \cite{han2018high-pde}. 
%
% However, we stress again that modern machine learning methods such as PINN, WAN, and DeepRitz all suffer from needing to be retrained whenever the initial condition $u_0$ changes \cite{raissi2019physics-informed,e2018deep,zang2020weak}. We now seek to explain how we overcome this issue.
%
However, these methods are computationally expensive if \eqref{eq:pde} must be solved with many different initial values, and they need to start from scratch for every new $g$.
In a sharp contrast, our goal is to find an approximation to the solution operator $\Scal_{F}$ which, once found, can help us to compute $u^{g}$ for any given $g$ at relatively much lower computational cost.

For ease of presentation, we use autonomous, second-order nonlinear differential operators $F[u] = F(x, u, \nabla_{x} u, \nabla_{x}^{2} u)$ as an example and take $\Omega = (0,1)^{d}$ in \eqref{eq:pde} to describe our main idea below.
Extensions to general non-autonomous nonlinear differential operators and PDEs defined on open bounded set $\Omega \subset \Rbb^{d}$ with given boundary values will be discussed in Section \ref{sec:variation}. 
% \ye{How about say ``and take $\Omega =(0,1)^d \subset \Rbb^{d}$ with zero Dirichlet boundary condition in \eqref{eq:pde} to describe our main idea below.''}
%
%In \eqref{eq:pde}, we considered the case where the spatial domain $\Omega$ is set to $\Rbb^{d}$.
%
%In practice, we may often encounter evolution PDEs defined on open, bounded set $\Omega \subset \Rbb^{d}$ with some regularity conditions on its boundary geometry, and certain boundary value functions are imposed as well. 
%
%We will also extend our approach to these cases in Section \ref{subsec:variation}. For now, we focus on \eqref{eq:pde} with $\Omega = \Rbb^{d}$.

To approximate the solution operator $\Scal_{F}$ in \eqref{eq:so}, we propose \emph{a control mechanism in the parameter space $\Theta$ of a prescribed reduced-order model $u_{\theta}$}.
Specifically, we first determine a reduced-order model $u_{\theta}$ to represent solutions of the IVP. 
%
% Here $\theta \in \Theta \subset \Rbb^m $ represents the $m$ trainable parameters of the neural network. Commonly, the parameter space contains weight and bias vectors of the neural network. 
%
We allow any parametric form of $u_{\theta}$ but only assume that $u_{\theta}(x) = u(x;\theta)$ is $C^1$ smooth with respect to $\theta$. 
This is a mild condition satisfied by the commonly used reduced-order models: if $\ut$ is a linear combination of basis functions and $\theta$ represents the combination coefficients, then $\ut$ is linear and hence smooth in $\theta$;
and if $\ut$ is a DNN as in \eqref{eq:u_param}, then $\ut$ is smooth in $\theta$ as long as all activation functions $\sigma$ are smooth.
Suppose there exists a trajectory $\{\theta_t\, | \, 0 \le t \le T\}$ in the parameter space $\Theta$ such that its corresponding $\utt$ approximates the solution of the IVP, we must have
\begin{equation}
\begin{cases}
\partial_t \utt(x) = \nabla_{\theta}u(x;\theta_{t}) \cdot \dot{\theta}_t = F[\utt](x), &\ \forall\, x \in \Omega,\ t \in (0,T],\\
u_{\theta_{0}}(x)=g(x), & \ \forall\,x \in \Omega.
\end{cases}
\label{eq:nn-pde}
\end{equation}
To compute $\utt$, it is sufficient to find a control vector (velocity) field $V_{F}:\Theta \to \Rbb^{m}$, in the sense of $\dot{\theta}_t = V_{F}(\theta_t)$, that steers the trajectory $\theta_{t}$ along the correct direction starting from the initial $\theta_{0}$ satisfying $u_{\theta_{0}}(x)=g(x)$. 

This observation suggests a new approach to solve the IVP with a fixed evolution PDE but varying initial values $g$: for the evolution equation in \eqref{eq:nn-pde} to hold, it suffices to find a vector field $V_{F}$ such that
\begin{equation}
\label{eq:V_condition}
    \nabla_{\theta} \ut \cdot V_{F}(\theta) = F[\ut]
\end{equation}
for all $\theta \in \Theta$. 
It is important to note that $V_{F}$ only depends on the nonlinear differential operator $F$ of the original evolution PDE, but not any actual initial value $g$ of the IVP.
Once this is achieved, we can effectively approximate the solution of the IVP with any initial value $g$: we first set $\theta_{0} = \theta^{g}$, where $\theta^{g}$ denotes the parameters such that $u_{\theta^{g}}$ fits $g$, then we numerically solve the following ODE in the parameter space $\Theta$ (which can be fast) using the control vector field $V_{F}$:
\begin{equation}
\begin{cases}
\dot{\theta}_t = V_{F}(\theta_{t}),& \ \forall\, t \in (0,T],\\
\theta_{0} = \theta^{g}.
\end{cases}
\label{eq:ode}
\end{equation}
The solution trajectory $\{\theta_t\, | \,  0 \le t \le T\}$ of the ODE \eqref{eq:ode} induces a path $\{\utt\, | \, 0 \le t \le T\}$ in $\Mcal$ as an approximation to the solution of the IVP. 
%
% If the initial value $g$ changes, we only need to find the initial parameter $\theta^{g}$ by fitting $u_{\theta^{g}}$ to $g$, and then solve \eqref{eq:ode} numerically which requires very low computation complexity.
%
The computational cost is thus composed of two parts: finding the parameters $\theta^{g}$ of $\ut$ to fit $g$ and numerically solving the ODE \eqref{eq:ode}, both of which are substantially cheaper than solving the IVP \eqref{eq:pde}.

\begin{figure}
\centering
\includegraphics[width=.75\textwidth]{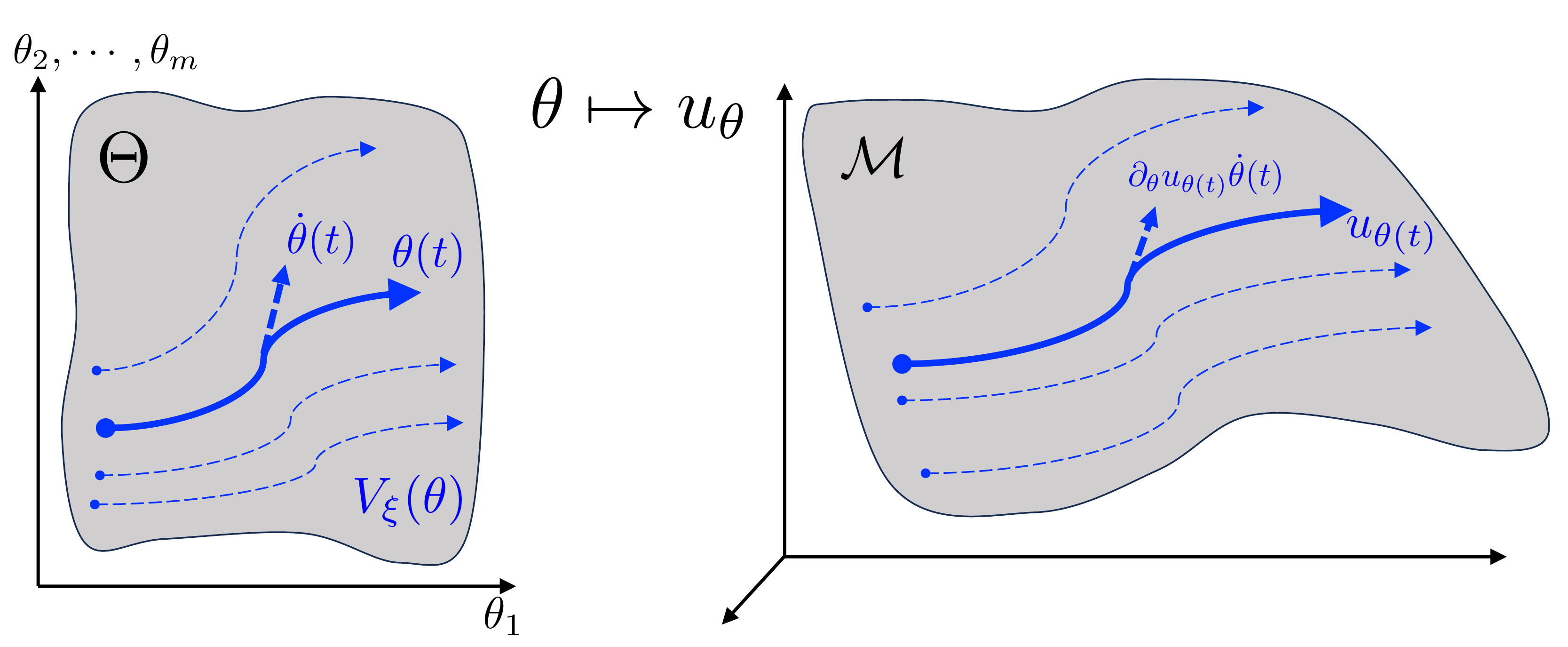}
\caption{Schematic plot of pulling back trajectories (solid and dashed blue curves) in $\Mcal=\{u_{\theta}:\theta \in \Theta\}$ to trajectories in the parameter space $\Theta$. Here each trajectory in $\Mcal$ represents the reduced-order model (e.g., DNN) $u_{\theta(t)}(\cdot)$ approximating the PDE solution $u^{*}(t,\cdot)$ starting from a given initial, and it is pulled back to the trajectory $\theta(t)$ (we use $\theta(t):=\theta_t$ as a trajectory here to avoid confusion with components $\theta_1,\dots,\theta_m$) in $\Theta$; and $V_{\xi}$ is a DNN approximating the control vector field $V_{F}$ in $\Theta$.}
\label{fig:control-pde}
\end{figure}

The main question is how to get the control vector field $V_{F}$ in \eqref{eq:ode}. 
%
%As finding an explicit form of $V$ is complicated in general, we choose to parameterize $V$ by another neural network $V_{\xi}$. The network $V_{\xi}$ is what we call our Neural Control. The training of this neural network represents the optimization problem we seek to solve.
%
As an explicit form of $V_F$ is unknown, we choose to express $V_{F}$ in a general parametric form $V_{\xi}$ with parameters $\xi$ to be determined. 
%
% The training of this neural network is explained in the next section, and represents the optimization problem we seek to solve.
%
Specifically, we propose to set $V_{\xi}$ as another DNN where $\xi$ represents the set of learnable network parameters in $V_{\xi}$.
A schematic plot of the pullback mechanism and the control vector field in $\Theta$ is provided in Figure \ref{fig:control-pde}.
We call $V_{\xi}$ the \emph{neural control} field. 
We learn the parameters $\xi$ by minimizing the following loss function:
\begin{equation}
\ell(\xi):=\int_{\Theta} \|\nabla_{\theta}\ut \cdot V_{\xi}(\theta)- F[\ut]\|_{2}^{2}\, d\theta.
% \ell(\xi):=\int_{\Theta}\int_{\Omega} |\nabla_{\theta}u_{\theta}(x,\theta) \cdot F_{\xi}(\theta)+N(u_{\theta}(x,\theta))|^2dxd\theta
\label{eq:total_loss}
\end{equation}
In practice, we approximate the integral in $\ell$ by Monte Carlo integration. We sample $K$ points $\{\theta_k\,|\, k = 1,\dots,K\}$ uniformly from $\Theta$ (here the subscript $k$ in $\theta_k$ stands for the $k$th point among the $K$ points sampled in $\Theta$) and form the empirical loss function 
\begin{equation} 
\hat{\ell}(\xi) = K^{-1} \cdot \sum_{k=1}^{K} \|\nabla_{\theta} u_{\theta_{k}} \cdot V_{\xi}(\theta_{k})- F[u_{\theta_{k}}]\|_{2}^{2}
\end{equation}
Then we minimize $\hat{\ell}(\xi)$ with respect to $\xi$, where the $L^2$ norm is also approximated by Monte Carlo integration on $\Omega$.
%
% To train the network $V_{\xi}$ we aim to minimize the loss function
% \begin{equation}
% \ell(\xi):=\int_{\Theta}\int_{\Omega} |\nabla_{\theta}u(x,\theta) \cdot V_{\xi}(\theta)-F(u(x,\theta))|^2dxd\theta
% \label{eq:total_loss}
% \end{equation}
% In the case where $F [\ut] \in \text{span} \{ \partial_{\theta_i} \ut \ : \ i=1,\ldots, m\}=:T_{\theta} u$ then the integrand of \eqref{eq:total_loss} can be 0. As such, the loss function \eqref{eq:total_loss} is trying to measure the total projection error of the function $F [\ut]$ onto the $T_{\theta}u$ spaces. In general, as the integral in \eqref{eq:total_loss} will not be computable, we can use Monte Carlo integration techniques to attempt to approximate it. We can sample points in the domain $\Omega$ to form an empirical loss which we call $\hat{\ell}$. We further flesh out how the training is accomplished in Section \ref{sec:numerical-results}. Generally, we will only seek to have an approximately low error which we set by the tolerance $\bar{\varepsilon}$, then we set a max iteration number $K$ and a fixed learning rate $\zeta$ for use in a standard gradient descent approach. 
%
The training of $V_{\xi}$ is summarized in Algorithm \ref{alg:neural-control}.

\begin{algorithm}
\caption{Training neural control $V_{\xi}$}
\label{alg:neural-control}
\begin{algorithmic}[1]
\Require{Reduced-order model structure $\ut$ and parameter set $\Theta$. Control vector field structure $V_{\xi}$. Error tolerance ${\varepsilon}$.}
\Ensure{Optimal control parameters ${\xi}$.}
\State Sample $\{\theta_{k}\}_{k=1}^{K}$ uniformly from $\Theta$ and $\{x_{n}\}_{n=1}^{N}$ from $\Omega$.
\State Form empirical loss $\hat{\ell}(\xi)$ as in \eqref{eq:empirical_loss}.
\State Minimize $\hat{\ell}$ with respect to $\xi$ using any optimizer (e.g., ADAM or AdaGrad) until $\hat{\ell}(\xi) \le \varepsilon$.
% \While{$k \leq K$ or $\hat{\ell}(\xi_k)>\bar{\varepsilon}$}
%     \State{$\xi_{k+1}\gets \xi_k-\zeta\nabla_{\xi}\hat{\ell}(\xi)$}
%     \State{$k \gets k+1$}
% \EndWhile
\end{algorithmic}
\end{algorithm}

%\subsubsection{Solving IVPs using learned neural control}

Once we trained the vector field $V_{\xi}$, we can implement the solution operator $\Scal_{F}$ in the following two steps: we first find a $\theta_{0}$ such that $u_{\theta_0}$ fits $g$, i.e., find $\theta_{0}$ that minimizes
%
%the $L^2$ loss between our neural network $\ut$ and $g$ 
%
$\| u_{\theta} - g \|_2$. This can be done by sampling $\{x_{n}\}_{n=1}^{N}$ from $\Omega$ and minimizing the empirical squared $L^2$ norm $(1/N)\cdot \sum_{n=1}^{N} |\ut(x_n) - g(x_n)|^2$ with respect to $\theta$. Then we solve the ODE \eqref{eq:ode} using any numerical ODE solver (e.g., Euler, 4th order Runge-Kutta, predictor-corrector) with $\theta_{0}$ as the initial value. Both steps can be done efficiently and the total computational cost is substantially lower than that of solving the original IVP \eqref{eq:pde} again. We summarize how neural control solves IVPs in Algorithm \ref{alg:solving-ivp}. Further details on the practical implementation of Algorithm \ref{alg:neural-control} and \ref{alg:solving-ivp} are discussed in Section \ref{sec:numerical-results}.

\begin{algorithm}
\caption{Implementation of solution operator $\Scal_{F}$ of the IVP \eqref{eq:pde} using trained control $V_{\xi}$}
\label{alg:solving-ivp}
\begin{algorithmic}[1]
\Require{Initial value $g$ and tolerance $\varepsilon_{0}$. Reduced-order model $\ut$ and trained neural control $V_{\xi}$.}
\Ensure{Trajectory $\hatthetat$ such that $u_{\hatthetat}$ approximate the solution $\Scal_{F}[g]$ of the IVP \eqref{eq:pde}.}
\State{Compute initial parameters $\theta_{0}$ such that $\|u_{\theta_{0}} - g \|_2 \le \varepsilon_{0}$.}
% \While{$k \leq K$ or $\|\ut-g\|^2_{2}>\bar{\varepsilon}$}
%     \State{$\theta_{k+1}\gets \theta_k-\zeta \nabla_{\theta}\|\ut-g\|^2_{2}$}
%     \State{$k \gets k+1$}
% \EndWhile
% \State{$t \gets 0$}
% \State{$\theta(0)=\theta_k$}
\State{Use any ODE solver to compute $\hatthetat$ to solve \eqref{eq:ode} with approximate field $V_{\xi}$ and initial $\theta_{0}$.}
% \While{$t < T$}
%     \State{$\theta(t+h)\gets \theta(t) + h V_{\xi}(\theta(t)) $}
%     \State{$t \gets t+h$}

% \EndWhile
\end{algorithmic}
\end{algorithm}

% From our method it is not hard to see how one can return any desired time-step $t$ by terminating at that time, or simply collecting the desired parameters as the algorithm steps through.

\subsection{Error analysis}
\label{subsec:error-analysis}

In this subsection, we develop an error estimate of the proposed method.
We first focus on the error due to projection onto the tangent space $T_{\ut}\Mcal$ in the $L^2$ space in Section \ref{subsubsec:proj-error}.
Then we establish the solution approximation error for linear and semilinear parabolic PDEs in Section \ref{subsubsec:nonlinear-error}.
For ease of discussion, we again assume zero Dirichlet boundary condition $u(x,t) = 0$ for all $x\in \barOmega$ and $t \in [0,T]$, and we let $\Omega = (0,1)^{d} \subset \Rbb^{d}$ be the unit open cube in $\Rbb^{d}$ and $\Theta$ some open bounded set in $\Rbb^{m}$ (note that our analysis below applies as long as $\Omega$ is open and bounded). 
%\gaby{ I propose adding this here if not added in \eqref{eq:pde}, ``We shall assume the PDE has a zero Dirichlet boundary condition and our DNN satisfies $\ut(x)=0$ for all $x \in \partial \Omega$."} %A reduced-order model $\ut(x) = u(x;\theta) \in C( \bar{\Omega} \times \bar{\Theta})$.
%
% In all that follows we suppose $\Omega$ is $(0,1)^d$, and will analyze the following PDE problem
% \[ (PDE)
% \begin{cases}
% \partial_t u(x,t)+\mathcal{F}(u(x,t))=0, & \ x \in \Omega,\ t \in (0,T]\\
% u(x,t)=0, & \ x \in \partial \Omega,\ t \in [0,T]\\
% u(x,0)=g(x), & \ x \in \Omega
% \end{cases}
% \] 
We let $F[u]:=F(u,\nabla u, \nabla^2 u)$ be a nonlinear differential operator with necessary regularity conditions to be specified later and allows for a unique solution to the PDE for each initial. 
Additional requirements on the regularity of $\ut$ will be given when needed.
%
% If we wish for $F$ to be explicitly linear then we shall use $\Lcal$ instead. 
%
% We let $u(x,\theta)$ represent some fixed neural network structure with parameter space $\theta \in \Theta\subset \Rbb^m$ such that $u(x,\theta)=0$ for all $x \in \partial \Omega$ and $\Theta$ is compact.
%
% We will denote by $\|\cdot\|_{2}$ the $L^2(\Omega)$ norm, $|\cdot |$ the Euclidean norm in $\mathbb{R}^d$ or the spectral norm in the case of matrices, and $\| \cdot \|_{\infty}$ as the $L^{\infty}(\Omega)$ norm.

\subsubsection{Approximation error of control vector field}
\label{subsubsec:proj-error}
We first investigate the main source of error when using a reduced-order model to approximate the time-evolving solution of the given PDE. We show that this error is due to the imperfect representation of $F[\ut]$ using $\nabla_{\theta} \ut$ in \eqref{eq:V_condition}.
Specifically, due to the approximation of reduced-order models, $T_{\ut}\Mcal$ is only a finite-dimensional subspace of $L^{2}$, and thus we can only approximate the projection of $F[\ut]$ onto this tangent space.
We will need the following assumptions on the regularity of $\ut$ and $F$.
\begin{assumption}
The reduced-order model $u_{\theta}(\cdot ) \in C^{3}(\Omega) \cap C(\bar{\Omega})$ for every $\theta \in \bar{\Theta}$ and $u(x;\cdot) \in C^{2}(\Theta) \cap C(\bar{\Theta})$. Moreover, there exists $L>0$ such that for all $\theta \in \bar{\Theta}$
\begin{equation}
\label{eq:def_F}
F[\ut] \in \Fcal^{L} := \{f \in C^{1}(\Omega) \cap C(\bar{\Omega}): \|f \|_{\infty} \le L, \ \|\nabla f\|_{\infty} \le L \}. 
\end{equation}
%
% $\partial_{\theta_i}u_{\theta}\in C^{2,1}(\Omega\times \Rbb^m)$ for all $\theta = (\theta_{1},\dots,\theta_{m}) \in \Theta$.
\label{assump:regularity}
\end{assumption}

Assumption \ref{assump:regularity} provides some sufficient regularity conditions on the reduced-order model $\ut$ and boundedness of $F[\ut]$ and its gradient to be used in our error estimates.
Notice that we consider $F$ as second-order differential operator here and therefore the assumption $\ut \in C^{3}(\Omega)$ ensures that $\ut(x),\nabla \ut(x), \nabla^2 \ut(x)$ are all sufficiently smooth. The regularity condition on $F$ in Assumption \ref{assump:regularity} requires that the mapping $F[\ut](x)$ is a $C^1$ function and have magnitudes and gradients bounded by $L$ over $\bar{\Omega}$.
These assumptions are generally mild as we will use reduced-order models smooth in $(x,\theta)$, e.g., a DNN with smooth activation functions, and the operator $F$ is sufficiently regular.

\begin{assumption}
For any $\bar{\varepsilon} >0$, there exist a reduced-order model $\ut$ and a bounded open set $\Theta \subset \Rbb^{m}$, such that for every $\theta \in \barTheta$ there exists a vector $\alpha_{\theta} \in \Rbb^{m}$ satisfying
\[
\|\alpha_{\theta} \cdot \nabla_{\theta}\ut-F[\ut]\|_{2}\leq \bar{\varepsilon}.
\]
\label{assump:tangent_error}
\end{assumption}
% \ye{$v_{\theta}$ is also used in Lemma 3.2 but has different meaning. Better change $v_{\theta}$ to something else (such as $\alpha_{\theta}$ which is used in Lemma 3.2) in Assumption 2 and Examples 3.1 and 3.2.}
% \gaby{Ammended}

Assumption \ref{assump:tangent_error} provides an upper bound on the error when projecting $F[\ut]$ onto the tangent space $T_{\ut}\Mcal$, which is spanned by the functions in $\nabla_{\theta} \ut$. This error bound is determined by the choice of the reduced-order model $\ut$ and the parameter set $\Theta$. As will be demonstrated in our numerical experiments, a small projection error can be achieved by using a standard DNN as reduced-order model $\ut$. As such error is difficult to analyze due to the complex structures of general DNNs. 
We provide an example reduced-order model with special structure to justify the reasonableness of Assumption \ref{assump:tangent_error}.

\begin{example}
\label{ex:galerkin}
Let $\barvarepsilon>0$ and $\{\varphi_{j}\}_{j=1}^{\infty}$ be a complete smooth orthonormal basis (e.g., generalized Fourier basis) for $L^2(\Omega)$. Suppose there exist $C>0$, $\gamma>1$, and $C_0>0$ such that for all $u \in C^3(\Omega) \cap C(\barOmega)$ and $\|u \|_{2}^2 \le C_0$ we have
%
\begin{equation}
\label{eq:ex-funclass-G}
F[u] \in \Gcal^{C,\gamma} := \set{f \in C^{1}(\Omega) \cap C(\bar{\Omega}):  |\langle f, \varphi_j \rangle|^{2} \le C j^{-\gamma},\ \forall\,j \ge 1}.  
\end{equation}
%
Then there exists $m = m(\barvarepsilon,C,\gamma) \in \Nbb$ such that $\sum_{j=m+1}^{\infty} C j^{-\gamma} < \barvarepsilon^2$. 
%
Consider $\ut = \theta \cdot \varphi = \sum_{j=1}^{m} \theta_{j} \varphi_{j}$. 
We denote $f_{\theta}:= F[\ut]$ for short. 
%
Then $\nabla_{\theta} \ut = \varphi= (\varphi_1,\dots,\varphi_{m})$ and for $\alpha^{f_{\theta}} = (\alpha_1^{f_{\theta}},\dots, \alpha_{m}^{f_{\theta}})$ with $\alpha_{j}^{f_{\theta}} := \langle f_{\theta}, \varphi_{j} \rangle$, there is
\begin{equation*}
    \| \alpha^{f_{\theta}} \cdot \nabla_{\theta} \ut - F[\ut] \|_2^2 = \Big\| \sum_{j=1}^{m} \alpha_{j}^{f_{\theta}} \varphi_{j} - f_{\theta} \Big\|_2^2 = \sum_{j=m+1}^{\infty} |\langle f_{\theta}, \varphi_{j} \rangle |^2 \le \barvarepsilon^2. 
\end{equation*}
Therefore, the reduced-order model $\ut = \theta \cdot \varphi$ with $\Theta = \{\alpha \in \Rbb^m : |\alpha|^2<C_0\}$ and $\alpha_{\theta} = \alpha^{f_{\theta}}$ satisfy Assumption \ref{assump:tangent_error}.

This example can be modified to use a more general form of reduced-order model $\ut$, such as a DNN. To see this, we first repeat the procedure above but with $\barvarepsilon$ replaced by $\barvarepsilon/2$. Then the universal approximation theorem \cite{hornik1991approximation, yarotsky2017error} and the continuity of DNNs in its parameters imply that there exist DNNs $\{\hatvarphi_j: 1 \le j \le m\}$, whose network parameters are collectively denoted by $\eta \in \Rbb^{m'}$, satisfy $\| \hatvarphi_j - \varphi_j\|_{\infty} \le \barvarepsilon/(2\sqrt{m C_0 |\Omega|})$ and hence $\| \hatvarphi_j - \varphi_j\|_{2} \le \barvarepsilon/(2\sqrt{m C_0})$ for all $\eta$ in an open set $H \subset \Rbb^{m'}$. 
%
% Further, there exists an open set $H \subset \Rbb^{|\eta|}$ such that 
%
Consider the DNN $\ut = c\cdot \hatvarphi$ with parameters $\theta = (c,\eta) \in \Rbb^{n}$ where $n = m+m'$. Then $\nabla_{c} \ut(x) = (\hatvarphi_1,\dots,\hatvarphi_{m})$. Using the example above, we know for any $f_{\theta} := F[\ut] \in \Gcal^{C,\gamma}$, there exists $\alpha^{f_{\theta}} \in \Rbb^{m}$ such that $\| \alpha^{f_{\theta}} \cdot \varphi - F[\ut] \|_{2} \le \barvarepsilon/2$. 
Therefore, we use $(\alpha^{f_{\theta}},0)$ which concatenates $\alpha^{f_{\theta}}$ and $0 \in \Rbb^{m'}$ as the combination coefficients of $\nabla_{\theta}\ut$ to obtain
\begin{align*}
    \| (\alpha^{f_{\theta}},0) \cdot \nabla_{\theta} \ut - F[\ut]\|_{2}
    & = \| \alpha^{f_{\theta}} \cdot \nabla_{c} \ut - F[\ut]\|_{2} \\
    & \le \|\alpha^{f_{\theta}} \cdot \hatvarphi - \alpha^{f_{\theta}} \cdot \varphi \|_{2} + \| \alpha^{f_{\theta}} \cdot \varphi - F[\ut] \|_{\infty} \\
    & \le \sum_{j=1}^{m} |\alpha_{j}^{f_{\theta}}| \| \hatvarphi_j - \varphi_j\|_{2} + \frac{\barvarepsilon}{2} \\
    & \le \sqrt{m C_0}\cdot \frac{\barvarepsilon}{2\sqrt{m C_0}} + \frac{\barvarepsilon}{2} \\
    & = \barvarepsilon.
\end{align*}
Therefore, the DNN $\ut = c\cdot \hatvarphi$ with $\Theta = \{(c,\eta): |c_{j}|^2 < C_0,\ \eta \in H\}$ and $\alpha_{\theta} = (\alpha^{f_{\theta}},0)$ satisfy Assumption \ref{assump:tangent_error}.
%
% For any $m>0$ we define the neural network $\ut^{(m)}: \Omega \times \Rbb^m \to \Rbb$ by $\ut^{(m)}=\sum_{i=1}^m \theta_i e_i$ where $\theta=(\theta_1,\cdots, \theta_m) \in \Rbb^m$. Then We select $\Theta:=\{\theta \in \Rbb^m \ \ : \ \sum_{i=m+1}^{\infty}|\langle e_i, F[\ut^{(m)} ] \rangle| < \varepsilon\}$. Essentially, these are the functions $F[\ut^{(m)}]$ whose coefficients in the $\{e_i\}_{i=1}^{\infty}$ basis decay sufficiently fast. Notice as $m \to \infty$ then we have $\Theta \to \Rbb^{m}$. Then for all $\theta \in \Theta$ if we set $\alpha=(\alpha_1,\cdots,\alpha_m)$ by $\alpha_i=\langle e_i, F[\ut^{(m)}] \rangle$ then
% \[
% \|\alpha \cdot \nabla_{\theta}\ut^{(m)}-F[\ut^{(m)}]\|_{2}=\|\sum_{i=1}^{m}\alpha_i e_i - \sum_{i=1}^{\infty}\langle e_i, F[\ut^{(m)}] \rangle e_i\|_2\leq\sum_{i=m+1}^{\infty}|\langle e_i, F[\ut^{(m)}] \rangle | < \varepsilon.
% \]
% This example is closely related to the many traditional reduced-order methods, such as the spectral or Galerkin-type methods, of which our method can be seen as a generalization.
\end{example}

% \begin{remark}
% While the above examples are somewhat arbitrary. As we see in the numerical experiments, a far more general setting than the ones used in the examples also provide good results. This suggests Assumption \ref{assump:tangent_error} can be satisfied in more general settings as well.
% \end{remark}

Before proving the main proposition of this section we will need the following lemma.

\begin{lemma}
Suppose Assumption \ref{assump:regularity} and \ref{assump:tangent_error} are satisfied. For all $\varepsilon>\bar{\varepsilon}$ there exists $v: \bar{\Theta} \to \Rbb^m$ such that $v$ is bounded over $\bar{\Theta}$ and the value of $v$ at $\theta$, denoted by $v_{\theta}$, satisfies
\[
\|v_{\theta} \cdot \nabla_{\theta}\ut-F[\ut]\|_{2}\leq \varepsilon, \qquad \forall\, \theta \in \barTheta.
\]
\label{lem:bounded}
\end{lemma}
\begin{proof}
Let $\varepsilon>\bar{\varepsilon}$ and $\delta \in (0,\varepsilon-\bar{\varepsilon})$. By Assumption \ref{assump:tangent_error}, for all $\theta \in \Theta$ there exists $\alpha_{\theta}\in \Rbb^m$ coefficient such that
\[
\| \alpha_{\theta}\nabla_{\theta}\ut-F [\ut]\|_{2} \leq \bar{\varepsilon}.
\] 
As $F [\ut]$ and $\nabla_{\theta}\ut$ are continuous in $\theta$ and $\Omega$ is bounded, we associate to each $\theta$ and coefficient $\alpha_{\theta}$ the open set $U_{\theta}$ containing $\theta$, small enough, such that for all $\theta' \in U_{\theta}$ we have
\begin{equation}
    \| \alpha_{\theta}\nabla_{\theta}u_{\theta'}-\alpha_{\theta}\nabla_{\theta}\ut\|_{2} +\|F [\ut]-F [u_{\theta'}]\|_{2} \leq \delta
\end{equation} 
and hence 
\begin{equation}
\begin{aligned}
    \|\alpha_{\theta}\cdot \nabla_{\theta}u_{\theta'}-F [u_{\theta'}]\|_{2} & \leq \| \alpha_{\theta}\nabla_{\theta}u_{\theta'}-\alpha_{\theta}\nabla_{\theta}\ut\|_{2}+ \| \alpha_{\theta}\nabla_{\theta}\ut-F [\ut]\|_{2}+\|F [\ut]-F [u_{\theta'}]\|_{2}
    % &\leq \delta + \| \alpha_{\theta}\nabla_{\theta}\ut-F [\ut]\|_{2}\\
    \leq \delta + \bar{\varepsilon}.
\end{aligned}
\label{eq:thm-delta-bound}
\end{equation} 
Therefore $\cup_{\theta \in \bar{\Theta}}U_{\theta}$ is an open cover of $\bar{\Theta}$. As $\bar{\Theta}$ is compact 
%\gaby{can show result on in the case $\Theta$ is just open bounded as well. Just would use the closure here.}
this open cover has a finite subcover $\cup_{i=1}^N U_{\theta_i}$ for particular $\theta_i$'s. Define $v:\bar{\Theta} \to \Rbb^m$ such that $v_{\theta}:=v(\theta)=\alpha_{\theta_i}$ if $\theta \in U_{
\theta_i}$ (if $\theta$ is in the intersection of multiple $U_{\theta_i}$'s we choose a single $\alpha_{\theta_i}$ arbitrarily). We see from this construction that
$v_{\theta}$ is uniformly bounded over $\bar{\Theta}$ as the range of $v_{\theta}$ is finite. From \eqref{eq:thm-delta-bound} we have 
\[
\|v_{\theta} \cdot \nabla_{\theta}\ut-F[\ut]\|_{2}\leq \delta + \barvarepsilon \le \varepsilon.
\]
\end{proof}

% Realize that Lemma \ref{lem:bounded} only shows the existence of a bounded vector field. As we want to approximate a vector field on $\Theta$ by a Neural Network we need to show there exists a continuous vector field that accomplishes the same. This is difficult to prove in general as assumption \ref{assump:tangent_error} does not guarantee a continuous way to select the coefficients of $\nabla_{\theta}\ut$. This is because  the mapping $\ut : \Theta \to L^2$ which creates the parameter submanifold is in general only a submersion and not an immersion. This can clearly be seen from the dimension-dropping issue as was discussed in Section \ref{subsec:model}. With this in mind, we will present the following relaxation for the control vector field.

With Assumptions \ref{assump:regularity} and \ref{assump:tangent_error}, and Lemma \ref{lem:bounded} we can prove the existence of an accurate neural control field $V_{\xi}$ parameterized as a neural network, as shown in the next proposition.

\begin{proposition}
Suppose Assumption \ref{assump:regularity} and \ref{assump:tangent_error} hold. Then for any $\varepsilon > 0$, there exists a differentiable vector field parameterized as a neural network $V_{\xi}: \bar{\Theta} \to \Rbb^{m}$ with parameters $\xi$, such that 
\[
\|V_{\xi}(\theta) \cdot \nabla_{\theta}\ut-F[\ut]\|_{2} \leq \varepsilon,
\]
for all $\theta \in \barTheta$.
\label{prop:F_exists}
\end{proposition}
\begin{proof}
We first show that there exists a differentiable vector-valued function $V: \bar{\Theta} \to \Rbb^{d}$ such that 
\begin{equation}
\label{eq:V_bound}
    \|V(\theta) \cdot \nabla_{\theta}\ut-F[\ut]\|_{2} \leq \frac{\varepsilon}{2}
\end{equation}
for all $\theta \in \bar{\Theta}$. 
% Then the claim follows immediately by invoking the universal approximation theorem \cite{hornik1991approximation,petersen2018optimal}.
%
To this end, we choose $\bar{\varepsilon}_0 \in (0, \varepsilon/2)$ and $\bar{\varepsilon} \in (\bar{\varepsilon}_0, \varepsilon/2) $, then by Assumption  \ref{assump:tangent_error} and Lemma \ref{lem:bounded} we know that there exist a reduced-order model $\ut$, a bounded open set $\Theta \subset \Rbb^{m}$, and $M_v>0$ such that there is a vector-valued function $\theta \mapsto v_{\theta}$, where for any $\theta \in \bar{\Theta}$, we have $|v_{\theta}|<M_v$ and
\[
\|v_{\theta} \cdot \nabla_{\theta}\ut-F[\ut]\|_{2} \leq \barvarepsilon.
\]
Note that $v_{\theta}$ is not necessarily differentiable with respect to $\theta$.
To obtain a differentiable vector field $V(\theta)$, for each $\theta \in \bar{\Theta}$, we define the function $\psi_{\theta}$ by
\begin{equation*}
    \psi_{\theta}(w) : = \| w \cdot \nabla_{\theta} \ut - F[\ut] \|_2^2 = w^{\top} G(\theta) w - 2 w^{\top} p(\theta) + q(\theta),
\end{equation*}
%
% \begin{align*}
% \|a \cdot \nabla_{\theta}\ut-F[\ut]\|_{2}^2&=\| a \cdot \nabla_{\theta}\ut\|^2_{2}-2\langle a \cdot \nabla_{\theta}\ut,F[\ut]\rangle_{2} + \|F[\ut]\|_{2}^2\\
% &=a^{\top}\left(\int_{\Omega}\nabla_{\theta}\ut\nabla_{\theta}\ut^{\top} dx\right)a-2a^{\top}\left(\int_{\Omega}\nabla_{\theta}\ut F[\ut]dx\right)\\
% &\ \ \ \ \ \ \ \  + \|F[\ut]\|_{2}^2.\\
% &=a^{\top}A(\theta)a-2a^{\top}b(\theta)+c(\theta),
% \end{align*}
where
\begin{equation}
\label{eq:def-G}
G(\theta):=\int_{\Omega}\nabla_{\theta}\ut(x) \nabla_{\theta}\ut(x)^{\top}\, dx, \quad 
p(\theta):=\int_{\Omega}\nabla_{\theta}\ut(x) F[\ut](x)\,dx, \quad 
q(\theta):=\int_{\Omega}F[\ut](x)^2\,dx.
\end{equation}
%Define $f_{\theta}(a):\Rbb^m \to \Rbb$ by $f_{\theta}(a)=a^{\top}A(\theta)a-2a^{\top}b(\theta)+c(\theta)$.
%
Then we know 
\begin{equation}
\label{eq:psi_theta_star}
    \psi_{\theta}^{*} : = \psi_{\theta}(v_{\theta}) = \| v_{\theta} \cdot \nabla_{\theta} \ut - F[\ut] \|_2^2 \le \barvarepsilon^2.
\end{equation}
It is also clear that $G(\theta)$ is symmetric and positive semi-definite. Moreover, due to the compactness of $\bar{\Omega}$ and $\bar{\Theta}$, as well as that $\nabla_{\theta} u \in C(\barOmega \times \barTheta)$, we know there exists $\lambda_{G} > 0$ such that 
\begin{equation*}
    \| G(\theta) \|_2 \le \lambda_{G}
\end{equation*}
for all $\theta \in \bar{\Theta}$.
Therefore, $\psi_{\theta}$ is a convex function and the Lipschitz constant of $\nabla \psi_{\theta}$ is uniformly upper bounded by $\lambda_{G}$ over $\bar{\Theta}$.
Now for any $w \in \Rbb^{m}$, $h>0$, and $K \in \Nbb$ (we reuse the letter $K$ as the iteration counter instead of the number of sampling points in this proof), we define
\begin{equation*}
    \Ocal_{\theta}^{K,h}(w) := w_{K}, \quad \mbox{where} \quad w_k = w_{k-1} - h \nabla \psi_{\theta}(w_{k-1}), \quad w_0 = w, \quad k = 1,\dots,K.
\end{equation*}
Namely, $\Ocal_{\theta}^{K,h}$ is the oracle of executing the gradient descent optimization scheme on $\psi_{\theta}$ with step size $h>0$ for $K$ iterations.
%
% , and so $f_{\theta}$ is a convex function for all fixed $\theta$. Denote $a^*(\theta)=G(\theta)$ where $G(\theta)$ is as in Lemma \ref{lem:bounded} for constant $\varepsilon'$ where $\bar{\varepsilon}<\varepsilon'<\varepsilon$
%
% Let $a_k(\theta)=a_{k,\theta}$ where $a_{k,\theta}$ is defined by $a_{0,\theta}=0$ and $a_{k,\theta}=a_{k-1,\theta}-h\nabla f_{\theta}(a_{k-1,\theta})$ if $f_{\theta}(a_{i,\theta})>(\varepsilon')^2$, otherwise $a_{k,\theta}=a_{k-1,\theta}.$ This is just gradient descent on a convex function with stopping criteria. 
%
% As $A(\theta)$ is continuous and hence bounded over $\Theta$, we know the Hessian of $f_{\theta}(a)$ is bounded for all $a$ and $\theta \in \Theta$. 
%

Next, we slightly modify the standard convergence result of gradient descent in convex optimization \cite[Theorem 2.1.14]{nesterov1998introductory} and obtain Lemma \ref{lem:appendix} in Appendix \ref{appsec:proof}.
%\textcolor{red}{[Ye: Put the lemma and proof in an appendix (after conclusion and before references)]}\gaby{I have done this}. 
Notice that $\psi_{\theta}$ is convex, differentiable, and $\nabla \psi_{\theta}$ is Lipschitz continuous with Lipschitz constant upper bounded by $\lambda_{G}$.
Therefore, applying 
Lemma \ref{lem:appendix} with $y=v_{\theta}$, $f=\psi_{\theta}$, and the gradient descent scheme for $K$ iterations ($K$ to be determined soon) with initial 0 and any fixed step size $h \in (0, 1/\lambda_{G})$ to $\psi_{\theta}$ directly yields an error bound for $\psi_{\theta}(\Ocal_{\theta}^{K,h}(0))$:
%fixed step size $h \in (0, 1/\lambda_{G})$:
\begin{equation}
\label{eq:psiK_bound}
\psi_{\theta}(\Ocal_{\theta}^{K,h}(0)) - \psi_{\theta}^{*}  \le \frac{|0-v_{\theta}|^2}{2Kh}.  
\end{equation}
Combining this with the bound $|v_{\theta}|<M_v$, we choose any
\begin{equation*}
    K \ge \frac{ M_v^2}{2h((\varepsilon/2)^2 - \bar{\varepsilon}^2)},
\end{equation*}
and there is
\begin{equation}
\label{eq:psi_opt}
    \psi_{\theta}(\Ocal_{\theta}^{K,h}(0)) - \psi_{\theta}^{*}  \le \frac{ M_v^2}{2Kh} \le \Big( \frac{\varepsilon}{2} \Big)^2 - \bar{\varepsilon}^2.
\end{equation}

Notice that $\Ocal_{\theta}^{K,h}$ is a differentiable vector-valued function of $\theta$ because $K$ and $h$ are fixed.
Therefore, combining \eqref{eq:psi_theta_star} and \eqref{eq:psi_opt} yields
\begin{equation*}
    0\leq \psi_{\theta}(\Ocal_{\theta}^{K,h}(0)) = (\psi_{\theta}(\Ocal_{\theta}^{K,h}(0)) - \psi_{\theta}^{*}) + \psi_{\theta}^{*} \le (\varepsilon/2)^2 -\barvarepsilon^2 + \barvarepsilon^2 = (\varepsilon/2)^2.
\end{equation*}
As this inequality holds $\forall \theta \in \barTheta$, we set $V(\theta) = \Ocal_{\theta}^{K,h}(0)$ which is a differentiable function of $\theta$ satisfying \eqref{eq:V_bound}.

By the universal approximation theorem of neural networks \cite{guhring2020error} (see also Remark \ref{remark:guhring}), we know there exists a differentiable vector-valued function parameterized as a neural network $V_{\xi}$ with parameters $\xi$ such that 
\begin{equation*}
    |V_{\xi}(\theta) - V(\theta)|_{\infty} \le \varepsilon/(2B)
\end{equation*}
for all $\theta \in \barTheta$,
where $B:= \max_{\theta \in \barTheta} \|\nabla_{\theta} \ut\|_2 < \infty$ and $|\cdot|_{\infty}$ stands for the $\infty$-norm of vectors. Hence we know 
\begin{equation*}
    \|V_{\xi}(\theta)\cdot \nabla_{\theta} \ut - F[\ut] \|_2 \le \|V_{\xi}(\theta)\cdot \nabla_{\theta} \ut - V(\theta) \cdot \nabla_{\theta} \ut \|_2 + \|V(\theta) \cdot \nabla_{\theta} \ut - F[\ut] \|_2 \le B\cdot \frac{\varepsilon}{2B} + \frac{\varepsilon}{2} = \varepsilon.
\end{equation*}
This completes the proof.
\end{proof}

\begin{remark}
It is important to note the geometry of $\Mcal$, especially its dimensionality, is complex and highly dependent on the structure of $\ut$ and the parameter space $\Theta$. 
In particular, we can show that the tangent space $T_{\ut}\Mcal = \text{span}(\nabla_{\theta} \ut)$ at any $\ut \in \Mcal$ is in the $L^2$ space, where $\nabla_{\theta} \ut = (\partial_{\theta_1} \ut, \dots, \partial_{\theta_m} \ut)$ for $\theta = (\theta_1,\dots,\theta_m)$. 
(Here we use discrete indices $1,\dots,m$ as subscripts of $\theta$ to indicate its components for notation simplicity. This is to be distinguished from the subscript $t$ in $\theta_{t}$ which stands for time of the trajectory $\theta_{t}$.)
However, $\dim(T_{\ut}\Mcal)$ may vary across different $\ut$ on $\Mcal$.
For example, consider the reduced-order model $\ut$ parameterized as a DNN as in \eqref{eq:u_param}: when $w=0$, we have $\theta=(0,b,\cdots)$ and hence $\partial_{W_l} u_{\theta} = 0$ and $\partial_{b_l} u_{\theta} = 0$ for all $l=1,\dots,L$.
In this case, the $m$ components of $\nabla_{\theta} \ut$ are \emph{not} linearly independent, and $\dim(T_{\ut}\Mcal) < m$ for such $\theta$'s.
This distinguishes our parameter submanifold from existing ones, such as \cite{anderson2022evolution}, which assumes that the tangent space is always of full dimension $m$ at any point of the submanifold.
In our case, however, challenges and complications in dealing with the parameter submanifold $\Mcal$ can be avoided if we made such an assumption, but it will lead to incorrect analysis and error estimation, which poses a major technical challenge for the proposed framework.
%
% \zhou{This paragraph is harder to understand for readers. Many points discussed here seem to be easier after the proposed method and/or the analysis. Shall we move them to a later place?}
%
Specifically, we note that the rank of $G(\theta)$ varies across $\Theta$, and therefore the pseudoinverse $G(\theta)^+$ may be discontinuous.
% otherwise the proof of Proposition \ref{prop:F_exists} will be significantly simplified since $V(\theta)$ can be simply set to $G(\theta)^{+} \nabla_{\theta} \ut F[\ut]$.
%
A major theoretical merit of Proposition \ref{prop:F_exists} is that we can still ensure the existence of a differentiable control vector field in such case.

% . If the rank was full for all $\theta \in \Theta$ then $A^{-1}(\theta)$ exists and is differentiable. This would render the proof of Proposition \ref{prop:F_exists} to be trivial as the minimum of the $f_{\theta}$ function would be differentiable in $\theta$. Further in this case we can show $\ut$ defines an immersion in the sense of \eqref{eq:M}. However, as this may not always be the case, we have to resort to the relaxation in Proposition \ref{prop:F_exists} instead.

% If we considered our neural network as a Galerkin projection where only the output layer of the network changes, then the worrisome discontinuities naturally disappear. In our case, the more general neural networks will likely lack this same property. That being said, it is the opinion of the authors that application-specific neural network structures exist such that this issue is overcome without sacrificing the neural network expressiveness, we do not consider such structures here to keep our analysis as general as possible.
\end{remark}

\subsubsection{Error analysis in solving (semi-)linear parabolic PDEs}
\label{subsubsec:nonlinear-error}

Now we are ready to provide error bounds of our method in solving a large class of linear and semilinear parabolic PDEs. This class of PDEs covers many types of reaction-diffusion equations, such as heat equations, Fisher's equation or the Allen-Cahn equation. 
The differential operator $F$ in linear and semilinear parabolic PDEs has the form
\[
F [u]=\nabla \cdot (A \nabla u) + b \cdot \nabla u + f(u)
\]
where $A:\Omega \to \Rbb^{d\times d}$ and $b : \Omega \to \Rbb^d$ are continuous, $f : \Rbb \to \Rbb$ is $L_f$-Lipschitz and acts on $u(x)$ for each $x$. Moreover we assume that there exist $\lambda \ge 0$ and $B\ge 0$ such that
\begin{equation}
\label{eq:A-bound}
z^{\top}A(x) z \geq \lambda|z|^2, \quad \forall\, z \in \Rbb^d, \  x \in \Omega,    
\end{equation}
and
\begin{equation}
\label{eq:b-bound}
\|\nabla \cdot b \|_{\infty} \leq B.    
\end{equation}
% \label{assump:non-linear}
% \end{assumption}

Furthermore, due to the smoothness of $V_{\xi}$ and compactness of $\barTheta$, we know there exist $M_{V}>0$ and $L_{V}>0$ such that
\begin{equation}
\label{eq:MV}
    \max_{\theta \in \barTheta} |V_{\xi}(\theta)| \le M_{V} \qquad \mbox{and} \qquad \max_{\theta \in \barTheta} |\nabla_{\theta} V_{\xi}(\theta)| \le L_{V}.
\end{equation}

\begin{theorem}
% Suppose Assumptions \ref{assump:regularity} and \ref{assump:} hold. Further, assume
% for all  $\theta \in \Theta \subset \Rbb^m$ that
% \[
% \| \nabla_\theta \ut\cdot V(\theta)-F [\ut]\|_{2} \leq \bar{\varepsilon}.
% \] 
Suppose Assumptions \ref{assump:regularity} and \ref{assump:tangent_error} hold. Then there exist control field $V_{\xi}$ such that for any $u^*$ satisfying the evolution PDE in \eqref{eq:pde} there is
%Let $g \in C(\barOmega)$ be any initial value of the IVP \eqref{eq:pde} such that $\|u_{\theta_{0}}-g \|_{2} \leq \varepsilon_{0}$ for some $\theta_{0} \in \Theta$. Then there exists $V_{\xi}$ such that for $\theta_t$ solving the ODE \eqref{eq:ode} with $V_{\xi}$ and initial $\theta_0$, we have
\begin{equation}
\label{eq:utheta-ustar-L2}
    \|u_{\theta_t}(\cdot)-u^*(\cdot,t)\|_{2}\leq e^{(L_f+B/2-\lambda/C_p)t}(\varepsilon_{0} + {\varepsilon} t)
\end{equation}
for all $t$ as long as $\theta_{t} \in \barTheta$, where $\theta_t$ is solved from the ODE \eqref{eq:ode} with $V_{\xi}$ and initial $\theta_0$ satisfying $\|u_{\theta_0}(\cdot)-u^*(\cdot,0)\|_{2}\leq \varepsilon_0$. Here $C_p$ is a constant depending only on $\Omega$.
\label{thm:ode-sol-approx-error}
\end{theorem}
\begin{proof}
We denote the residual
\[
r(x,t):=\nabla_{\theta}u_{\theta_t}(x) \cdot V_{\xi}(\theta_t)-F [u_{\theta_t}](x).
\]
Then by Proposition \ref{prop:F_exists}  we have  $\|r(\cdot,t)\|_{2} \le  \varepsilon$ for all $t$.
Furthermore, we denote
\[
\delta(x,t) := \utt(x) - u^{*}(x,t)
\]
for all $(x,t) \in \barOmega \times [0,T]$ and $D(t):= \|\delta(\cdot,t)\|_{2}$, then there is
\begin{equation}
\label{eq:dD}
    D'(t) = \Big\langle \frac{\delta(\cdot,t)}{\|\delta(\cdot,t)\|_2}, \partial_{t}\delta(\cdot,t) \Big\rangle.
\end{equation}
Here we use the convention that $\delta(\cdot,t)/\|\delta(\cdot,t)\|_{2}=0$ if $\delta(\cdot,t) = 0$ a.e.
By the definition of $\delta$, we have
\begin{align*}
\partial_{t} \delta(x,t)
& = \partial_{t} \utt(x) - \partial_{t} u^{*}(x,t) \\
& = \nabla_{\theta} \utt(x) \cdot \dot{\theta}_{t} - F[u^{*}](x,t) \\
& = \nabla_{\theta} \utt(x) \cdot V_{\xi}(\theta_{t}) - F[u^{*}](x,t) \\
& = F[\utt](x) - F[u^{*}](x,t) + r(x,t) \\
& = \nabla \cdot (A(x) \nabla \delta(x,t)) + b(x) \cdot \nabla \delta(x,t) + f(\utt(x)) - f(u^{*}(x,t)) + r(x,t).
\end{align*}
Therefore, we have
\begin{align}
    \langle \delta(\cdot,t), \partial_{t} \delta(\cdot,t) \rangle 
    &=\int_{\Omega} \delta(x,t) \left(\nabla \cdot (A(x) \nabla \delta(x,t)) + b(x) \cdot \nabla \delta(x,t)\right)\,dx \nonumber \\
    & \qquad + \int_{\Omega} \delta (x,t)(f(u_{\theta_t}(x))-f(u^*(x,t)) + r(x,t))\,dx \label{eq:IJ}\\
    & =: I(t) + J(t). \nonumber 
\end{align}
Because $\utt(\cdot)|_{\partial\Omega} = u^{*}(\cdot,t)|_{\partial\Omega} = 0$, we know $\delta(\cdot,t)|_{\partial\Omega} = 0$. Thus, we have
\begin{align}
I(t) & = \int_{\Omega}\delta(x,t)\left(\nabla \cdot (A(x) \nabla \delta(x,t)) + b(x) \cdot \nabla \delta(x,t)\right)\,dx \nonumber \\
&=-\int_{\Omega}\nabla \delta(x,t)^{\top} A(x) \nabla \delta(x,t) \,dx - \frac{1}{2}\int_{\Omega}(\nabla \cdot b(x))  \delta(x,t)^2 \, dx \label{eq:I}\\
&\leq - \lambda \int_{\Omega}|\nabla \delta(x,t)|^2 \,dx - \frac{1}{2}\int_{\Omega}(\nabla \cdot b(x))  \delta(x,t)^2 \, dx \nonumber \\
&\leq -\frac{\lambda}{C_p}\int_{\Omega}|\delta(x,t)|^2 \,dx+\frac{B}{2}\int_{\Omega}|\delta(x,t)|^2 \,dx,\nonumber 
\end{align}
where the first equality is just by the definition of $I(t)$, the second equality is obtained by integrating by parts on both terms and using $\delta(\cdot,t)|_{\partial \Omega}=0$, the first inequality is due to \eqref{eq:A-bound}, and the last inequality is due to the Poincare's inequality 
\[
\|\delta(\cdot,t)\|_{2} \le C_p \|\nabla \delta(\cdot,t)\|_2
\]
as $\delta(\cdot,t)|_{\partial\Omega}=0$ for all $t$ (here $C_p$ the Poincare's constant depending on $\Omega$ only) and the bound \eqref{eq:b-bound}.
We can also obtain
\begin{align}
    J(t) 
    & = \int_{\Omega}\delta(x,t)(f(u_{\theta_t}(x))-f(u^*(x,t))-r(x,t))\,dx \nonumber \\
    & \leq \int_{\Omega}|\delta(x,t)| \cdot |f(u_{\theta_t}(x))-f(u^*(x,t))-r(x,t)| \,dx \nonumber \\
    &\leq \int_{\Omega}|\delta(x,t)| \cdot (L_f|\delta(x,t)|+|r(x,t)|) \,dx \label{eq:J}\\
    &\leq L_f\|\delta(x,t)\|_2^2+\|r(\cdot,t)\|_{2}\|\delta(\cdot,t)\|_{2} \nonumber \\
    &\leq L_f\|\delta(x,t)\|_2^2+  \varepsilon \|\delta (\cdot,t)\|_{2}, \nonumber 
    %&\leq L_f\|\delta(x,t)\|_2^2+ \varepsilon \|\delta (\cdot,t)\|_{2}, \nonumber 
\end{align}
where the first identity is by the definition of $J(t)$, the second inequality is due to the Lipschitz condition of $f$. 
Combining \eqref{eq:dD}, \eqref{eq:IJ}, \eqref{eq:I} and \eqref{eq:J}, we obtain
\[
D'(t) \leq  \Big(L_f+\frac{B}{2}-\frac{\lambda}{C_p}\Big) D(t) +  {\varepsilon}.
\] 
By Gr\"{o}nwall's inequality we deduce that
\begin{equation*}
    D(t)\leq e^{(L_f+B/2-\lambda/C_p)t}(D(0) + {\varepsilon} t).
\end{equation*}
Recall that
\[
D(0)=\|\delta(\cdot,0)\|_{2}= \|u_{\theta_{0}}(\cdot) - u^{*}(\cdot,0)\|_2=\|u_{\theta_{0}}(\cdot)-g(\cdot)\|_{2}\leq \varepsilon_{0},
\] 
we thus have
\[
\|u(\cdot,\theta(t))-u(\cdot,t)\|_{2}=D(t)\leq e^{(L_f+B/2-\lambda/C_p)t}(\varepsilon_{0} + {\varepsilon} t )
\]
for all time $t$, which completes the proof.
\end{proof}

The error estimate in Theorem 3.5 indicates that the approximation error is determined by three factors: (i) the approximation error $\varepsilon_0$ of the reduced order model to the initial value $g$, (ii) the local approximation error $\varepsilon$ of the projection of $F[u_\theta]$ onto the tangent space of $\mathcal{M}$ at $u_{\theta}$; and (iii) the irregularity of the differential operator $F$ itself. While the error from (iii) is determined by the given PDE, we can make an effort to suppress (i) and (ii) in practice by robust architecture of $u_{\theta}$ and the training of $V_{\xi}$. We note the error estimate provided in Theorem 3.5 is an upper bound of the approximation error.

\begin{remark}
While we assumed $f$ to be globally Lipschitz, the result in Theorem \ref{thm:ode-sol-approx-error} still holds locally with local Lipschitz condition of $f$. 
%for some trajectory $u^*(x,t)$, then the error bounds of Theorem \ref{thm:ode-sol-approx-error} still apply for trajectories $\theta_t$ that follow these particular trajectories using the local Lipschitz constant. 
For example, in the case of the Allen-Cahn example, we know if our initial function is bounded by 1 the true trajectories will remain bounded allowing the results of Theorem \ref{thm:ode-sol-approx-error} to apply.
\end{remark}

\begin{corollary}
Suppose the conditions in Theorem \ref{thm:ode-sol-approx-error} hold. Let $\hatthetat$ be the numerical solution to the ODE \eqref{eq:ode} obtained by using the Euler scheme with step size $h>0$. 
Then
\begin{equation}
\label{eq:uhattheta-ustar-L2}
\|u_{\hatthetat}(\cdot)-u^*(\cdot,t)\|_{2} \leq \frac{L_{V}M_{V}|\Omega|h}{2}(e^{L_{V}t} - 1) + e^{(L_f-B/2+\eta/C_p)t}(\varepsilon_{0} + \bar{\varepsilon} t)   
\end{equation}
for all $t$ as long as $\theta_{t} \in \barTheta$.
\label{cor:ode-numerical-approx-error}
\end{corollary}

\begin{proof}
Given the estimate provided in Theorem \ref{thm:ode-sol-approx-error}, we only need to show
\begin{equation}
\label{eq:uhatttheta-utheta-L2}
    \| u_{\hatthetat}(\cdot) - \utt(\cdot) \|_{2} \le \frac{L_{V}M_{V}|\Omega|h}{2}(e^{L_{V}t} - 1),
\end{equation}
since combined with \eqref{eq:utheta-ustar-L2} it yields the claimed estimate \eqref{eq:uhattheta-ustar-L2}.
To show \eqref{eq:uhatttheta-utheta-L2}, we notice that 
\[
\ddot{\theta}_t = \frac{d}{dt}V_{\xi}(\theta_{t}) =\nabla_{\theta}V_{\xi}(\theta_t) \cdot \dot{\theta}_{t} = \nabla_{\theta}V_{\xi}(\theta_t) \cdot V_{\xi}(\theta_t).
\] 
Therefore we have 
\[
|\ddot{\theta}_t| = |\nabla_{\theta}V_{\xi}(\theta_t) \cdot V_{\xi}(\theta_t)| \leq L_{V} M_{V}
\]
where $L_{V}$ and $M_{V}$ are defined in \eqref{eq:MV}.
Hence, by the standard results for the Euler's method \cite[pp.~346]{atkinson1989introduction}), we know the numerical solution $\hatthetat$ satisfies 
\begin{equation}
\label{eq:theta-ode-error}
|\hatthetat-\theta_t| \leq \frac{h M_V}{2} \left(e^{L_V t}-1\right)    
\end{equation}
for all $t$. 
Therefore, we obtain
\begin{align*}
    \| u_{\hatthetat} - \utt \|_{2} 
    & = \Big( \int_{\Omega} |u_{\hatthetat}(x) - \utt(x)|^2\,dx \Big)^{1/2} = \Big( \int_{\Omega} |\nabla_{\theta} u_{\tildethetat} (x) \cdot (\hatthetat - \theta_{t})|^2\,dx \Big)^{1/2} \\
    & \le L_{V} |\Omega| |\hatthetat - \theta_{t}| \le \frac{L_{V}M_{V}|\Omega|h}{2}(e^{L_{V}t} - 1),
\end{align*}
where the second equality is due to the fact that $\ut$ is $C^{1}$ in $\theta$ and hence the mean value theorem applies to $\ut$ (here $\tildethetat$ is some point on the line segment between $\hatthetat$ and $\theta_{t}$).
% \[
% \begin{aligned}
% \|u_{\theta_t}(\cdot)-u^*(\cdot,t)\|_{\infty}& \leq \|u_{\theta_t}(\cdot)-u_{\theta^*_t}(\cdot)\|_{\infty}+\|u_{\theta^*_t}(\cdot)-u^*(\cdot,t)\|_{\infty}\\
% & \leq  M_{u}|\theta_t-\theta^*_t|+C_{F}\bar{\varepsilon}+\varepsilon_{0}\\
% & \leq \left(\frac{C_{V}M_{u}}{2}\right)\left(e^{M_V t}-1\right) h + C_{F}\bar{\varepsilon}+\varepsilon_{0}.
% \end{aligned}
% \] Define $ c_1(t) = (\frac{C_{V} M_{u}}{2})(e^{M_V t}-1)$ and we are done.
\end{proof}
The proof above can be modified if a different numerical ODE solver is employed. In that case one can obtain improved upper bound and order in step size $h$ in \eqref{eq:theta-ode-error}.

\section{Numerical Results}
\label{sec:numerical-results}

\subsection{Implementation of the training process of control field $V_{\xi}$}
\label{subsec:training}
In Section \ref{subsec:method}, we have showed that the neural control field $V_{\xi}$ is parameterized as a deep network, and its parameters $\xi$ can be learned by solving
\begin{align*}
\min_{\xi} \cbr[2]{\ell(\xi) :=\int_{\Theta} \| V_{\xi}(\theta) \cdot \nabla_{\theta} \ut - F[\ut] \|_{2}^2 \, d \theta }.
\end{align*}
The first-order optimality condition of this minimization problem is given by $G(\theta)V_{\xi}(\theta)=p(\theta)$ where $G(\theta)$ and $p(\theta)$ are defined in \eqref{eq:def-G}.
The objective function $\ell(\xi)$ above shares the same minimizers as the following one:
\begin{equation}
 \bar{\ell}(\xi):=\int_{\Theta}|G(\theta)V_{\xi}(\theta)-p(\theta)|^2 \, d\theta.
\label{eq:proj_loss}
\end{equation}
In our numerical experiments, we use $\bar{\ell}$ defined in \eqref{eq:proj_loss}, as we can train to towards the optimal solution $V_{\xi}=G^{+}(\theta)p(\theta)$ as the optimal value which seems to produce lower error empirically. 
% \textcolor{blue}{However, where the dimension of $\theta$ is large, using \eqref{eq:proj_loss} may be infeasible.  In such cases, it is best to return to \eqref{eq:total_loss} so that $G(\theta)$ is not needed at any point during training or sample generation.}  
%\ye{I feel it's better to stay focused on $\ell$ and move $\bar{\ell}$ to experiment and explain there that appears to be more effective in training empirically.} 
Moreover, we know the minimum loss value of \eqref{eq:proj_loss} is $0$, which contrasts to \eqref{eq:total_loss} where the minimum loss value is often unknown.

%This is in contrast to the case for \eqref{eq:total_loss} where the minimum loss value is less likely to be 0, which makes it trickier to know how close to the minimum you are just by looking at the loss 
% If so we can guarantee that for any trajectory $\theta(t)$ that remains in $\Theta$ we have a good approximate solution. 

In practice, as the dimension of $\theta$ and $\Omega$ could be large, we have to approximate \eqref{eq:proj_loss} using techniques such as Monte-Carlo integration. This leads to the approximate forms
\[
\tilde{G}(\theta) = \frac{1}{N_x} \sum_{i=1}^{N_x} \nabla_{\theta}\ut(x_i)\nabla_{\theta}\ut(x_i)^{\top}, \quad \tilde{p}(\theta) =\frac{1}{N_x} \sum_{i=1}^{N_x}\nabla_{\theta}\ut (x_i)F [\ut](x_i),
\]
where $x_i$, $i=1,\ldots, N_x$ are sampled from $\Omega$. 
%
% As discussed in \cite{du2021evolutional}, we know that $N_x \to \infty$ gives us $\tilde{A}(\theta) \to A(\theta)$ and $\tilde{b}(\theta) \to b(\theta)$. 
%
By also drawing samples from $\Theta$, we arrive at our empirical loss function defined by
\begin{equation}
\ell_1(\xi):=\frac{1}{N_{\theta}}\sum_{j=1}^{N_{\theta}} |\tilde{G}(\theta_j) \cdot V_{\xi}(\theta_j)-\tilde{p}(\theta_j)|^2.
\label{eq:empirical_loss}
\end{equation}

To improve the training of $V_{\xi}$, we also augment the loss function $\ell_1$ in \eqref{eq:empirical_loss} with an additional term following a data-driven approach. Specifically, we follow the methods in \cite{du2021evolutional, bruna2022neural-galerkin} to generate multiple sample trajectories starting from randomly sampled initial values $\{\theta_{0}^{(i)}:i\in [M]\}$ in $\Theta$. For the $i$th trajectory, a sequence of directions $\{v_{j}^{(i)}: j=0,1\dots,N_{t}\}$ are solved from linear systems $\tilde{G}(\theta_{j}^{(i)}) v_{j}^{(i)} = \tilde{p}(\theta_{j}^{(i)})$ and the discrete-time points on the trajectory are obtained by $\theta_{j+1}^{(i)} = \theta_{j}^{(i)} + h v_{j}^{(i)}$ for $j=0,1,\dots,N_t-1$. We add the augment loss term
\begin{equation}
\ell_2(\xi):=\frac{1}{N_tM}\sum_{i=1}^M\sum_{j=1}^{N_t}|V_{\xi}(\theta_{j}^{(i)})-v_{j}^{(i)}|^2
\label{eq:data-driven}.
\end{equation}
%The use of \eqref{eq:data-driven} is to add further time horizon information to $V_{\xi}$ and improve the longer range prediction power. 
Combining with \eqref{eq:empirical_loss}, we obtain our final loss function
\begin{equation}
    \ell_{\text{total}}(\xi)=\ell_1(\xi)+\zeta \ell_2(\xi),
\label{eq:final-loss}
\end{equation}
where $\zeta$ is a weight parameter. In our experience for parabolic linear PDEs using only $\ell_1$ is sufficient to generate a good result. 
For the nonlinear case adding $\ell_2$ substantially improves training results empirically as network parameters may move far away from those we sampled near the initial parameters.

\subsection{Experimental setting}
To demonstrate the performance of the proposed method, we test it on three different PDEs: a 10-dimensional (10D) transport equation, a 10D heat equation, and a 2D Allen-Cahn equation. 
Both of the transport equation and heat equation are linear PDEs, while the Allen-Cahn is a highly nonlinear PDE.
In fact, we also tested 10D Allen-Cahn equation but only present the result of the 2D one here. This is because the true solution of Allen-Cahn equation does not have closed-form, and we have to employ a classical finite difference method, which does not scale to 10D case, to produce a reference solution for comparison. In contrast, we have closed-form solutions of the IVPs with transport and heat equations and hence we can use them as the true solution for direct comparison.
In our tests, we employ the following structure of our reduced-order model
\begin{align}
\label{eq:nn-zero-bdry}
\ut(x)= \alpha(x)z_{L}(x,\theta)
\end{align}
for the heat equation and Allen-Cahn equation. 
We use the following network structure
\begin{align}
\label{eq:nn-periodic-bdry}
\ut(x) =z_{L}(\beta(x),\theta)
\end{align}
for the transport equation.
In \eqref{eq:nn-zero-bdry}, $\alpha(x)$ is a distance function of $\partial \Omega$ such that it satisfies the zero boundary condition, and in \eqref{eq:nn-periodic-bdry} $\beta(x)$ is a function chosen to satisfy a periodic boundary condition as in \cite{du2021evolutional}. 
This aligns with our choice of $\ut$ in \eqref{eq:nn-zero-bdry} and \eqref{eq:nn-periodic-bdry} as the IVP with heat and Allen-Cahn equations have zero boundary value whereas the IVP with transport equation has periodic boundary value in our experiments.
In both \eqref{eq:nn-zero-bdry} and \eqref{eq:nn-periodic-bdry}, $z_{L}$ is the neural network and is defined by
\begin{equation}
z_L= w_L z_{L-1}, \quad z_{l}=z_{l-1}+\sigma(W_{l}z_{l-1}+b_{l}),\quad l=1,\ldots, L-1
\label{eq:neural-network}
\end{equation}
and $z_0 = \sigma(W_0x+b_0)$.
Here $\sigma$ is a user-chosen activation function (we use $\tanh$ or ReLU in our experiments) $W_{l} \in \Rbb^{d' \times d'}$ are the weight matrices and $b_{l} \in \Rbb^{d'}$ are the bias vectors, and $W_0 \in \Rbb^{d' \times d}$ and $w_L \in \Rbb^{1 \times d}$, all of these matrices and vectors make up the parameters vector $\theta$. Networks such as in \eqref{eq:neural-network} are often called \textit{residual neural networks} (ResNet), and have been shown performing better than basic feed forward networks in function approximation \cite{tabuada2022resnet}. 
The values of $L$ and $d'$ in our experiments are shown in Table \ref{tab:parameters}. 
%For the activation function, we shall use $\tanh$. 
%
They are selected manually to balance the depth $L$ and width $d'$ so that $\ut$ does not have too many neurons but still remains expressive. 
We use a similar structure for the vector field $V_{\xi}$, but adjust the layers to be $\eta_{l}=\eta_{l-1}+\text{GeLU}(\bar{U}_{l}\theta+\bar{b}_{l})\tanh(U_{l}\eta_{l-1}+b_{l})$. Here $\text{GeLU}(x)=x\Phi(x)$ where $\Phi(x)$ is the standard Gaussian cumulative distribution function. This is a slight modification of the network architecture proposed in \cite{srivastava2015highway} for improved effectiveness in training by gradient descent. 
We selected this network structure by starting with a ResNet with small width and depth and ReLU activation, then we gradually increased the width and depth until the improvement in the final loss value became insignificant. Finally, we attempted a few different activation functions and network architectures for this width and depth and selected the aforementioned structure which appeared to perform slightly better. This process was by no means exhaustive.
% We selected this network structure by training multiple architectures with different activations, starting with a small width and depth and slowly incrementing the width and depth up until the training loss no longer decreases. We then selected the architecture that yielded the lowest training loss. This does not constitute an exhaustive search of the possible architectures, but we found the given structure to perform well in our cases.

%
Other network architectures can be used as well. The width and depth of our network are reported in Table \ref{tab:parameters}.
Information about the number of trajectories used for \eqref{eq:data-driven} is also collected in Table \ref{tab:parameters}. For all of the experiments, we set the weight $\zeta=0.1$ in \eqref{eq:final-loss} to reflect the scale difference of the two loss terms and use the standard ADAM optimizer with learning rate 0.001, $\beta_1=0.9$, $\beta_2=0.999$.
%\ye{also provide values of parameters like learning rate, $\beta_1$, $\beta_2$.} to train $V_{\xi}$.
We terminate the training process when the empirical loss $\ell_{\text{total}}(\xi) < 0.1$ or when the percent decrease of the empirical loss is less than $0.1 \%$ averaged over the past 100 steps. Once $V_{\xi}$ is learned, we use the $4$th-order Runge-Kutta method with a step size of $T/200$ ($T$ is determined from the problem)
%\ye{step size?}
to solve $\theta_{t}$ from the ODE in \eqref{eq:ode} in Algorithm \ref{alg:solving-ivp} and compare the corresponding $\ut$ with the reference solutions. All the implementations and experiments are performed using PyTorch in Python 3.9 in Windows 10 OS on a desktop computer with an AMD Ryzen 7 3800X 8-Core Processor at 3.90 GHz, 16 GB of system memory, and an Nvidia GeForce RTX 2080 Super GPU with 8 GB of graphics memory. 
Total computational time is split between three unique activities: (i) the generation of $N_{\theta}$ samples for $\ell_1$ in \eqref{eq:empirical_loss}; (ii) the generation of the $M$ trajectories for $\ell_2$ in \eqref{eq:data-driven}; and (iii) the training of the network $V_{\xi}$. Parts (i) and (ii) can be parallelized offline to speed up the process. We discuss the specific time cost of the implementation of our method in the examples below.

We also provide a few remarks on the sampling strategy in $\Theta$. While one can draw $\theta$ uniformly from $\Theta$, adding samples $\theta$ corresponding to some example solutions to the PDE may further improve training efficiency. In practice, we use both uniformly sampled $\theta$'s and those close to the $\theta$'s corresponding to some randomly chosen initial functions. These initial functions are only used to help the loss function weigh more on the regions that are potentially more important than others in $\Theta$; but they are not among the randomly chosen initial functions used for any testing. 
%These samples may help the training to weigh more on the regions that are potentially more important than others in $\Theta$, thus making the training of $V_{\xi}$ more effective. 
Details on samplings are given below.

\begin{table}
\caption{Problem settings, network structures, and the number of training trajectories/samples in numerical experiments. Here $M$ is the number of trajectories used from $\Theta$ and $N_{\theta}$ is the total number of samples from $\Theta$.}
%\ye{Switch rows 1 and 2}
\begin{center}
\begin{tabular}{cccccc}
\toprule
\textbf{Problem} &\textbf{Dim.\ $d$} &$\ut$ \textbf{Width/Depth}&$V_{\xi}$ \textbf{Width/Depth} & \textbf{$M$} & \textbf{$N_{\theta}$} \\
\midrule
Transport Equation& 10 & 12/4 & 1,500/4& 0 & 160,000\\

Heat Equation &10 & 12/5& 2,000/10& 600& 200,000\\

Allen-Cahn Equation &2& 10/3 & 2,000/5& 200 & 200,000\\
\bottomrule
\end{tabular}
\end{center}
\label{tab:parameters}
\end{table}

% \begin{table}
% \caption{The time taken (hours) to complete each of the steps in the training of $V_{\xi}$ for each of the considered problems.}
% \begin{center}
% \begin{tabular}{ccccc}
% \toprule
% \textbf{Problem} &\textbf{Total Time} & \textbf{Time for (i)}& \textbf{Time for (ii)} & \textbf{Time for (iii)} \\
% \midrule
% Heat Equation &4.10 hrs & 2.64 hrs& 1.33 hrs& 0.13 hrs\\

% Transport Equation& 1.76 hrs& 1.68 hrs& 0 hrs& 0.08 hrs\\

% Allen-Cahn Equation &1.04 hrs& 0.25 hrs& 0.72 hrs & 0.07 hrs\\
% \bottomrule
% \end{tabular}
% \end{center}
% \label{tab:times}
% \end{table}

\subsection{Numerical results on transport equation}
We first consider the initial value problem defined by a 10D transport equation with periodic boundary conditions as follows:
\begin{align}
\label{eq:ivp-transport-eq}
\begin{cases}
\partial_t u(x,t)=-\mathbf{1} \cdot \nabla_x u(x,t), & \quad \forall\, x \in \Omega,\ t \in [0,T],\\
 u(x,0) = g(x), & \quad \forall\, x \in \barOmega,    
\end{cases}
\end{align} 
where $\Omega=(0,1)^{10}$, $T=1$, $\mathbf{1}$ is the vector whose components are all ones, and the boundary value $u(x,t)=0$ for all $x \in \partial \Omega$ and $t \in [0,T]$. This IVP has the true solution $u^{*}(x,t)=g(x-\mathbf{1}\cdot t)$. 
We obtain the solution operator of the IVP \eqref{eq:ivp-transport-eq}, we use \eqref{eq:nn-periodic-bdry} as the reduced-order model $\ut$. Although our error analysis requires certain regularity on initial and solution of PDEs, we test on the case where both are only Lipschitz continuous but not differentiable for this transport equation. To this end, we set the activation of $u_{\theta}$ to ReLU. Further, define $\beta(x)=(\cos(2\pi(x-b)),\sin(2\pi(x-b)))^{\top}$ where $b \in \Rbb^{10}$ is a trainable parameter with $\sin$ and $\cos$ acting component-wise to $x$. 
%\ye{what is $b$ inside cos and sin?}
This means that the first hidden layer uses $W_0\in\Rbb^{12 \times 20}$. For this example, we shall set $\Theta = [-1,1]^{m}$ where $m$ are the number of parameters in $\ut$.
%We set the structure of $V_{\xi}$ as \eqref{eq:neural-network} with ReLU activation. Here the depth is 2 and the width is a constant 32. %\ye{specify $V_{\xi}$ here} 
Then we train the neural control vector field $V_{\xi}$ by minimizing \eqref{eq:total_loss} with the number of sampled $\theta$ drawn uniformly from $\Theta$ shown in Table \ref{tab:parameters}. We note that this equation performed equally well with or without the loss $\ell_2$ in \eqref{eq:data-driven}. As such we need not generate any trajectories and this is reflected in Table \ref{tab:parameters}.

\begin{figure}
\centering
\includegraphics[width=.75\textwidth]{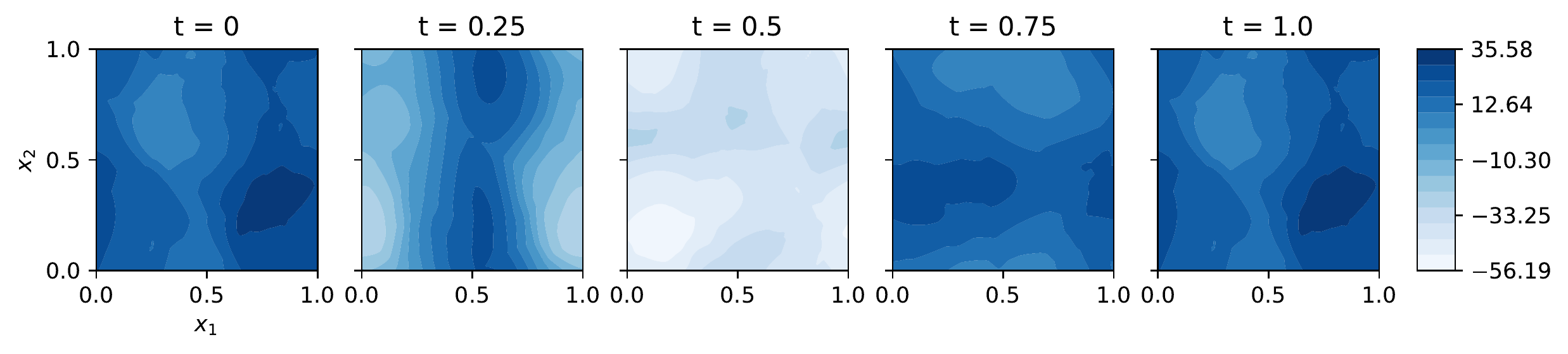}\vspace{-12pt}
\includegraphics[width=.75\textwidth]{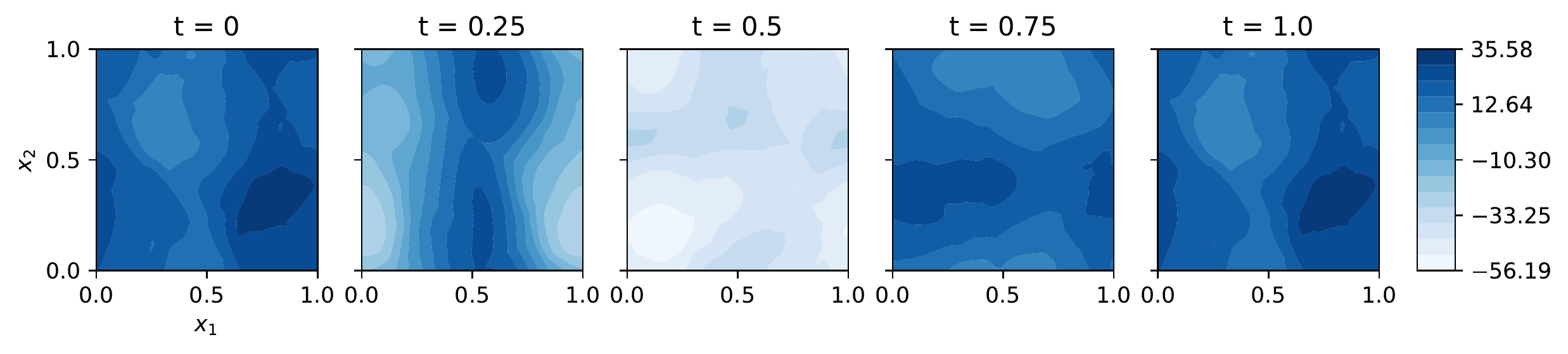}\vspace{-12pt}
\includegraphics[width=.75\textwidth]{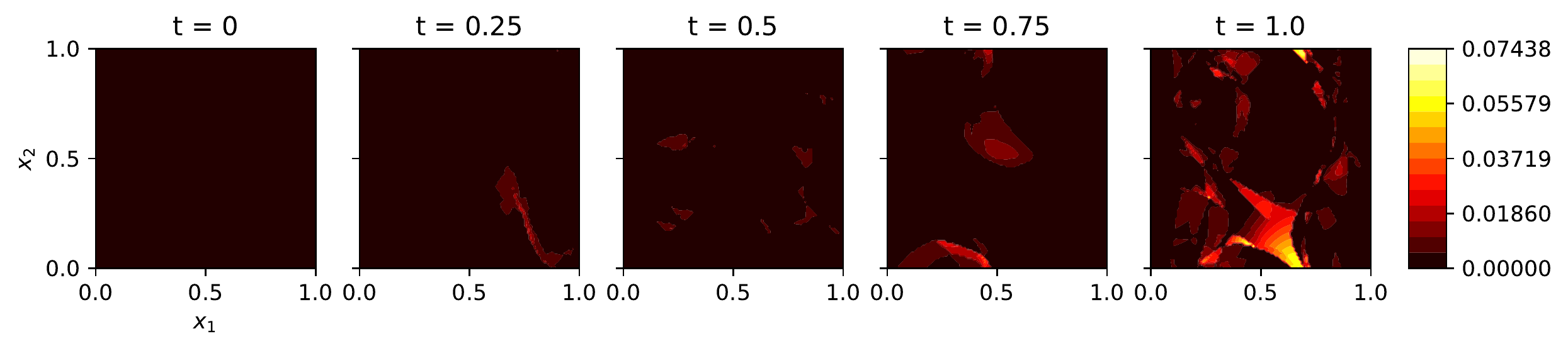}\vspace{-12pt}
\includegraphics[width=.75\textwidth]{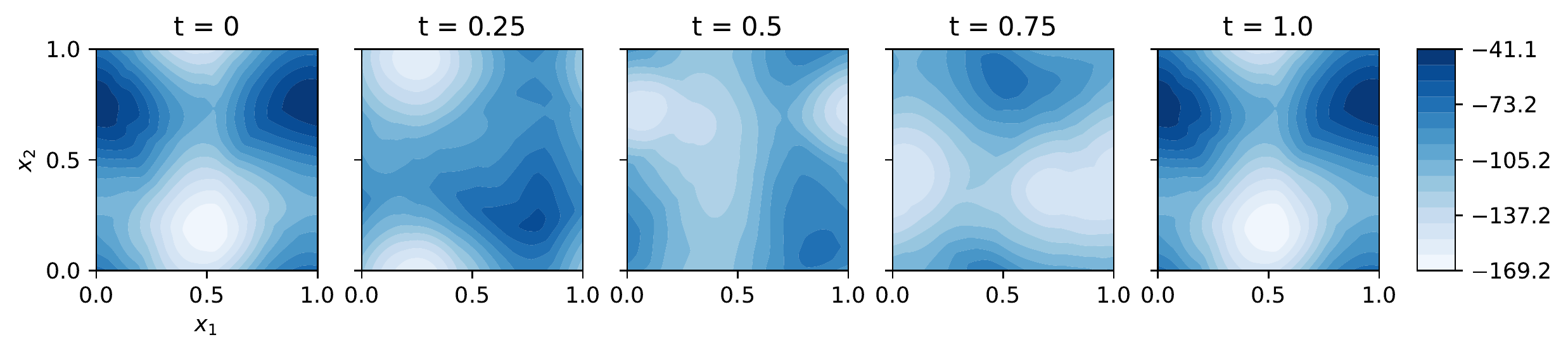}\vspace{-12pt}
\includegraphics[width=.75\textwidth]{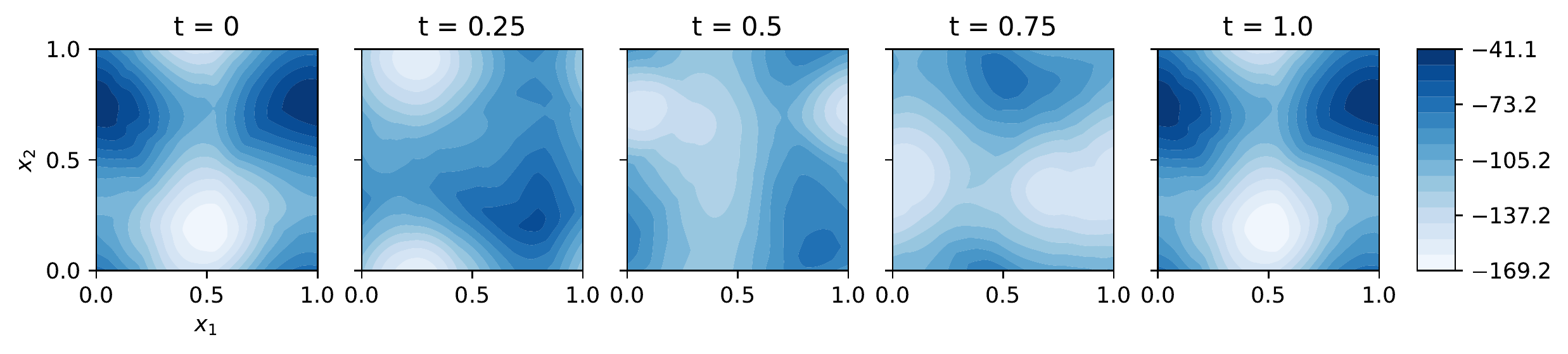}\vspace{-12pt}
\includegraphics[width=.75\textwidth]{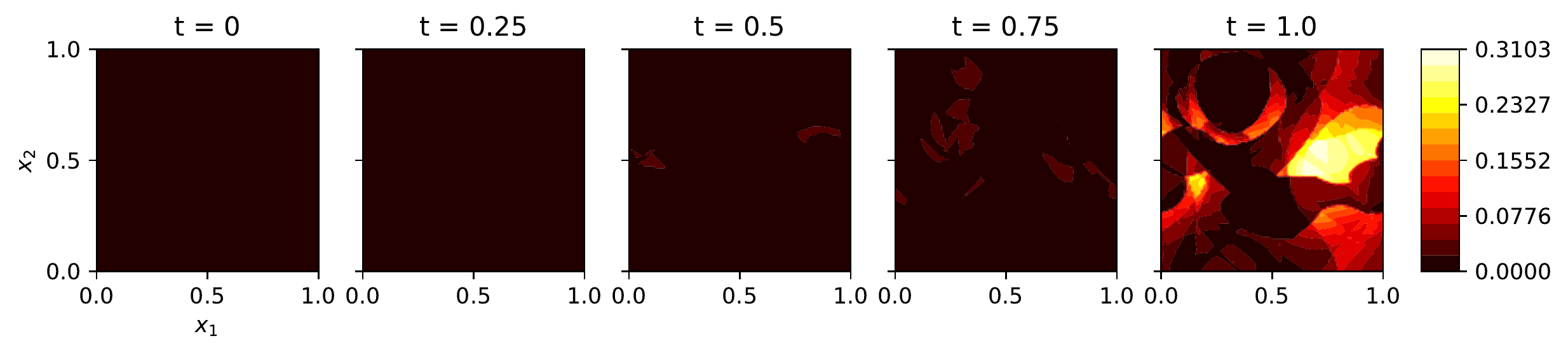}\vspace{-12pt}
\includegraphics[width=.75\textwidth]{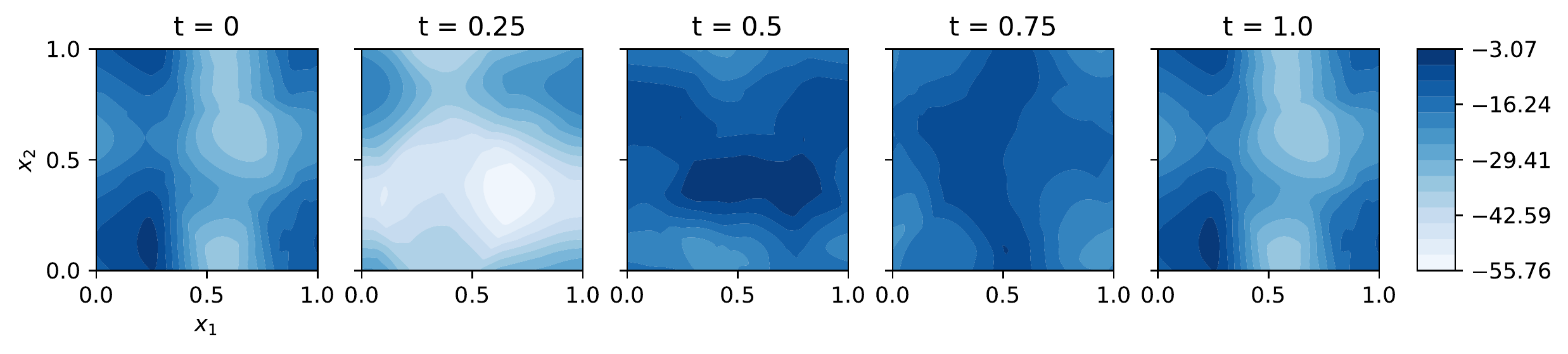}\vspace{-12pt}
\includegraphics[width=.75\textwidth]{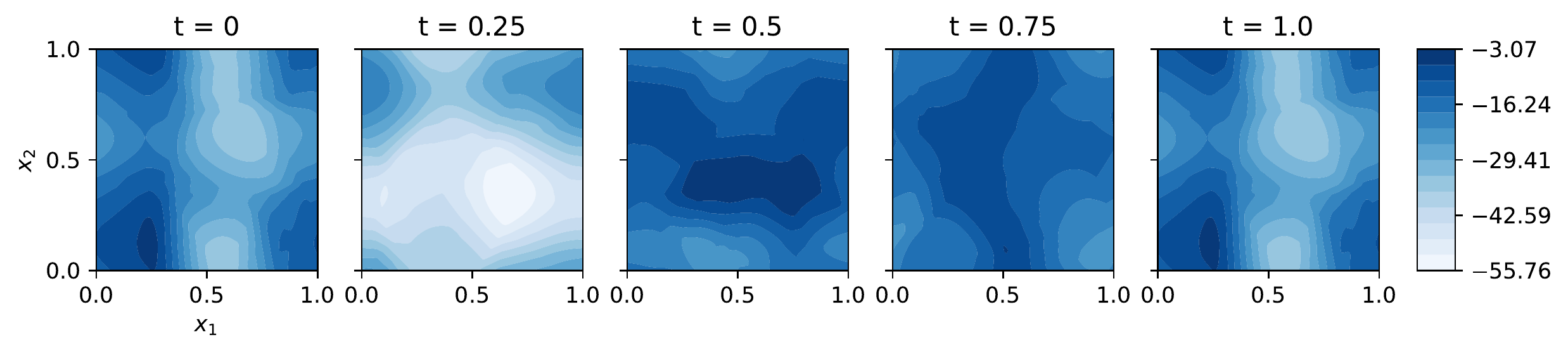}\vspace{-12pt}
\includegraphics[width=.75\textwidth]{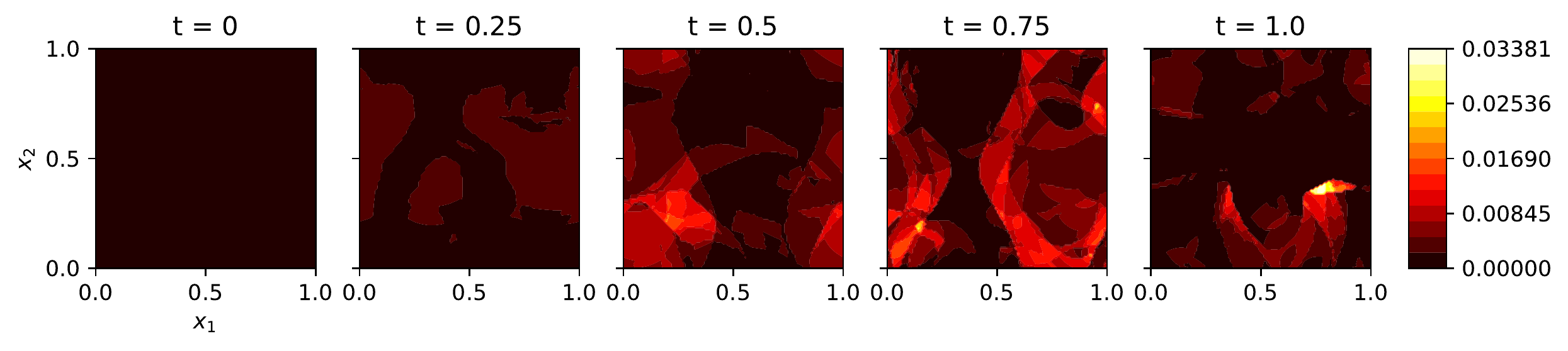}\vspace{-12pt}
\caption{(Transport equation). Comparison between true solution $u^{*}(\cdot,t)$, the approximation $\utt(\cdot)$ and their pointwise absolute difference $|\utt(x) - u^{*}(x,t)|$ for times $t=0,0.15,0.5,0.85,1$ for IVPs with the first initial (rows 1--3),  second initial (rows 4--6) and third initial (rows 7--9) given by $u_{\theta}$ with $\theta$ randomly drawn from $[-1,1]^m$.}
\label{fig:transport_example}
\end{figure}

\begin{figure}
    \begin{minipage}{0.3\textwidth}
    \centering
    \begin{subfigure}[b]{0.98\textwidth}
        \includegraphics[width=\textwidth]{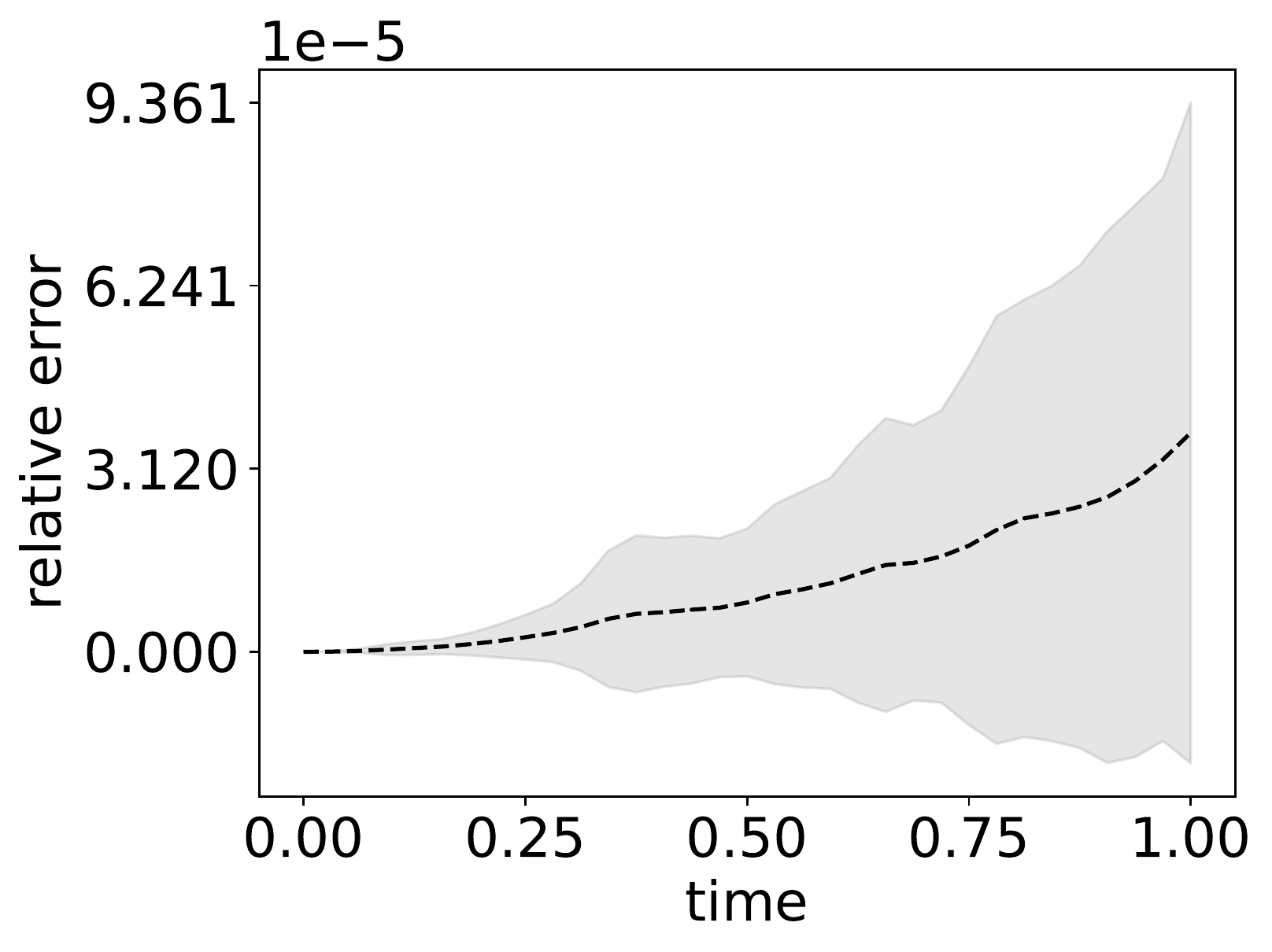}
        \caption{\label{fig:errors-a}}
    \end{subfigure}
    \begin{subfigure}[b]{\textwidth}
        \includegraphics[width=\textwidth]{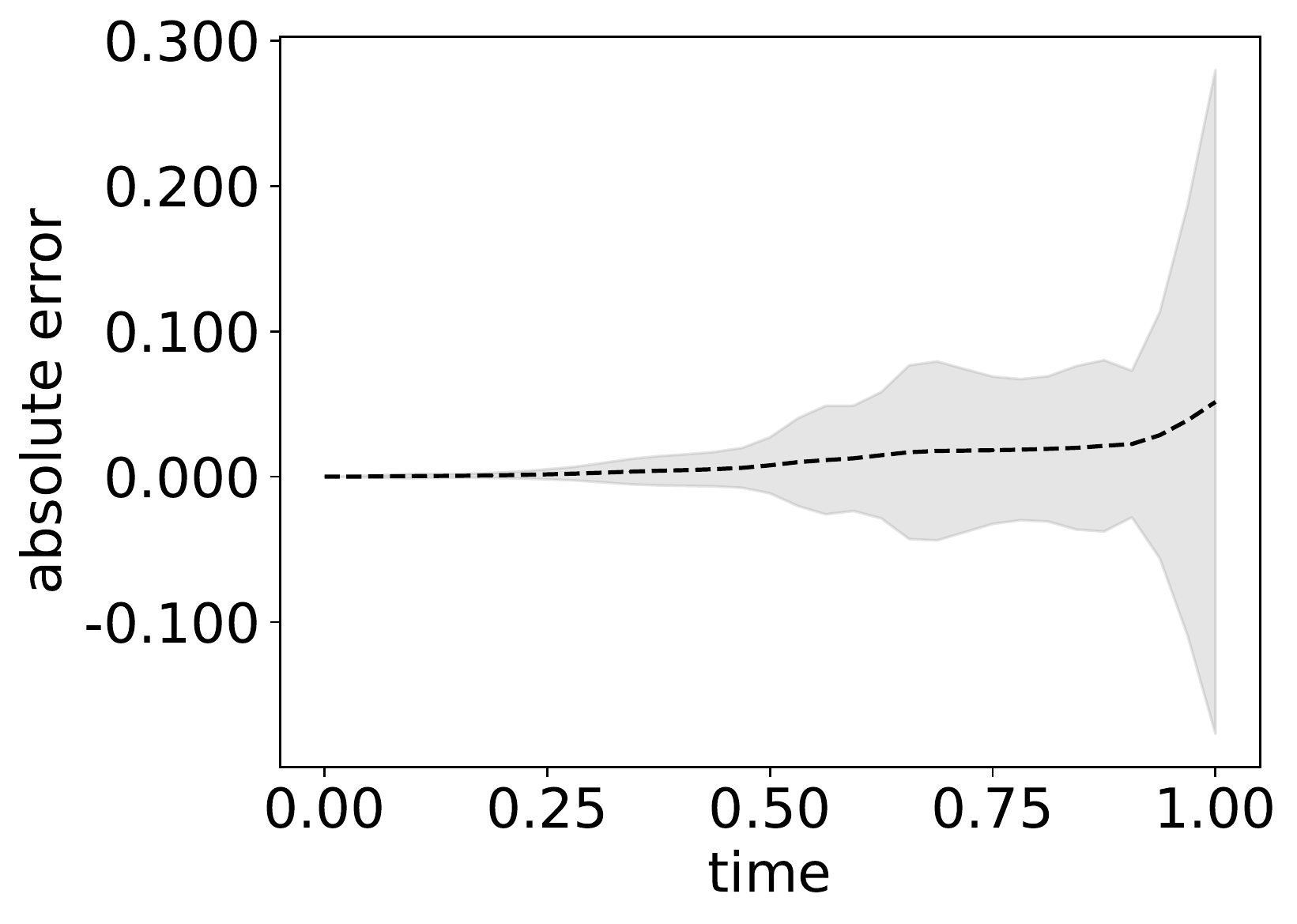}
        \caption{\label{fig:errors-b}}
    \end{subfigure}
    \vspace{-0.1cm}
    \end{minipage}
    \begin{minipage}{0.305\textwidth}
    \centering
    \begin{subfigure}[b]{\textwidth}
        \includegraphics[width=\textwidth]{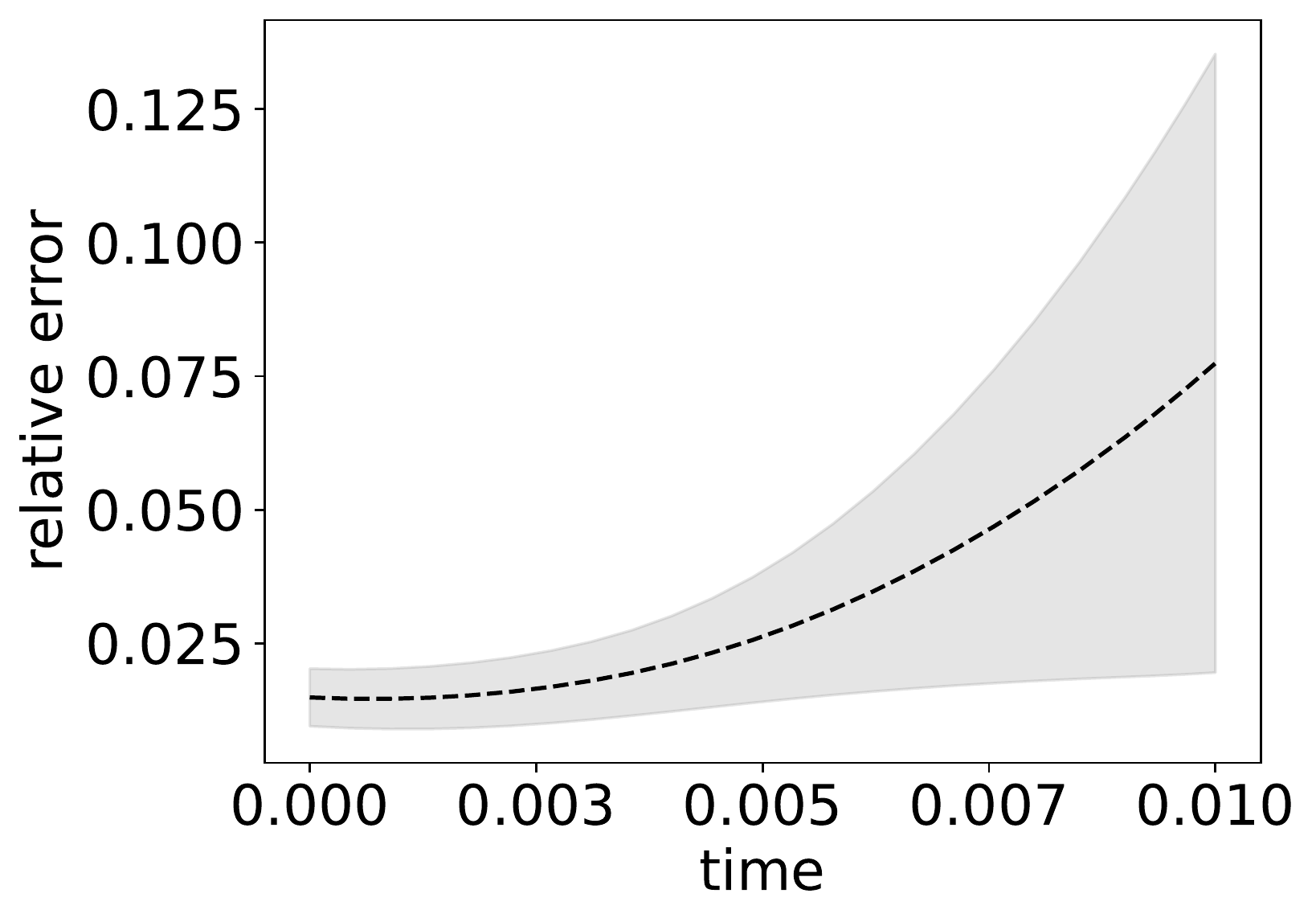}
        \caption{\label{fig:errors-c}}
        % \vspace{-0.2cm}
    \end{subfigure}
    \begin{subfigure}[b]{0.99\textwidth}
    % \vspace{-0.2cm}
        \includegraphics[width=\textwidth]{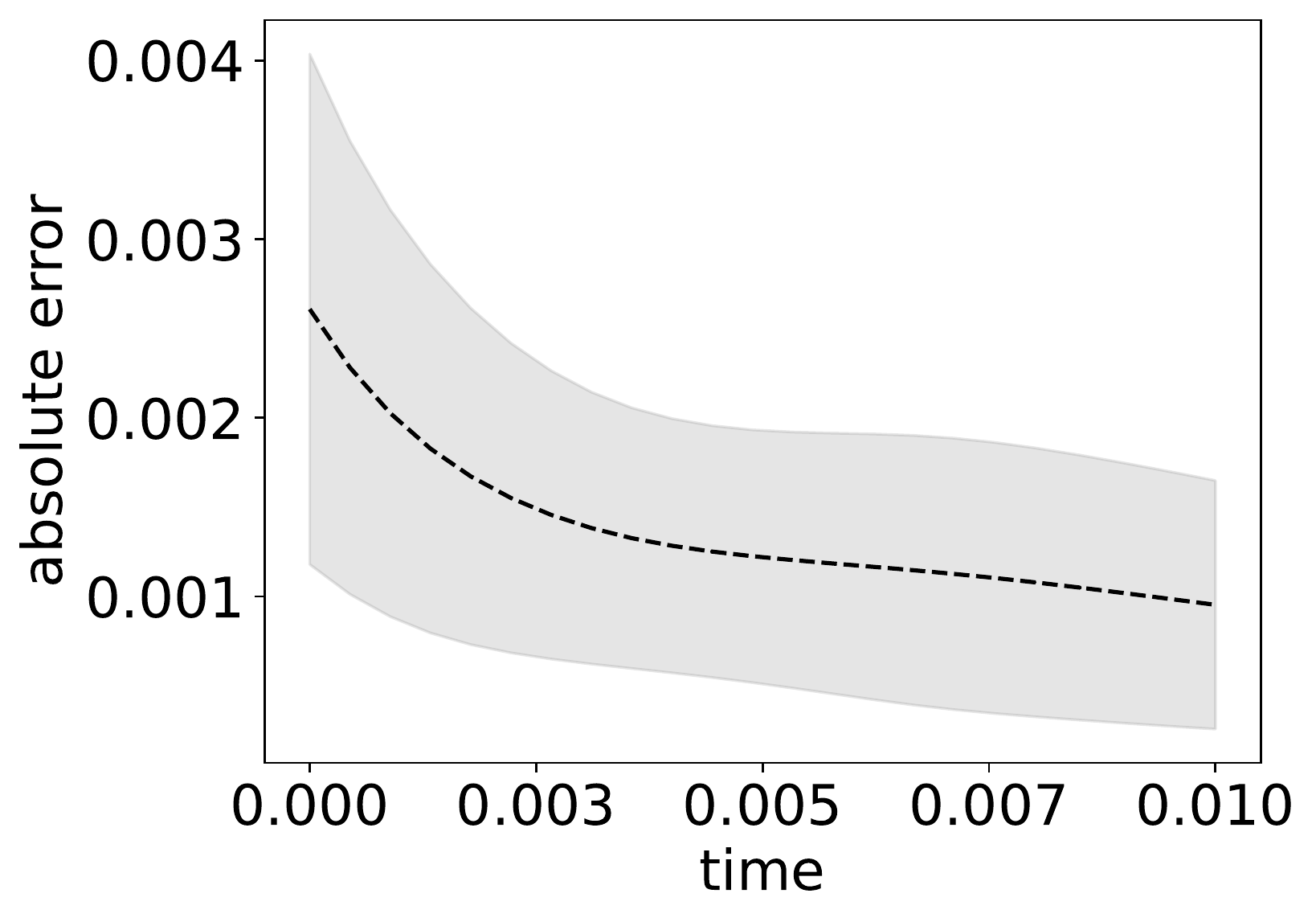}
        \caption{\label{fig:errors-d}}
    \end{subfigure}
    \end{minipage}
    \begin{minipage}{0.3\textwidth}
    \centering
    \begin{subfigure}[b]{0.99\textwidth}
        \includegraphics[width=\textwidth]{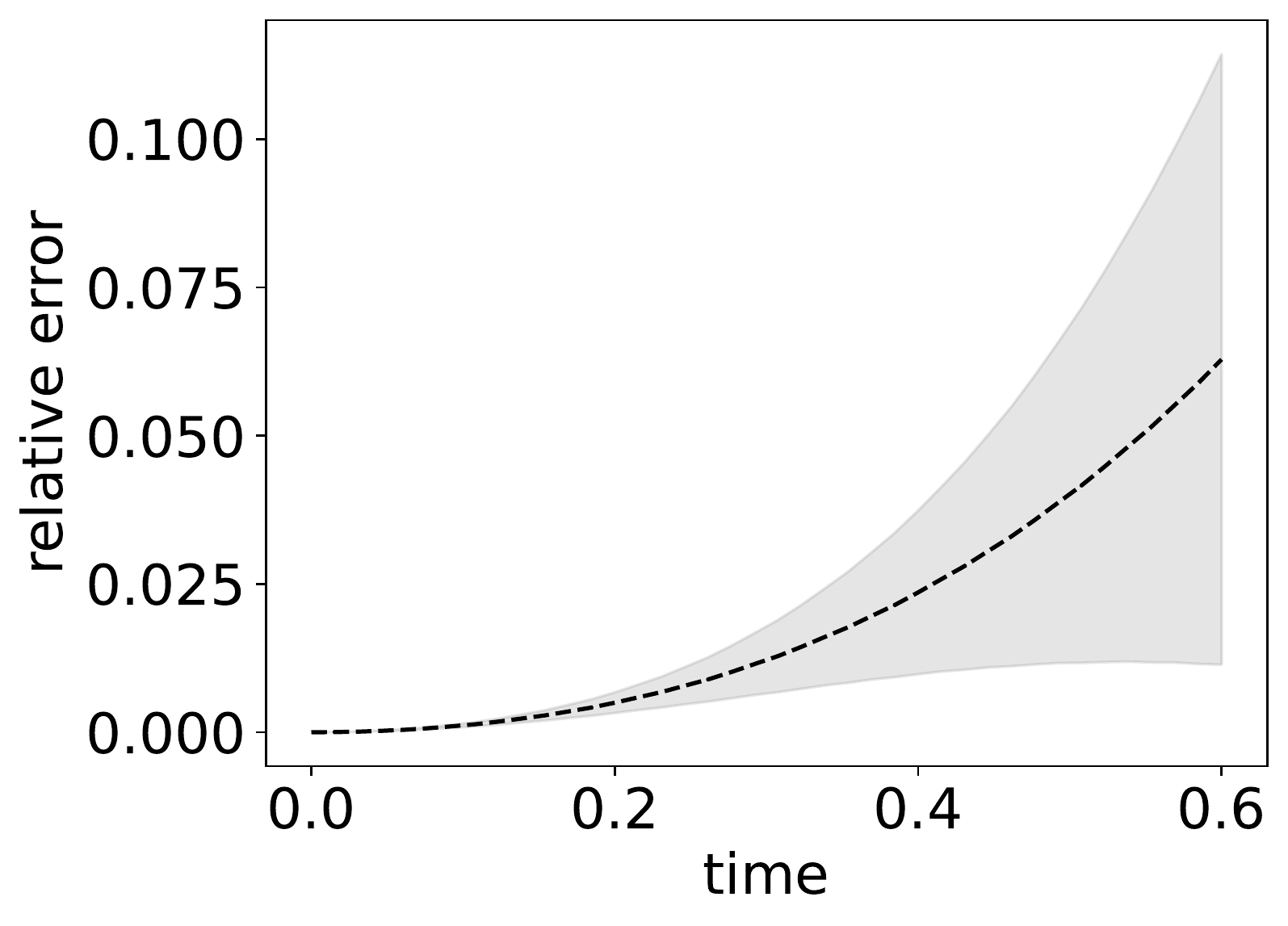}
        \caption{\label{fig:errors-e}}
    \end{subfigure}
    % \vspace{-0.2cm}
    \begin{subfigure}[b]{0.99\textwidth}
        \includegraphics[width=\textwidth]{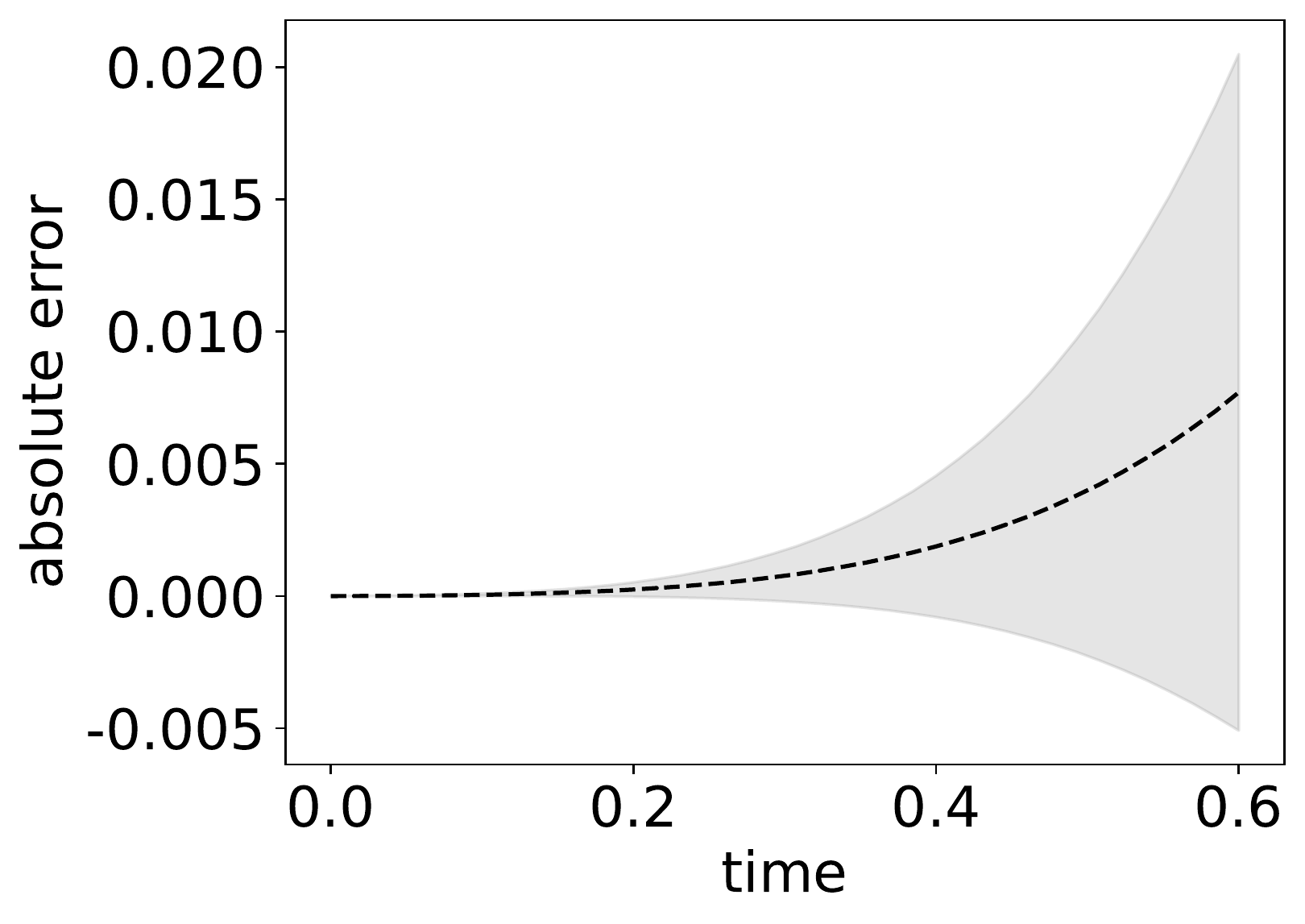}
        \caption{\label{fig:errors-f}}
    \end{subfigure}
    \end{minipage}
    \caption{Comparison of the mean relative error $\|u^{*}(\cdot,t)-\utt(\cdot)\|^2_2/\|u^{*}(\cdot,t)\|^2_2$ (top) and mean absolute $\|u^{*}(\cdot,t)-\utt(\cdot)\|^2_2$ (bottom) versus time $t$ for 100 different initial conditions of the transport (a)-(b), heat (c)-(d), and Allen-Cahn (e)-(f) equations. Shaded areas indicate the standard deviation over the 100 results.}
    %\ye{Usually (a) (b) etc are put under the plots.}}
    \label{fig:errors}
\end{figure}

After the control $V_{\xi}$ is obtained, we test the performance of $V_{\xi}$ on a variety of initial values $g$ by uniformly sampling $\theta_0 \in \Theta$.
We emphasize that the corresponding $\theta_0$'s to these initial values are not used in the training process. We show three approximate solutions for three random initials in Figure \ref{fig:transport_example}. For the first random initial $g_{1}$ determined by the random $\theta_0$, we plot the corresponding true solution $u^*(\cdot,t)$, the approximate solution $\utt(\cdot)$ obtained by Algorithm \ref{alg:solving-ivp}, and their pointwise absolute difference $|\utt(x) - u^*(x,t)|$ from row 1 to row 3 in Figure \ref{fig:transport_example} respectively for $t=0,0.15, 0.5, 0.85,1$. The plots for the second and third $g_2$ and $g_3$ random initials are shown in rows 4--6 and 7--9 in Figure \ref{fig:transport_example} respectively.
%\ye{If I remember correctly you also have plots to show $\|\ut - u(\cdot,t)\|_2$ vs $t$ for each initial? We should include these plots as the provide quantitative comparison.}
%
From Figure \ref{fig:transport_example}, we can see that the reduced-order model $\utt$ with $\theta_{t}$ controlled by the trained vector field $V_{\xi}$ closely approximates the true solution $u^*(\cdot,t)$ with low absolute errors (note that the scale of the error is different from that of $u^{*}(\cdot,t)$ and $\utt(\cdot)$). Figure \ref{fig:errors-a} and \ref{fig:errors-b} plots the mean of the absolute error $\|u^{*}(\cdot,t) -\utt(\cdot)\|^2_2$, and the relative error $\|u^{*}(\cdot,t) -\utt(\cdot)\|^2_2/\|u^{*}(\cdot,t)\|^2_2$ respectively over 100 randomly chosen initials, while the standard deviation is shaded in. We see mean errors $<1\%$ even though the initial functions considered are not smooth. This suggests that the proposed model can generalize to the case where the initial and solution of the PDE are not sufficiently smooth. 

We now discuss the computational cost of the method. In our tests, it took 1.78 hours to generate $\tilde{G}$ and $\tilde{p}$ from the samples in $\Theta$ used for training. 
%\ye{$\tilde{G}$ and $\tilde{p}$? If so better explicitly state it}\gaby{Yes stated}. 
Once generated, the training of $V_{\xi}$ (i.e., minimizing the loss function $\ell_{\text{total}}$ in \eqref{eq:final-loss}) took 5 minutes to complete. Testing each initial condition by solving \eqref{eq:ode} using a 4th-order Runge-Kutta (RK4) solver with step size $0.005$ took an average of 2.1 seconds per initial. We note that no time is needed in this case to fit an initial, as $\theta_0$ is chosen randomly.

The proposed method has evidient improvement on computational cost over existing methods that only solve specific instances of the PDEs. 
In this test, we compare the computational cost with PINN \cite{raissi2019physics-informed} and a time marching (TM) \cite{du2021evolutional} method. We use the same structure of $u_{\theta}$ for PINN and time marching as used by our method. We sample 10,000 points $(x,t) \in (0,1)^{10}\times[0,1]$ for PINN and 10,000 points $x \in (0,1)^{10}$ for each step of the time marching method. We train PINN using its default parameters until convergence. For TM, we use RK4 with the same step size $0.005$ and its default linear system solver for each step. 
%We report the times for these and our method in Table \ref{tab:pinn}. 
We follow all other implementation steps of both PINN and TM as described in their original papers. For a single initial $g$, PINN, TM, and the proposed method took 116.5s, 16.7s, and 2.1s respectively to obtain the solution. This significant time reduction is due to the fact that the proposed method has learned the control field in the parameter space and thus can compute the solution of the PDE by solving an ODE which has very low computation complexity. The improvement is more significant for higher-order PDEs because PINN and TM require more time to compute the differential operator whereas the computation complexity of the proposed method remains the same.
% \begin{table}[h]
%     \centering
%     \begin{tabular}{|c|c|}
%     \hline
%         Method & Time\\
%         \hline
%          PINN \cite{raissi2019physics-informed}&  116.5s\\
%          Time Marching \cite{du2021evolutional} &  16.7s\\
%          Neural Control& 2.1s \\
%          \hline
%     \end{tabular}
%     \caption{Computational time comparison of the Neural control method compared to standard PINN, and Time marching implementations using the $u_{\theta}$ network reported in Table \ref{tab:parameters}. Times are for a single initial $g$.}
%     \label{tab:pinn}
% \end{table}
% From Table \ref{tab:pinn} we see the favorable time benefit the method of Neural Control provides when one wishes to test multiple different $g$ initial functions for a particular PDE. We note as the transport equation is first-order the time to compute both PINN and time marching increases for the following each of the following examples due to their higher order.

The proposed method is capable of approximating solution operators of high-dimensional PDEs whereas existing methods cannot. This is because existing solution operator learning methods, such as DeepONet, require spatial discretization, and thus the network size and sampling amounts grow exponentially fast in problem dimension. For example, for a one-dimensional ($d=1$) evolution PDE, DeepONet \cite{lu2019deeponet} requires 100 sample solutions (which must be generated by another numerical method) each evaluated at $10^4$ grid points in the $(x,t)$ domain in $\Rbb\times \Rbb_{+}$. Thus the size of their trunk network alone is already 10 times larger than our $V_\xi$ in the 10-dimensional case. When the problem dimension $d$ becomes over 3, DeepONet will be infeasible computationally. In addition, our method does not require sample solutions which could be unavailable or difficult to obtain in practice.

% \begin{figure}
% \centering
% \includegraphics[scale=0.4]{pde_results/transport_errors.pdf}
% \caption{Error between the true and neural network solution for different initial conditions. To ensure these conditions were contained in the $\Theta$ region, the parameters were determined ahead of time, and $\Theta$ was designed to cover them. The points were not in the training data.}
% \label{fig:transport_errors}
% \end{figure}

\subsection{Heat equation}
Next we consider an initial value problem with heat equation in 10D:
\begin{equation}
\begin{cases}
\partial_t u(x,t) = \Delta u(x,t), & \quad \forall\, x \in \Omega, t \in [0,T]\\
% u(x,t)&=0, & \forall x \in \partial \Omega \\
u(x,0)=g(x), & \quad \forall\, x \in \barOmega,
\end{cases}
\label{eq:ivp-heat-eq}
\end{equation}
where $\Omega=(0,1)^{10}$ and the boundary value $u(x,t)=0$ for all $x \in \partial \Omega$ and $t \in [0,T]$. 
As most of the initial conditions we consider have rapid evolution in a short time, we use $T=0.01$ in this test. For neural network we use \eqref{eq:nn-zero-bdry}, with $\alpha (x)=\Pi_{i=1}^{10}4(x_i-x_i^2)$ and $\tanh$ activation. 
%We designed the region $\Theta$ to contain a few different initial conditions representing frequencies of sine functions. Once designed, we proceeded to train as discussed in Section \ref{subsec:training}. 

\begin{figure}[ht]
\centering
\includegraphics[width=.75\textwidth]{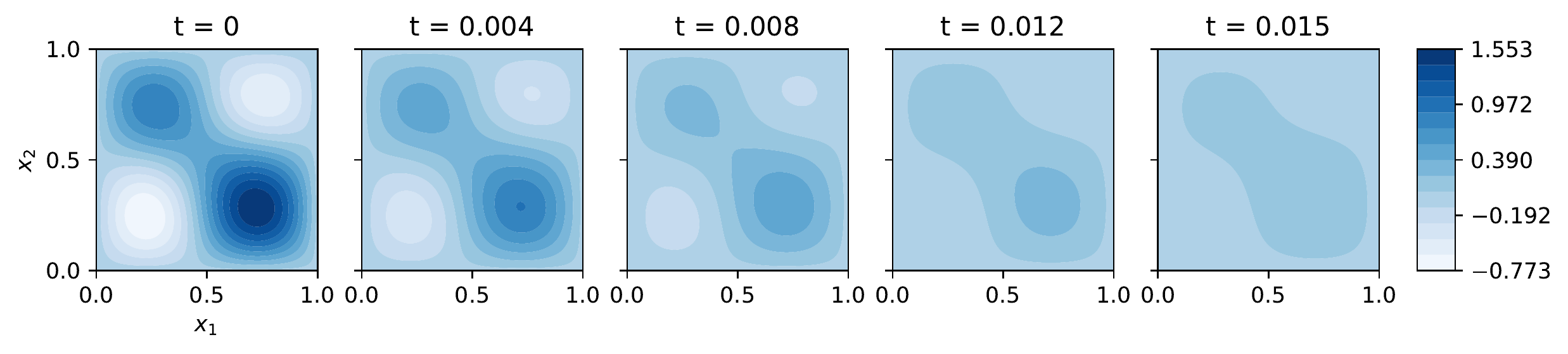}\vspace{-12pt}
\includegraphics[width=.75\textwidth]{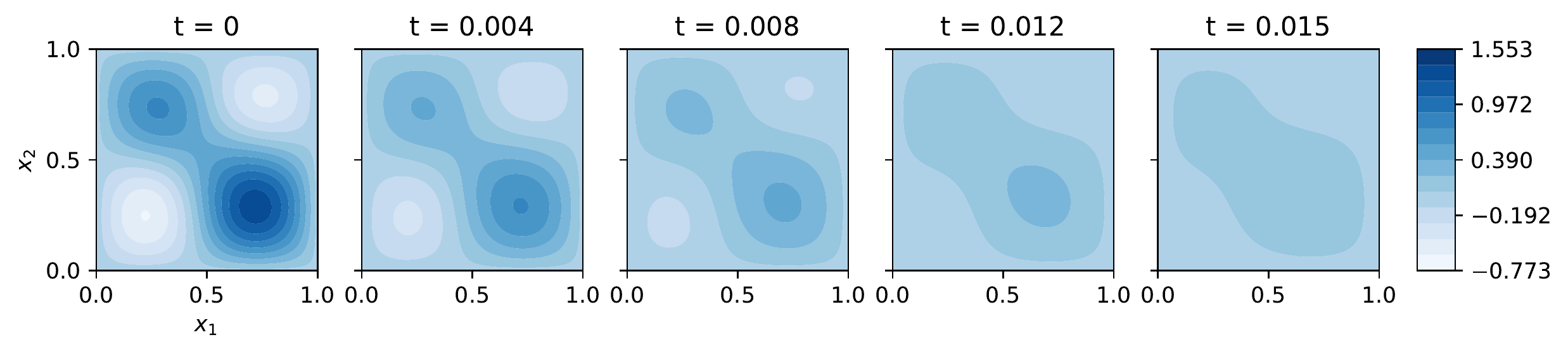}\vspace{-12pt}
\includegraphics[width=.75\textwidth]{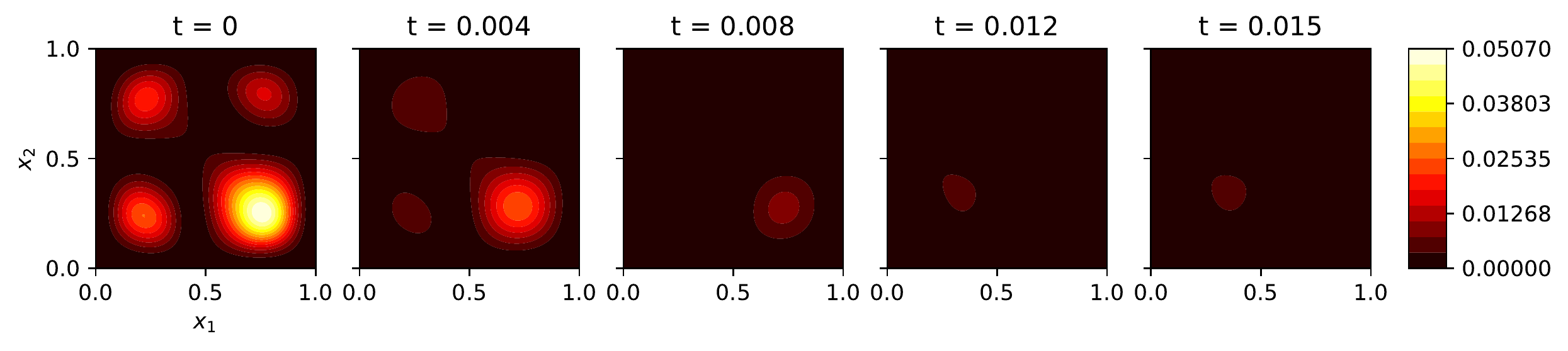}\vspace{-12pt}
\includegraphics[width=.75\textwidth]{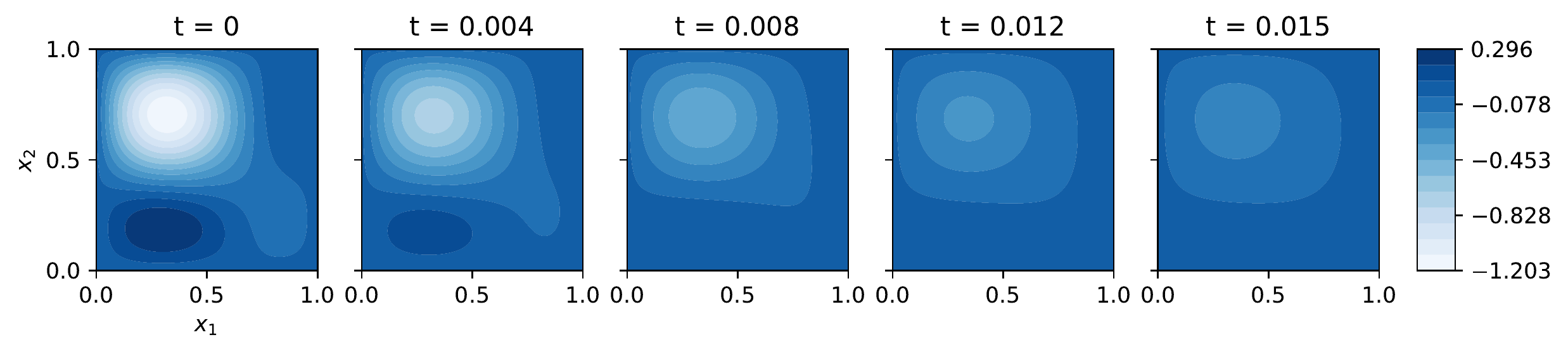}\vspace{-12pt}
\includegraphics[width=.75\textwidth]{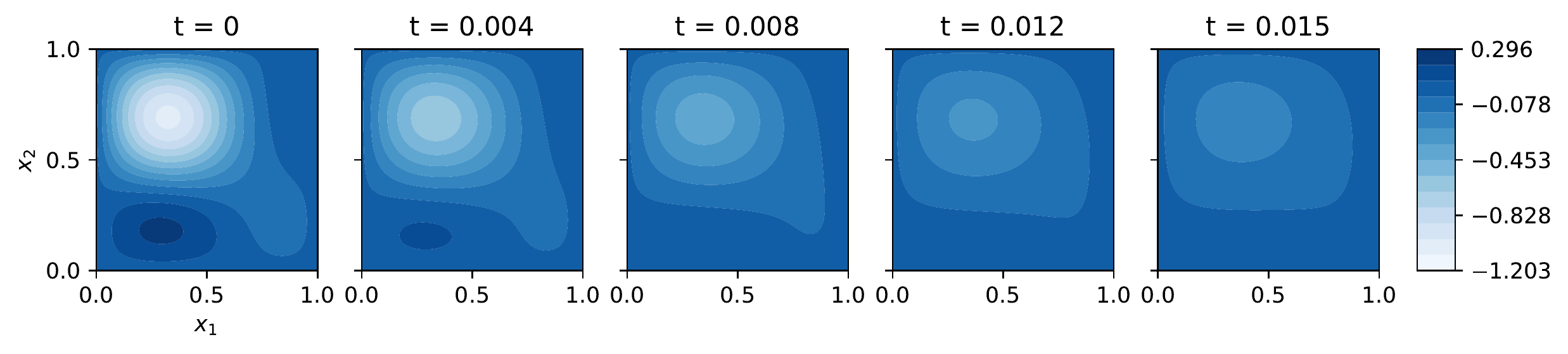}\vspace{-12pt}
\includegraphics[width=.75\textwidth]{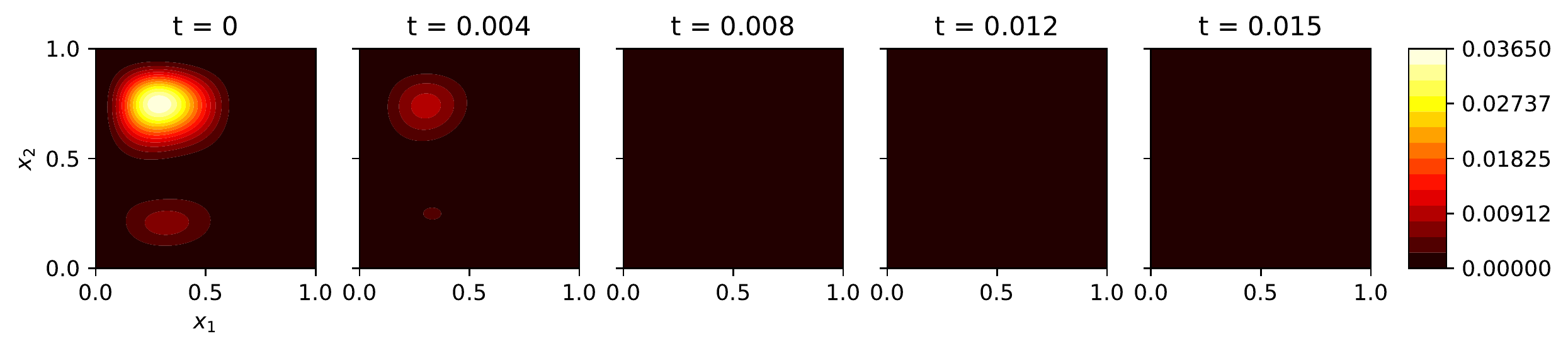}\vspace{-12pt}
\includegraphics[width=.75\textwidth]{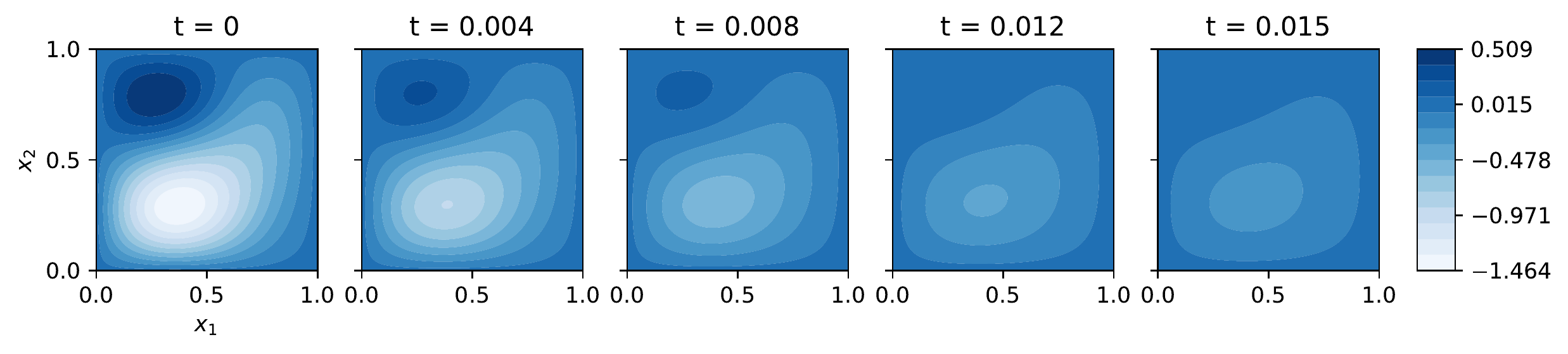}\vspace{-12pt}
\includegraphics[width=.75\textwidth]{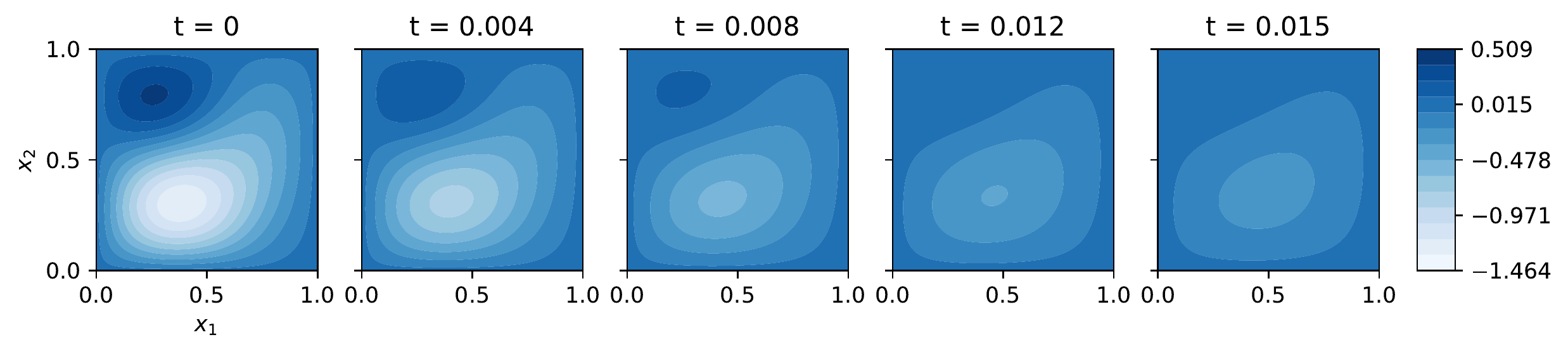}\vspace{-12pt}
\includegraphics[width=.75\textwidth]{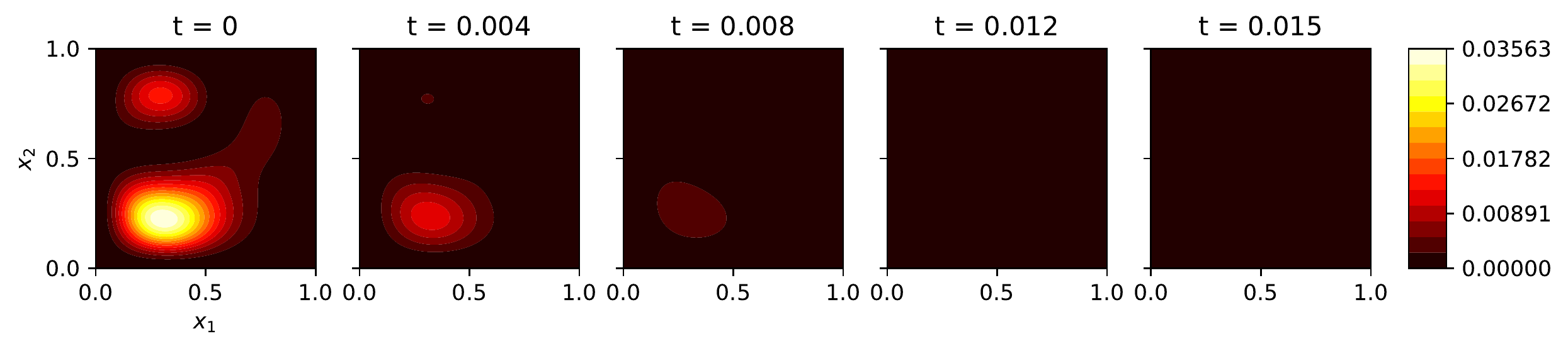}\vspace{-12pt}
\caption{(Heat equation). Comparison between true solution $u^{*}(\cdot,t)$, the approximation $\utt(\cdot)$ and their pointwise absolute difference $|\utt(x) - u^{*}(x,t)|$ for times $t=0,0.004,0.008,0.012,0.015$ for IVPs with the first (rows 1--3), second (rows 4--6) and third initial (rows 7--9) drawn from the set $\Gcal := \{\sum_{i=1}^4 c_i g_i \ : \ c_i \in [-1,1]\}$ where $g_i$ is defined in \eqref{eq:ivp-heat-initials}}
\label{fig:heat_example}
\end{figure}

% \begin{figure}
% \centering
% \includegraphics[scale=0.5]{pde_results/True_heat.pdf}
% \includegraphics[scale=0.5]{pde_results/NN_heat.pdf}
% \caption{Comparison between true and neural network evolution for $g(x)=\sin (2\pi x_1)\sin(2\pi x_2)\prod_{i=3}^{10}\sin(2\pi x_i)+0.5\Pi_{i=1}^10\sin(\pi x_i)$ at different time steps. This represents a perturbation contained within the region.}
% \label{fig:heat_error}
% \end{figure}

In order to have a class of analytical examples to compare against, we use the base functions
\begin{equation}
\begin{aligned}
g_1(x)&=\Pi_{i=1}^{10}\sin(\pi x_i),\\
g_2(x)&=\sin (2\pi x_1)\Pi_{i=1}^{10}\sin(\pi x_i),\\
g_3(x)&=\sin (2\pi x_2)\Pi_{i\neq 2}^{10}\sin(\pi x_i)\\
g_4(x)&=\sin (2\pi x_1)\sin(2\pi x_2)\Pi_{i=3}^{10}\sin(\pi x_i).
\end{aligned}
\label{eq:ivp-heat-initials}
\end{equation} 
to generate a class of initial conditions 
$\Gcal := \{\sum_{i=1}^4 c_i g_i \ : \ c_i \in [-1,1]\}$. 
To train our method, we drew 600 samples from $\Gcal$ and found a corresponding $\theta_0^{(j)}$ for each sample. We set the parameter space to be $\Theta:=\{\theta_0^{(j)} +\delta \ : \ |\delta|\leq 3, \ j=1,\ldots,600 \}$. We then uniformly sampled 200,000 points from this set $\Theta$ and generated paths for \eqref{eq:empirical_loss} from the 600 centers to train $V_{\xi}$. 
%\ye{uniformly sampled 200,000?}\gaby{Yes, clarified.} 
We then tested the method on a new set of 100 initials randomly drawn from $\Gcal$ by following the method outlined in Algorithm \ref{alg:solving-ivp}. 
We randomly select three from the test set containing the 100 initials and plot the result using our method in Figure \ref{fig:heat_example}. 
In addition, Figure \ref{fig:errors-c} and \ref{fig:errors-d} show the mean and standard deviations of the relative and absolute errors versus time $t$. 
We notice that the relative error increases while absolute error decreases: this is because the true solution $u^*(t,\cdot)$ gradually vanishes in time and hence it is easy to cause large relative error even when the absolute error is small.

In this test, it took 2.64 hours to generate $\tilde{G}$ and $\tilde{p}$ for \eqref{eq:empirical_loss} and 1.33 hours to generate the trajectories for \eqref{eq:data-driven}. This time cost is significantly higher than the transport equation as the heat equation requires the computation of the Laplacian which is second-order. Once the samples were generated, training $V_{\xi}$ took approximately 10 minutes. For testing, it took an average of 25 seconds to train a $\theta_0$ to a sampled $g$ and an average of 2.6 seconds to then solve \eqref{eq:ode} using a 4th order Runge-Kutta solver with step size 0.0001. This amounts to less than 30 seconds in time per initial for the testing stage.

\subsection{Allen-Cahn equation}
\begin{figure}
\centering
\includegraphics[width=.75\textwidth]{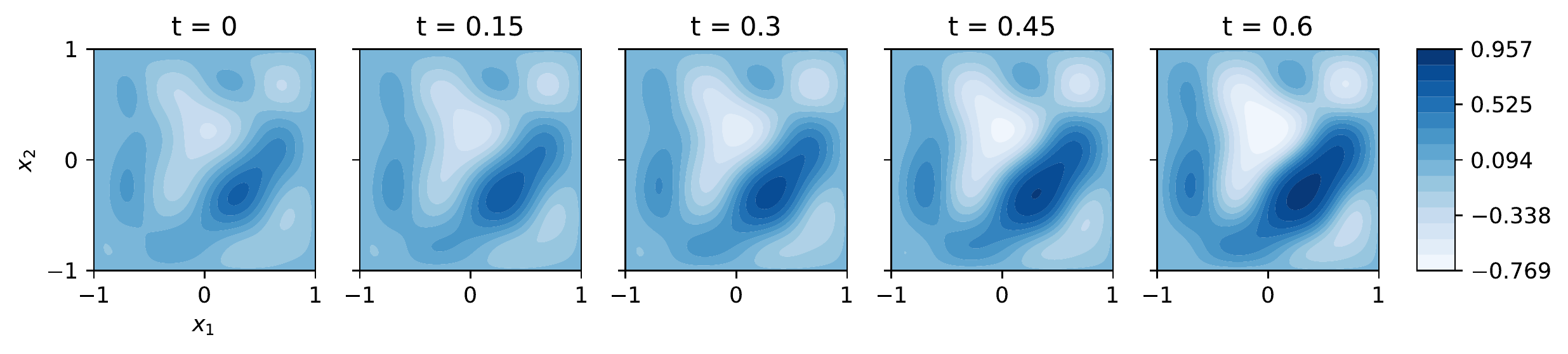}\vspace{-12pt}
\includegraphics[width=.75\textwidth]{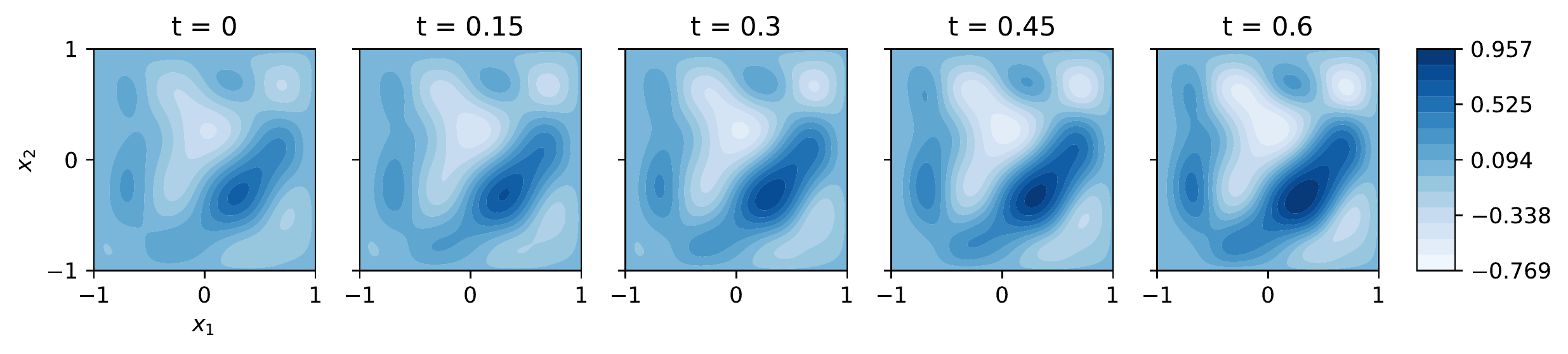}\vspace{-12pt}
\includegraphics[width=.75\textwidth]{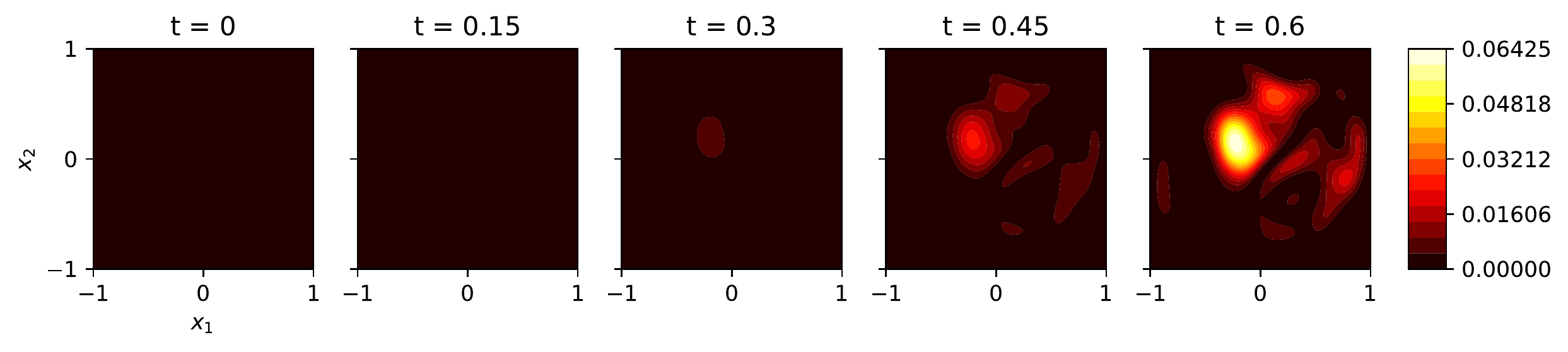}\vspace{-12pt}
\includegraphics[width=.75\textwidth]{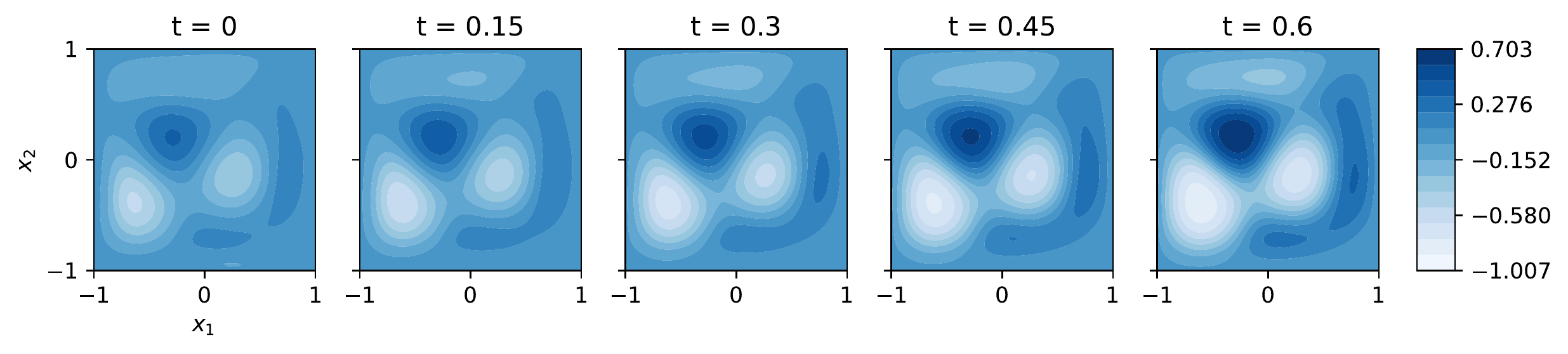}\vspace{-12pt}
\includegraphics[width=.75\textwidth]{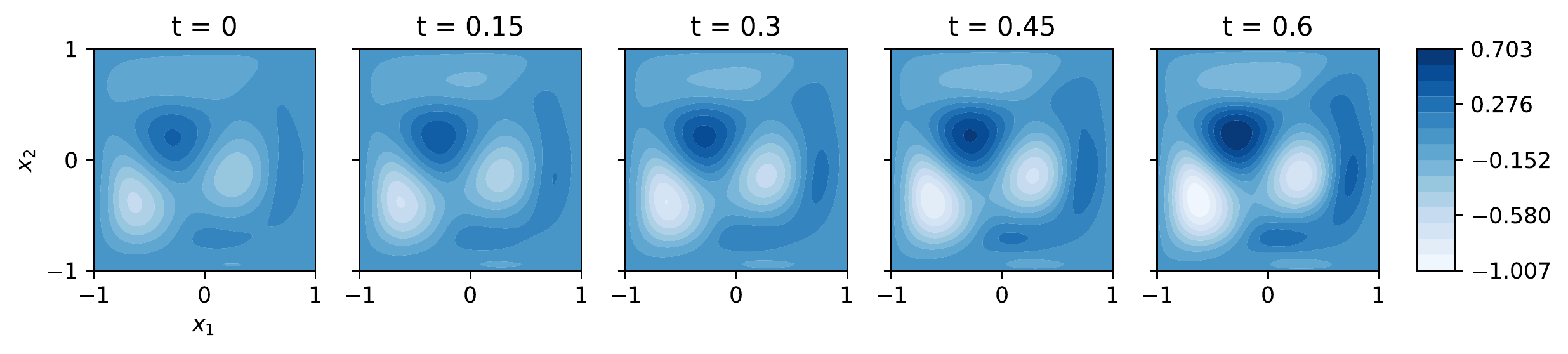}\vspace{-12pt}
\includegraphics[width=.75\textwidth]{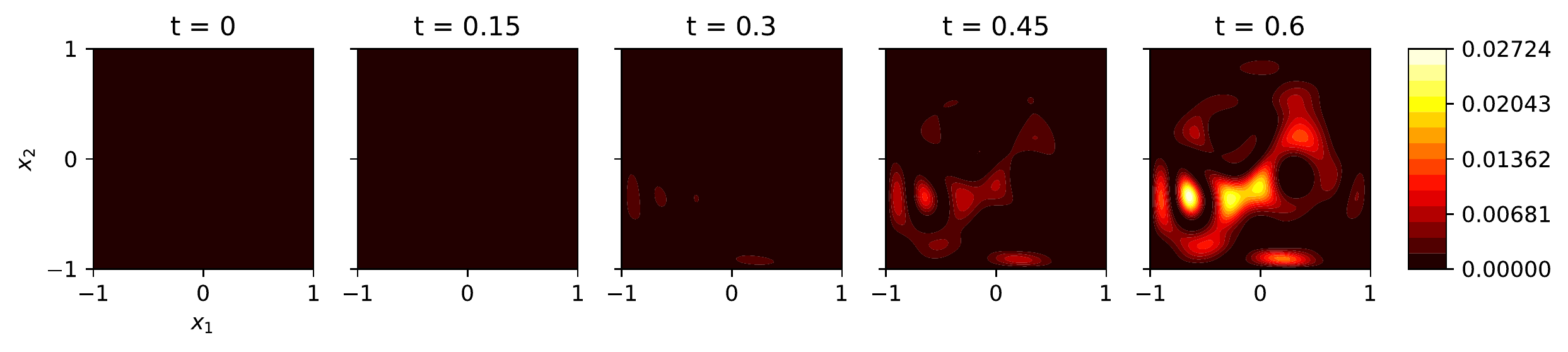}\vspace{-12pt}
\includegraphics[width=.75\textwidth]{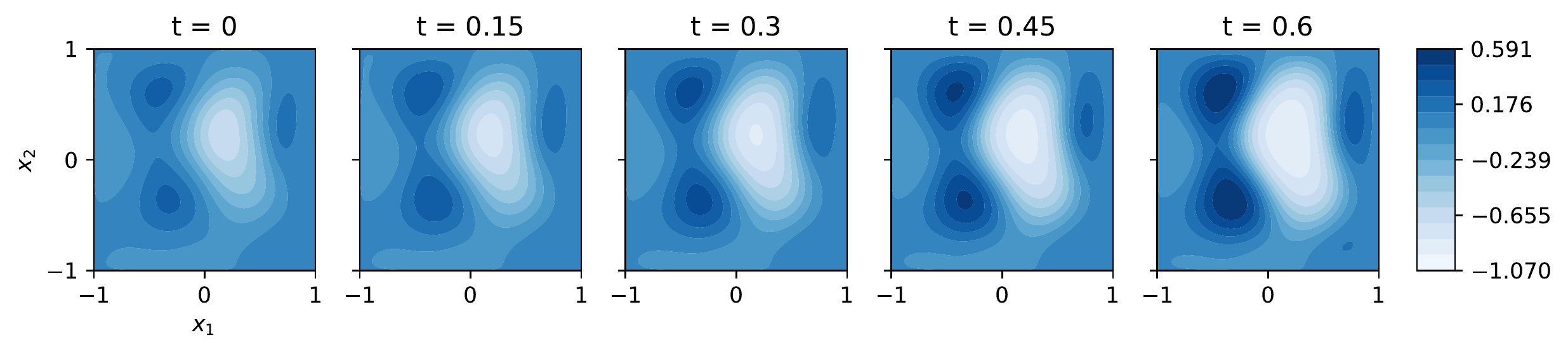}\vspace{-12pt}
\includegraphics[width=.75\textwidth]{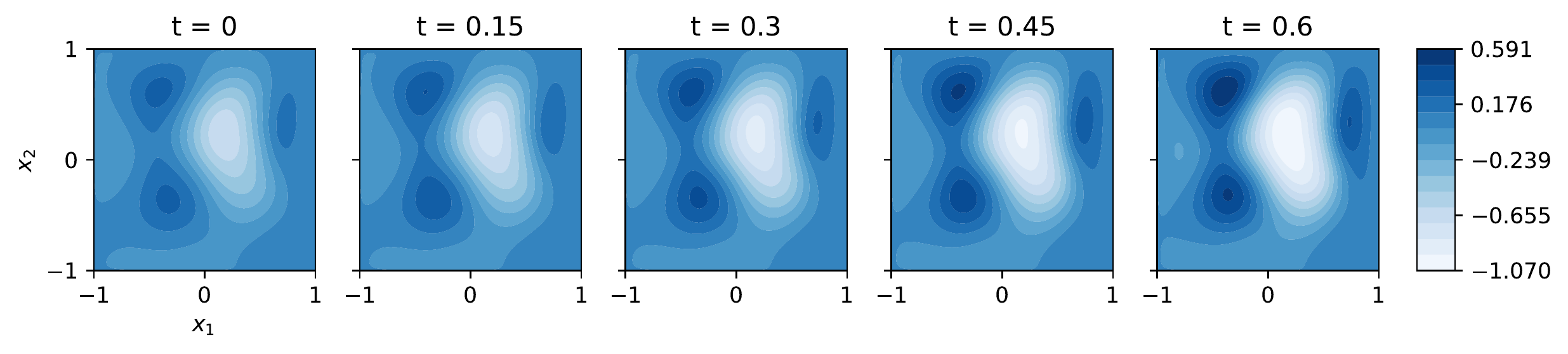}\vspace{-12pt}
\includegraphics[width=.75\textwidth]{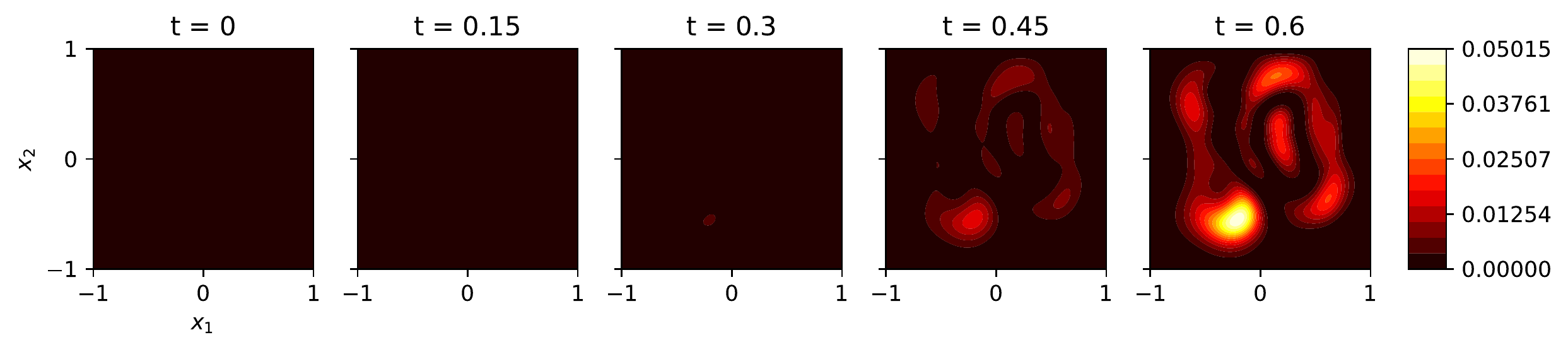}\vspace{-12pt}
\caption{(Allen-Cahn equation). Comparison between true solution $u^{*}(\cdot,t)$, the approximation $\utt(\cdot)$ and their pointwise absolute difference $|\utt(x) - u^{*}(x,t)|$ for times $t=0,0.004,0.008,0.012,0.015$ for IVPs with the first (rows 1--3), second (rows 4--6) and third initial (rows 7--9) drawn from the set $\Gcal$ defined in \eqref{eq:ivp-AC-initials}.}
\label{fig:AC_example}
\end{figure}

In this test, we consider the IVP with nonlinear Allen-Cahn equation given by
\begin{equation}
\begin{cases}
\partial_t u(x,t)=\epsilon \Delta u(x,t) + \frac{3}{2}\left(u(x,t)-u(x,t)^3\right), &\quad \forall\, x \in \Omega, t \in (0,T]\\
% u(x,t)&=0 & x \in \partial \Omega, t \in [0,T]\\
 u(x,0)  = g(x), & \quad \forall\, x \in \barOmega,
\end{cases}
\label{eq:ivp-AC-eq}
\end{equation}
where $\Omega =(-1,1)^2$, $\epsilon=0.0001$, and the boundary value $u(x,t)=0$ for all $x \in \partial \Omega$ and $t \in [0,T]$. As the Allen-Cahn PDE does not have an analytical solution to compare against, we resort to the classical implicit-explicit scheme (see e.g. \cite{tanh2016implicit}) with a 100$\times$100 grid and 2000 time points to generate a reference solution for comparison in 2D case only, despite that our method can be applied to higher dimensional case. In this test, we use \eqref{eq:nn-zero-bdry} with $\alpha(x)=(1-x_1^2)(1-x_2^2)$ as our neural network.
We let $T_i : \Rbb \to \Rbb$ represent the $i$th order Chebyshev polynomial. We generate a class of initial conditions 
\begin{equation}
\Gcal := \cbr[2]{(1-x_1^2)(1-x_2^2)\sum_{k=1}^{m} c_k T_{i_k}(x_1)T_{j_k}(x_2) \ : \ i_k, j_k \in \{0,\ldots,6\}, \ m \leq 36, \ |c_k|\leq 1 }.
\label{eq:ivp-AC-initials}
\end{equation} We see that $\Gcal$ is a space of all combinations of Chebyshev polynomials up to degree 6 multiplied by a boundary function. This set is chosen to represent a diverse spread of initials that can be approximated by our neural network from \eqref{eq:nn-zero-bdry}. We drew 200 samples from $\Gcal$ and found a corresponding $\theta_0^{(j)}$ for each sample. Then, as in the case of heat equations above, we generated the parameter set $\Theta:=\{\theta_0^{(j)} +\delta \ : \ |\delta|\leq 3, \ j=1,\ldots,200 \}$. We  sampled from $\Theta$ uniformly and generated paths from the 200 centers to train $V_{\xi}$. We again tested the method on a new set of 100 initials from $\Gcal$.
The results of the proposed method at time $t=0,0.15,0.3,0.45,0.6$ for three random initials are shown in Figure \ref{fig:AC_example}. 
% Again, the plots for $g_{1}$ are in rows 1--3, while the plots for $g_2$ and $g_3$ are shown in rows 4--6 and 7--9 in Figure \ref{fig:AC_example} respectively.
%
In Figure \ref{fig:errors-e} and \ref{fig:errors-f} we again plot the mean relative and absolute errors versus time, which demonstrate promising approximation performance of our method.

Figure \ref{fig:errors-e} shows some challenges in the relative error as time advances. This is because the solution to the Allen-Cahn equation for this initial value has fast-increasing derivatives as time progresses, which poses a challenge to all numerical methods including ours in solving Allen-Cahn equations in general. 
%This becomes a problem for the neural network model. 
Specifically, such large derivatives force the parameters $\theta$ of the neural network to blow up quickly, and hence the trajectory $\theta_t$ may rapidly escape from the prescribed $\Theta$ over which we trained the vector field $V_{\xi}$. This is a challenge that remains to be overcome by using more adaptive training methods and sampling strategies.

For this experiment, generating $\tilde{G}$ and $\tilde{p}$ for \eqref{eq:empirical_loss} took 1.04 hours while the generation of the trajectories for \eqref{eq:data-driven} took only 15 minutes. The much lower dimension of this problem compared to the transport and heat equation examples accounted for the speed up in the generation of samples. Similar to the transport equation, training $V_{\xi}$ took only 7 minutes. For testing, it took an average of 21 seconds to train a $\theta_0$ to a sampled $g$ and an average of 2.1 seconds to then solve \eqref{eq:ode} using a 4th order Runge-Kutta solver with step size 0.002. This amounts to less than 24 seconds in total time per initial for the testing stage.

\section{Variations and Generalizations}
\label{sec:variation}
In this section, we briefly discuss modifications of the proposed approach so that it can be applied to some other problems involving evolution PDEs. In particular, we consider the following two cases: general time-dependent PDEs and initial value problems with time-varying boundary conditions.

\paragraph{Applications to general time-dependent PDEs}
Our approach can be readily applied to a large variety of time-dependent PDEs.
The reason is that these PDEs can be converted to the exact form of \eqref{eq:pde} for which our method is designed.
To avoid overloading the bracket notation, we temporarily use $F(u)$ and $F(t,u)$ to represent $F[u]$ and $F_{t}[u]$ (differential operator that explicitly depends on time $t$).
We first note that one can convert any non-autonomous evolution PDE into an autonomous one:
%\ye{cite \url{https://www.math.ucdavis.edu/~hunter/notes/nonlinev.pdf}}
\begin{equation}
    \label{eq:nonauto-auto}
    \partial_t u = F(t,u) \quad \Longleftrightarrow \quad \partial_{t} \tilde{u} = \tilde{F}(\tilde{u}),\ \ \mbox{where} \ \ \tilde{u} := [t;u], \ \ \tilde{F}(\tilde{u}) := \big[1; F(u)\big],
\end{equation}
and $[\cdot\,;\, \cdot]$ means to stack the two arguments vertically to form a single one.
%
% Therefore, our approach can be easily modified to tackle non-autonomous PDEs due to \eqref{eq:nonauto-auto}.
%
We can also consider PDEs involving higher order time derivatives and convert them to first-order PDE systems by noticing equivalency as follows:
\begin{equation}
    \label{eq:wave}
    \partial_{tt} u = F(u) \quad \Longleftrightarrow \quad \partial_{t} \tilde{u} = \tilde{F}(\tilde{u}),\ \ \mbox{where} \ \ \tilde{u} := [u;v], \ \ \tilde{F}(\tilde{u}) := \big[v; F(u)\big].    
\end{equation}
History-dependent %(\zhou{delay?}) \ye{Works for general history-dependent case. Delayed PDE is a special case.} 
PDEs can also be considered: denote $H_{u}(t):=\{u(\cdot,s)\,|\, 0\le s \le t\}$ the trajectory recording the path of $u$ up to time $t$ and $F$ a nonlinear operator on path $H_{u}$, then we can set $H_{u}(t)$ as an auxiliary variable and convert the problem $\partial_t u = F[H_{u}]$ to an autonomous evolution PDE of $[u;H_u]$.

\paragraph{Evolution PDEs with boundary conditions}
We can also modify our method to solve IVPs with different boundary conditions. 
%We can extend our method to handle this situation. 
%
Let $(g,\phi)$ be the pair of initial and boundary values of the IVP. 
That is, $u(x,0)|_{\bar{\Omega}} = g$ and $u(x,t)|_{\partial \Omega \times [0,T]} = \phi$.
%
% Hence both of the initial value $g$ and the boundary value $\phi$ affect the solution.
%
In this case, we can parameterize $\ut(x) = \varphi_{\eta}(x) + \alpha(x) \psi_{\zeta}(x)$ with $\theta = (\eta, \zeta)$, where $\varphi_{\eta}$ and $\psi_{\zeta}$ are two reduced-order models (e.g., neural nets) with parameters $\eta$ and $\zeta$, respectively, and $\alpha(x)$ is a prescribed smooth function such that $\alpha(x) > 0$ if $x \in \Omega$ and $\alpha(x) = 0$ if $x \in \partial \Omega$.
Here $\varphi_{\eta}$ is to fit the boundary value $\phi$ without interference from $\alpha\psi_{\zeta}$ as the latter vanishes on the boundary $\partial \Omega$.

\section{Conclusion and Future Work}
\label{sec:conclusion}

We have shown a novel strategy for solving linear and nonlinear evolution PDEs numerically. Specifically, we propose to use deep neural networks as nonlinear reduced-order models to represent PDE solutions, and learn a control vector field to steer the network parameters so that the induced time-evolving neural network can approximate the solution accurately. The proposed method allows a user to quickly solve an evolution PDE with different initial values without the need to retrain the neural network. Error estimates of the proposed approach are also provided.
%The method is certified by using error estimates suggesting similar bounds as those in many ODE numerical methods. It is worth mentioning that as this method aims to solve a large amount of PDEs at once, we cannot hope for the approach to be as accurate as solvers aimed to solve a single PDE at a time. However, the flexibility of being able to compare many different solution trajectories without retraining is a significant potential benefit.

We implemented the nonlinear reduced-order models as generic deep networks which yield promising results. We expect that the accuracy and effectiveness can be further improved by incorporating structural information and prior knowledge about the PDE and its solutions into the design of these networks. Training of control vector fields can also be made more efficient by integrating informative sample trajectories of $\theta_t$. These improvements can potentially make the proposed method very effective in solving evolution PDEs in specified application domains.

\appendix

\section{Proof of \eqref{eq:psiK_bound}}
In the proof of Proposition \ref{prop:F_exists}, we need \eqref{eq:psiK_bound}. This can be obtained by applying the lemma below, whose proof is a slight modification of the proof of \cite[Theorem 2.1.14]{nesterov1998introductory}.
\label{appsec:proof}
\begin{lemma}
  Let $f:\Rbb^d \to \Rbb$ be a differentiable convex function and $\nabla f$ is $L$-Lipschitz continuous for some $L>0$. Define the gradient descent iterates by
  \[
    x_i = x_{i-1}-h\nabla f(x_{i-1})
  \] with $x_0 \in \Rbb^d$. Let $y \in \Rbb^d$ and $0<h<\frac{1}{L}$, then for any $k \ge 1$ there is
  \[
    f(x_k)-f(y) \leq \frac{|x_0-y|^2}{2kh}.
  \]
\label{lem:appendix}
\end{lemma}
\begin{proof}
    Following the standard steps in the proof of \cite[Theorem 2.1.14]{nesterov1998introductory} and using $0<h<1/L$ we can derive the bound
\begin{equation}
f(x_{i})-f(x_{i-1})\leq -h \del[2]{1-\frac{1}{2}hL }|\nabla f(x_{i-1})|^2 \leq -\frac{h}{2}|\nabla f(x_{i-1})|^2.
\label{eq:sup-eq1}
\end{equation} 
Since $f$ is convex, there is
\[
f(x) \leq f(y) + \nabla f(x)^{\top}(x-y), \quad \forall x \in \Rbb^d.
\]
Combining this with $x = x_{i-1}$ and \eqref{eq:sup-eq1}, we derive
\begin{align*}
    f(x_i)-f(y) & \leq \nabla f(x_{i-1})^{\top}(x_{i-1}-y)-\frac{h}{2}|\nabla f(x_{i-1})|^2\\
    &= \frac{1}{2h}\left(2h\nabla f(x_{i-1})^{\top}(x_{i-1}-y)-h^2|\nabla f(x_{i-1})|^2+|x_{i-1}-y|^2-|x_{i-1}-y|^2\right)\\
    & =\frac{1}{2h}\left(|x_{i-1}-y|^2-|x_{i-1}-h\nabla f(x_{i-1})-y|^2\right)\\
    &=\frac{1}{2h}\left(|x_{i-1}-y|^2-|x_{i}-y|^2\right).
\end{align*}
We can now bound the telescoping sum
\[
\sum_{i=1}^k (f(x_i)-f(y))\leq \frac{1}{2h}\sum_{i=1}^k(|x_{i-1}-y|^2-|x_{i}-y|^2)\leq\frac{1}{2h}|x_0-y|^2.
\] By \eqref{eq:sup-eq1} we know $f(x_k)\leq f(x_{k-1}) \leq \cdots \leq f(x_0)$ and therefore
\[
f(x_k)-f(y) \leq \frac{1}{k}\sum_{i=1}^k (f(x_i)-f(y)) \leq \frac{|x_0-y|^2}{2hk}.
\]
\end{proof}

\bibliographystyle{abbrv}
\bibliography{library}

\end{document}

%% file: ex_shared.tex
\usepackage{lipsum}
\usepackage{amsfonts}
\usepackage{graphicx}
\usepackage{epstopdf}
\usepackage{bm,commath,enumitem,color,hyperref}
\usepackage{mathfont}
\usepackage[sort]{cite}
\usepackage{graphicx,subcaption}

\ifpdf
  \DeclareGraphicsExtensions{.eps,.pdf,.png,.jpg}
\else
  \DeclareGraphicsExtensions{.eps}
\fi

\newsiamremark{remark}{Remark}
\newsiamremark{hypothesis}{Hypothesis}
\crefname{hypothesis}{Hypothesis}{Hypotheses}
\newsiamthm{claim}{Claim}
\newsiamthm{example}{Example}

\headers{Neural Control of Parametric Solutions for Evolution PDEs}{N. Gaby, X. Ye, and H. Zhou}

\title{Neural Control of Parametric Solutions for High-Dimensional Evolution PDEs\thanks{Submitted to the editors DATE.
\funding{This work was supported in part by National Science Foundation under grants DMS-1925263, DMS-2152960, DMS-2307465, DMS-2307466, and ONR N00014-21-1-2891.}}}

\author{Nathan Gaby\thanks{Department of Mathematics and Statistics, Georgia State University, Atlanta, GA, USA (\email{ngaby1@gsu.edu}).}
\and Xiaojing Ye\thanks{Department of Mathematics and Statistics, Georgia State University, Atlanta, GA, USA (\email{xye@gsu.edu}).}
\and Haomin Zhou\thanks{School of Mathematics, Georgia Institute of Technology, Atlanta, GA, USA (\email{hmzhou@gatech.edu}).}
}

\usepackage{amsopn}

%% file: main.bbl
\begin{thebibliography}{10}

\bibitem{ames2014numerical}
W.~F. Ames.
\newblock {\em Numerical methods for partial differential equations}.
\newblock Academic press, 2014.

\bibitem{anderson2022evolution}
W.~Anderson and M.~Farazmand.
\newblock Evolution of nonlinear reduced-order solutions for pdes with
  conserved quantities.
\newblock {\em SIAM Journal on Scientific Computing}, 44(1):A176--A197, 2022.

\bibitem{anitescu2019artificial}
C.~Anitescu, E.~Atroshchenko, N.~Alajlan, and T.~Rabczuk.
\newblock Artificial neural network methods for the solution of second order
  boundary value problems.
\newblock {\em Computers, Materials \& Continua}, 59(1):345--359, 2019.

\bibitem{atkinson1989introduction}
K.~Atkinson.
\newblock {\em An Introduction to Numerical Analysis (2nd ed.)}.
\newblock John Wiley I\& Sons, 1989.

\bibitem{bao2020numerical}
G.~Bao, X.~Ye, Y.~Zang, and H.~Zhou.
\newblock Numerical solution of inverse problems by weak adversarial networks.
\newblock {\em Inverse Problems}, 36(11):115003, 2020.

\bibitem{beck2017machine}
C.~Beck, W.~E, and A.~Jentzen.
\newblock Machine learning approximation algorithms for high-dimensional fully
  nonlinear partial differential equations and second-order backward stochastic
  differential equations.
\newblock {\em Journal of Nonlinear Science}, pages 1--57, 2017.

\bibitem{berg2018unified}
J.~Berg and K.~Nystr{\"o}m.
\newblock A unified deep artificial neural network approach to partial
  differential equations in complex geometries.
\newblock {\em Neurocomputing}, 317:28--41, 2018.

\bibitem{nicolas2022data-driven}
N.~Boull{\'e}, C.~Earls, and A.~Townsend.
\newblock Data-driven discovery of {Green's} functions with
  human-understandable deep learning.
\newblock {\em Scientific Reports}, 12:4824, 03 2022.

\bibitem{boulle2022learning-green}
N.~Boull{\'e}, S.~Kim, T.~Shi, and A.~Townsend.
\newblock Learning {G}reen's functions associated with time-dependent partial
  differential equations.
\newblock {\em Journal of Machine Learning Research}, 23:1--34, 08 2022.

\bibitem{bruna2022neural-galerkin}
J.~Bruna, B.~Pherstorfer, and E.~Vanden-Eijnden.
\newblock Neural {Galerkin} scheme with active learning for high-dimensional
  evolution equations.
\newblock {\em arXiv preprint arXiv:2203.01360}, 2022.

\bibitem{cai2020deepm}
S.~Cai, Z.~Wang, L.~Lu, T.~A. Zaki, and G.~E. Karniadakis.
\newblock Deepm\&mnet: Inferring the electroconvection multiphysics fields
  based on operator approximation by neural networks.
\newblock {\em arXiv preprint arXiv:2009.12935}, 2020.

\bibitem{cai2020phase}
W.~Cai, X.~Li, and L.~Liu.
\newblock A phase shift deep neural network for high frequency approximation
  and wave problems.
\newblock {\em SIAM Journal on Scientific Computing}, 42(5):A3285--A3312, 2020.

\bibitem{cai2019multi}
W.~Cai and Z.-Q.~J. Xu.
\newblock Multi-scale deep neural networks for solving high dimensional pdes.
\newblock {\em arXiv preprint arXiv:1910.11710}, 2019.

\bibitem{chen2020meta}
Y.~Chen, B.~Dong, and J.~Xu.
\newblock Meta-mgnet: Meta multigrid networks for solving parameterized partial
  differential equations.
\newblock {\em arXiv preprint arXiv:2010.14088}, 2020.

\bibitem{clark2020deep}
P.~Clark Di~Leoni, C.~Meneveau, G.~Karniadakis, and T.~Zaki.
\newblock Deep operator neural networks ({DeepONet}s) for prediction of
  instability waves in high-speed boundary layers.
\newblock {\em Bulletin of the American Physical Society}, 2020.

\bibitem{cuomo2022scientific}
S.~Cuomo, V.~S. Di~Cola, F.~Giampaolo, G.~Rozza, M.~Raissi, and F.~Piccialli.
\newblock Scientific machine learning through physics-informed neural networks:
  Where we are and what's next.
\newblock {\em arXiv preprint arXiv:2201.05624}, 2022.

\bibitem{dissanayake1994neural-network-based}
M.~Dissanayake and N.~Phan-Thien.
\newblock Neural-network-based approximations for solving partial differential
  equations.
\newblock {\em Communications in Numerical Methods in Engineering},
  10(3):195--201, 1994.

\bibitem{dong2020method}
S.~Dong and N.~Ni.
\newblock A method for representing periodic functions and enforcing exactly
  periodic boundary conditions with deep neural networks.
\newblock {\em arXiv preprint arXiv:2007.07442}, 2020.

\bibitem{du2021evolutional}
Y.~Du and T.~A. Zaki.
\newblock Evolutional deep neural network.
\newblock {\em Phys. Rev. E}, 104:045303, Oct 2021.

\bibitem{e2017deep}
W.~{E}, J.~{Han}, and A.~{Jentzen}.
\newblock Deep learning-based numerical methods for high-dimensional parabolic
  partial differential equations and backward stochastic differential
  equations.
\newblock {\em arXiv preprint arXiv:1706.04702}, 5(4):349--380, 2017.

\bibitem{e2018deep}
W.~E and B.~Yu.
\newblock The deep {Ritz} method: a deep learning-based numerical algorithm for
  solving variational problems.
\newblock {\em Communications in Mathematics and Statistics}, 6(1):1--12, 2018.

\bibitem{evans2012numerical}
G.~Evans, J.~Blackledge, and P.~Yardley.
\newblock {\em Numerical methods for partial differential equations}.
\newblock Springer Science \& Business Media, 2012.

\bibitem{evans1998partial}
L.~C. Evans.
\newblock {\em Partial differential equations}, volume~19 of {\em Graduate
  Studies in Mathematics}.
\newblock American Mathematical Society, Providence, RI, 1998.

\bibitem{fan2019bcr}
Y.~Fan, C.~O. Bohorquez, and L.~Ying.
\newblock Bcr-net: A neural network based on the nonstandard wavelet form.
\newblock {\em Journal of Computational Physics}, 384:1--15, 2019.

\bibitem{fujii2017asymptotic}
M.~Fujii, A.~Takahashi, and M.~Takahashi.
\newblock Asymptotic expansion as prior knowledge in deep learning method for
  high dimensional bsdes.
\newblock {\em Asia-Pacific Financial Markets}, pages 1--18, 2017.

\bibitem{gu2020structure}
Y.~Gu, C.~Wang, and H.~Yang.
\newblock Structure probing neural network deflation.
\newblock {\em arXiv preprint arXiv:2007.03609}, 2020.

\bibitem{gu2020selectnet}
Y.~Gu, H.~Yang, and C.~Zhou.
\newblock Selectnet: Self-paced learning for high-dimensional partial
  differential equations.
\newblock {\em arXiv preprint arXiv:2001.04860}, 2020.

\bibitem{guhring2020error}
I.~Guhring, G.~Kutyniok, and P.~Peterson.
\newblock Error bounds for approximations with deep relu neural networks in
  $w^{s,p}$ norms.
\newblock {\em Analysis and Applications}, 18:803--859, 2020.

\bibitem{guhring2020approximation}
I.~G{\"u}hring and M.~Raslan.
\newblock Approximation rates for neural networks with encodable weights in
  smoothness spaces.
\newblock {\em Neural Networks}, 134:107--130, 11 2020.

\bibitem{guo2018data-driven}
M.~Guo and J.~S. Hesthaven.
\newblock Data-driven reduced order modeling for time-dependent problems.
\newblock {\em Computer Methods in Applied Mechanics and Engineering},
  345:75--99, 2019.

\bibitem{guo2016cnn-operator}
X.~Guo, W.~Li, and F.~Iorio.
\newblock Convolutional neural networks for steady flow approximation.
\newblock In {\em Proceedings of the 22nd ACM SIGKDD International Conference
  on Knowledge Discovery and Data Mining}, page 481–490, New York, NY, USA,
  2016. Association for Computing Machinery.

\bibitem{han2017overcoming}
J.~Han, A.~Jentzen, and W.~E.
\newblock Overcoming the curse of dimensionality: Solving high-dimensional
  partial differential equations using deep learning.
\newblock {\em arXiv preprint arXiv:1707.02568}, pages 1--13, 2017.

\bibitem{han2018high-pde}
J.~Han, A.~Jentzen, and W.~E.
\newblock Solving high-dimensional partial differential equations using deep
  learning.
\newblock {\em Proceedings of the National Academy of Sciences},
  115(34):8505--8510, 2018.

\bibitem{han2020solving}
J.~Han, J.~Lu, and M.~Zhou.
\newblock Solving high-dimensional eigenvalue problems using deep neural
  networks: A diffusion monte carlo like approach.
\newblock {\em Journal of Computational Physics}, 423:109792, 2020.

\bibitem{han2018deep}
Y.~Han, J.~Yoo, H.~H. Kim, H.~J. Shin, K.~Sung, and J.~C. Ye.
\newblock Deep learning with domain adaptation for accelerated
  projection-reconstruction mr.
\newblock {\em Magnetic resonance in medicine}, 80(3):1189--1205, 2018.

\bibitem{hornik1991approximation}
K.~Hornik.
\newblock Approximation capabilities of multilayer feedforward networks.
\newblock {\em Neural networks}, 4(2):251--257, 1991.

\bibitem{hornik1989multilayer}
K.~Hornik, M.~Stinchcombe, and H.~White.
\newblock Multilayer feedforward networks are universal approximators.
\newblock {\em Neural networks}, 2(5):359--366, 1989.

\bibitem{huang2020int}
J.~Huang, H.~Wang, and H.~Yang.
\newblock Int-deep: A deep learning initialized iterative method for nonlinear
  problems.
\newblock {\em Journal of Computational Physics}, 419:109675, 2020.

\bibitem{hure2020deep}
C.~Hur{\'e}, H.~Pham, and X.~Warin.
\newblock Deep backward schemes for high-dimensional nonlinear pdes.
\newblock {\em Mathematics of Computation}, 89(324):1547--1579, 2020.

\bibitem{hutzenthaler2020proof}
M.~Hutzenthaler, A.~Jentzen, T.~Kruse, and T.~A. Nguyen.
\newblock A proof that rectified deep neural networks overcome the curse of
  dimensionality in the numerical approximation of semilinear heat equations.
\newblock {\em SN partial differential equations and applications}, 1(2):1--34,
  2020.

\bibitem{jagtap2020adaptive}
A.~D. Jagtap, K.~Kawaguchi, and G.~E. Karniadakis.
\newblock Adaptive activation functions accelerate convergence in deep and
  physics-informed neural networks.
\newblock {\em Journal of Computational Physics}, 404:109136, 2020.

\bibitem{johnson2012numerical}
C.~Johnson.
\newblock {\em Numerical solution of partial differential equations by the
  finite element method}.
\newblock Courier Corporation, 2012.

\bibitem{kharazmi2020hp}
E.~Kharazmi, Z.~Zhang, and G.~E. Karniadakis.
\newblock hp-vpinns: Variational physics-informed neural networks with domain
  decomposition.
\newblock {\em arXiv preprint arXiv:2003.05385}, 2020.

\bibitem{kovachki2021universal}
N.~Kovachki, S.~Lanthaler, and S.~Mishra.
\newblock On universal approximation and error bounds for {Fourier} neural
  operators.
\newblock {\em Journal of Machine Learning Research}, 22:Art--No, 2021.

\bibitem{kovachki2023neural}
N.~B. Kovachki, Z.~Li, B.~Liu, K.~Azizzadenesheli, K.~Bhattacharya, A.~M.
  Stuart, and A.~Anandkumar.
\newblock Neural operator: Learning maps between function spaces with
  applications to pdes.
\newblock {\em J. Mach. Learn. Res.}, 24(89):1--97, 2023.

\bibitem{kumar2011multilayer}
M.~Kumar and N.~Yadav.
\newblock Multilayer perceptrons and radial basis function neural network
  methods for the solution of differential equations: a survey.
\newblock {\em Computers \& Mathematics with Applications}, 62(10):3796--3811,
  2011.

\bibitem{lagaris1998artificial}
I.~E. Lagaris, A.~Likas, and D.~I. Fotiadis.
\newblock Artificial neural networks for solving ordinary and partial
  differential equations.
\newblock {\em IEEE transactions on neural networks}, 9(5):987--1000, 1998.

\bibitem{lee1990neural}
H.~Lee and I.~S. Kang.
\newblock Neural algorithm for solving differential equations.
\newblock {\em Journal of Computational Physics}, 91(1):110--131, 1990.

\bibitem{li2020repu}
B.~Li, S.~Tang, and H.~Yu.
\newblock Better approximations of high dimensional smooth functions by deep
  neural networks with rectified power units.
\newblock {\em Communications in Computational Physics}, 27:379--411, 02 2020.

\bibitem{li2020fourier}
Z.~Li, N.~Kovachki, K.~Azizzadenesheli, B.~Liu, K.~Bhattacharya, A.~Stuart, and
  A.~Anandkumar.
\newblock {Fourier} neural operator for parametric partial differential
  equations.
\newblock {\em arXiv preprint arXiv:2010.08895}, 2020.

\bibitem{li2020neural}
Z.~Li, N.~Kovachki, K.~Azizzadenesheli, B.~Liu, K.~Bhattacharya, A.~Stuart, and
  A.~Anandkumar.
\newblock Neural operator: Graph kernel network for partial differential
  equations.
\newblock {\em arXiv preprint arXiv:2003.03485}, 2020.

\bibitem{liang2017deep}
S.~Liang and R.~Srikant.
\newblock Why deep neural networks for function approximation?
\newblock In {\em International Conference on Learning Representations (ICLR)},
  2017.

\bibitem{liang2022finite}
S.~Liang and H.~Yang.
\newblock Finite expression method for solving high-dimensional partial
  differential equations.
\newblock {\em arXiv preprint arXiv:2206.10121}, 2022.

\bibitem{lin2022bi-green}
G.~Lin, F.~Chen, P.~Hu, X.~Chen, J.~Chen, J.~Wang, and Z.~Shi.
\newblock Bi-greennet: Learning {Green's} functions by boundary integral
  network.
\newblock {\em arXiv preprint arXiv:2204.13247}, 2022.

\bibitem{liu2020multi}
Z.~Liu, W.~Cai, and Z.-Q.~J. Xu.
\newblock Multi-scale deep neural network (mscalednn) for solving
  poisson-boltzmann equation in complex domains.
\newblock {\em arXiv preprint arXiv:2007.11207}, 2020.

\bibitem{lotzsch2022learning}
W.~L{\"o}tzsch, S.~Ohler, and J.~S. Otterbach.
\newblock Learning the solution operator of boundary value problems using graph
  neural networks.
\newblock {\em arXiv preprint arXiv:2206.14092}, 2022.

\bibitem{lu2019deeponet}
L.~Lu, P.~Jin, and G.~E. Karniadakis.
\newblock {DeepONet}: Learning nonlinear operators for identifying differential
  equations based on the universal approximation theorem of operators.
\newblock {\em arXiv preprint arXiv:1910.03193}, 2019.

\bibitem{lu2019deepxde}
L.~Lu, X.~Meng, Z.~Mao, and G.~E. Karniadakis.
\newblock Deepxde: A deep learning library for solving differential equations.
\newblock {\em arXiv preprint arXiv:1907.04502}, 2019.

\bibitem{luo2020two}
T.~Luo and H.~Yang.
\newblock Two-layer neural networks for partial differential equations:
  Optimization and generalization theory.
\newblock {\em arXiv preprint arXiv:2006.15733}, 2020.

\bibitem{lyu2020enforcing}
L.~Lyu, K.~Wu, R.~Du, and J.~Chen.
\newblock Enforcing exact boundary and initial conditions in the deep mixed
  residual method.
\newblock {\em arXiv preprint arXiv:2008.01491}, 2020.

\bibitem{magill2018neural}
M.~Magill, F.~Qureshi, and H.~de~Haan.
\newblock Neural networks trained to solve differential equations learn general
  representations.
\newblock In {\em Advances in Neural Information Processing Systems}, pages
  4071--4081, 2018.

\bibitem{mao2020deepm}
Z.~Mao, L.~Lu, O.~Marxen, T.~A. Zaki, and G.~E. Karniadakis.
\newblock {DeepM\&Mnet} for hypersonics: Predicting the coupled flow and
  finite-rate chemistry behind a normal shock using neural-network
  approximation of operators.
\newblock {\em arXiv preprint arXiv:2011.03349}, 2020.

\bibitem{nabian2018deep}
M.~A. Nabian and H.~Meidani.
\newblock A deep neural network surrogate for high-dimensional random partial
  differential equations.
\newblock {\em arXiv preprint arXiv:1806.02957}, 2018.

\bibitem{nesterov1998introductory}
Y.~Nesterov.
\newblock Introductory lectures on convex programming.
\newblock {\em Lecture Notes}, pages 119--120, 1998.

\bibitem{nusken2021solving}
N.~N{\"u}sken and L.~Richter.
\newblock Solving high-dimensional {Hamilton--Jacobi--Bellman} pdes using
  neural networks: perspectives from the theory of controlled diffusions and
  measures on path space.
\newblock {\em Partial Differential Equations and Applications}, 2(4):1--48,
  2021.

\bibitem{pang2020npinns}
G.~Pang, M.~D'Elia, M.~Parks, and G.~E. Karniadakis.
\newblock {nPINNs}: nonlocal physics-informed neural networks for a
  parametrized nonlocal universal laplacian operator. algorithms and
  applications.
\newblock {\em arXiv preprint arXiv:2004.04276}, 2020.

\bibitem{pang2019fpinns}
G.~Pang, L.~Lu, and G.~E. Karniadakis.
\newblock fpinns: Fractional physics-informed neural networks.
\newblock {\em SIAM Journal on Scientific Computing}, 41(4):A2603--A2626, 2019.

\bibitem{petersen2018optimal}
P.~Petersen and F.~Voigtlaender.
\newblock Optimal approximation of piecewise smooth functions using deep relu
  neural networks.
\newblock {\em Neural Networks}, 108:296--330, 2018.

\bibitem{pham2021neural}
H.~Pham, X.~Warin, and M.~Germain.
\newblock Neural networks-based backward scheme for fully nonlinear pdes.
\newblock {\em SN Partial Differ. Equ. Appl.}, 2(1), 2021.

\bibitem{quarteroni2008numerical}
A.~Quarteroni and A.~Valli.
\newblock {\em Numerical approximation of partial differential equations},
  volume~23.
\newblock Springer Science \& Business Media, 2008.

\bibitem{raissi2017machine}
M.~Raissi and G.~E. Karniadakis.
\newblock Machine learning of linear differential equations using gaussian
  processes.
\newblock {\em arXiv preprint arXiv:1701.02440}, 2017.

\bibitem{raissi2017physics}
M.~Raissi, P.~Perdikaris, and G.~E. Karniadakis.
\newblock Physics informed deep learning (part i): Data-driven solutions of
  nonlinear partial differential equations.
\newblock {\em arXiv preprint arXiv:1711.10561}, 2017.

\bibitem{raissi2019physics-informed}
M.~Raissi, P.~Perdikaris, and G.~E. Karniadakis.
\newblock Physics-informed neural networks: A deep learning framework for
  solving forward and inverse problems involving nonlinear partial differential
  equations.
\newblock {\em Journal of Computational Physics}, 378:686--707, 2019.

\bibitem{ramabathiran2021spinn}
A.~A. Ramabathiran and P.~Ramachandran.
\newblock Spinn: Sparse, physics-based, and partially interpretable neural
  networks for pdes.
\newblock {\em Journal of Computational Physics}, 445:110600, 2021.

\bibitem{raonic2023convolutional}
B.~Raoni{\'c}, R.~Molinaro, T.~Rohner, S.~Mishra, and E.~de~Bezenac.
\newblock Convolutional neural operators.
\newblock {\em arXiv preprint arXiv:2302.01178}, 2023.

\bibitem{regazzoni2019machine}
F.~Regazzoni, L.~Dedè, and A.~Quarteroni.
\newblock Machine learning for fast and reliable solution of time-dependent
  differential equations.
\newblock {\em Journal of Computational Physics}, 397:108852, 2019.

\bibitem{shin2020convergence}
Y.~Shin, J.~Darbon, and G.~E. Karniadakis.
\newblock On the convergence and generalization of physics informed neural
  networks.
\newblock {\em arXiv preprint arXiv:2004.01806}, 2020.

\bibitem{sirignano2018dgm:}
J.~Sirignano and K.~Spiliopoulos.
\newblock Dgm: A deep learning algorithm for solving partial differential
  equations.
\newblock {\em Journal of Computational Physics}, 375:1339--1364, 2018.

\bibitem{srivastava2015highway}
R.~K. Srivastava, K.~Greff, and J.~Schmidhuber.
\newblock Training very deep networks.
\newblock In C.~Cortes, N.~Lawrence, D.~Lee, M.~Sugiyama, and R.~Garnett,
  editors, {\em Advances in Neural Information Processing Systems}, volume~28.
  Curran Associates, Inc., 2015.

\bibitem{tabuada2022resnet}
P.~Tabuada and B.~Gharesifard.
\newblock Universal approximation power of deep residual neural networks
  through the lens of control.
\newblock {\em IEEE Transactions on Automatic Control}, pages 1--14, 2022.

\bibitem{tanh2016implicit}
T.~Tang and J.~Yang.
\newblock Implicit-explicit scheme for the allen-cahn equation preserves the
  maximum principle.
\newblock {\em Journal of Computational Mathematics}, 34:471--481, 09 2016.

\bibitem{teng2022green}
Y.~Teng, X.~Zhang, Z.~Wang, and L.~Ju.
\newblock Learning {Green's} functions of linear reaction-diffusion equations
  with application to fast numerical solver.
\newblock In {\em Proceedings of Mathematical and Scientific Machine Learning},
  volume 190 of {\em Proceedings of Machine Learning Research}, pages 1--16.
  PMLR, 15--17 Aug 2022.

\bibitem{thomas2013numerical}
J.~W. Thomas.
\newblock {\em Numerical partial differential equations: conservation laws and
  elliptic equations}, volume~33.
\newblock Springer Science \& Business Media, 2013.

\bibitem{thomas2013fdm}
J.~W. Thomas.
\newblock {\em Numerical partial differential equations: finite difference
  methods}, volume~22.
\newblock Springer Science \& Business Media, 2013.

\bibitem{wang2020multi}
B.~Wang, W.~Zhang, and W.~Cai.
\newblock Multi-scale deep neural network (mscalednn) methods for oscillatory
  stokes flows in complex domains.
\newblock {\em arXiv preprint arXiv:2009.12729}, 2020.

\bibitem{wang2022isl2}
C.~Wang, S.~Li, D.~He, and L.~Wang.
\newblock Is $ l^2$ physics-informed loss always suitable for training
  physics-informed neural network?
\newblock {\em arXiv preprint arXiv:2206.02016}, 2022.

\bibitem{wang2022respecting}
S.~Wang, S.~Sankaran, and P.~Perdikaris.
\newblock Respecting causality is all you need for training physics-informed
  neural networks.
\newblock {\em arXiv preprint arXiv:2203.07404}, 2022.

\bibitem{wang2021learning}
S.~Wang, H.~Wang, and P.~Perdikaris.
\newblock Learning the solution operator of parametric partial differential
  equations with physics-informed {DeepONets}.
\newblock {\em Science advances}, 7(40):eabi8605, 2021.

\bibitem{wen2022u-fno}
G.~Wen, Z.~Li, K.~Azizzadenesheli, A.~Anandkumar, and S.~M. Benson.
\newblock U-fno---an enhanced {Fourier} neural operator-based deep-learning
  model for multiphase flow.
\newblock {\em Advances in Water Resources}, 163:104180, 2022.

\bibitem{yang2020physics}
L.~Yang, D.~Zhang, and G.~E. Karniadakis.
\newblock Physics-informed generative adversarial networks for stochastic
  differential equations.
\newblock {\em SIAM Journal on Scientific Computing}, 42(1):A292--A317, 2020.

\bibitem{yang2019adversarial}
Y.~Yang and P.~Perdikaris.
\newblock Adversarial uncertainty quantification in physics-informed neural
  networks.
\newblock {\em Journal of Computational Physics}, 394:136--152, 2019.

\bibitem{yarotsky2017error}
D.~Yarotsky.
\newblock Error bounds for approximations with deep relu networks.
\newblock {\em Neural Networks}, 94:103--114, 2017.

\bibitem{zang2020weak}
Y.~Zang, G.~Bao, X.~Ye, and H.~Zhou.
\newblock Weak adversarial networks for high-dimensional partial differential
  equations.
\newblock {\em Journal of Computational Physics}, page 109409, 2020.

\bibitem{zhang2020physics}
E.~Zhang, M.~Yin, and G.~E. Karniadakis.
\newblock Physics-informed neural networks for nonhomogeneous material
  identification in elasticity imaging.
\newblock {\em arXiv preprint arXiv:2009.04525}, 2020.

\bibitem{zhu2018bayesian}
Y.~Zhu and N.~Zabaras.
\newblock Bayesian deep convolutional encoder–decoder networks for surrogate
  modeling and uncertainty quantification.
\newblock {\em Journal of Computational Physics}, 366:415--447, 2018.

\end{thebibliography}
